%% file: double-shuffle.tex
\documentclass[12pt,twoside,leqno,openany]{amsart}
\usepackage{amssymb,amsbsy,amsmath,amsfonts,amscd,times}
\usepackage{graphics,color,xy,footmisc,fancyhdr,multicol}
\usepackage{fancybox,graphicx,mathrsfs,shuffle,wasysym}
\usepackage[T1]{fontenc} 
\definecolor{blue}{cmyk}{1.,1.,0.,0.21}
\definecolor{red}{cmyk}{0.,1.,1.,0.21}
\definecolor{green}{cmyk}{1.,0.,1.,0.32}
\definecolor{black}{cmyk}{1.,1.,1.,1.}

\newcommand{\explain}[1]{\text{\scriptsize\sf [#1]}}
\let\mathcal\mathscr

\newcommand{\N}{\mathbb{N}}
\newcommand{\Q}{\mathbb{Q}}
\newcommand{\Z}{\mathbb{Z}}

\newcommand{\blue}{\textcolor{blue}}
\newcommand{\green}{\textcolor{green}}
\newcommand{\red}{\textcolor{red}}
\newcommand{\bbf}[1]{\blue{\mathbf{#1}}}
\newcommand{\gbf}[1]{\green{\mathbf{#1}}}

\newcommand{\zero}[1]{\underline{#1}_{\circ}}
\newcommand{\zerozero}[1]{\underline{#1}_{\circ\!\circ}}
\newcommand{\vf}{\vfill\end{document}}
\newcommand{\mathmotsf}[1]{\text{\footnotesize\sf #1}}
\newcommand{\smallmathmotsf}[1]{\text{\scriptsize\sf #1}}
\newcommand{\bfdots}{{}_{{}^{\bullet\,\bullet\,\bullet\,\bullet\,\bullet}}}
\newcommand{\zsqov}[1]{\overbrace{\bbf{0}\,\bbf{0}\,\cdots\,\bbf{0}\,\bbf{0}}^{#1}}
\newcommand{\unsqov}[1]{\bbf{1}\,\overbrace{\bbf{1}\,\bbf{1}\,\cdots\,\bbf{1}\,\bbf{1}}^{#1}}
\newcommand{\upbfdots}{{}^{{}_{{}^{\bullet\,\bullet\,\bullet\,\bullet\,\bullet}}}}

\newtheorem{The}{Theorem}[section]

\newtheorem{Theorem}{Theorem}[section]
\newtheorem{DualityTheorem}{Duality Theorem}[section]
\newtheorem{Proposition}[The]{Proposition}
\newtheorem{Lemma}[The]{Lemma}
\newtheorem{Corollary}[The]{Corollary}

\theoremstyle{definition}
\newtheorem{Definition}[The]{Definition}
\newtheorem{Notation}[The]{Notation}
\newtheorem{Summaryaboutnotations}[The]{Summary about notations}

\newcommand{\HEAD}[2]{%
\pagestyle{fancy}
\fancyhead[RO]{\scriptsize\sf\thepage}
\fancyhead[CO]{{\scriptsize\sf \thesection.\,\,\,#1}}
\fancyhead[LE]{\scriptsize\sf\thepage}
\fancyhead[CE]{{\scriptsize\sf #2}}
\fancyfoot{}}

\begin{document}

\title[]{
Multizeta Calculus (I)
}

\address{D\'epartment de Math\'ematiques d'Orsay,
B\^atiment 425, Facult\'e des Sciences, 
Universit\'e Paris XI - Orsay, 
F-91405 Orsay Cedex, FRANCE}
\email{Joel.Merker@math.u-psud.fr}

\date{\number\year-\number\month-\number\day}

\maketitle

\begin{center}
Jo\"el {\sc Merker}
\end{center}

\bigskip

\begin{center}
\begin{minipage}[t]{11.75cm}
\baselineskip =0.35cm {\scriptsize

\centerline{\bf Table of contents}

\smallskip

{\bf \ref{Introduction}. Introduction
\dotfill~\pageref{Introduction}.}

{\bf \ref{weight-depth-height}. Weight, Depth, Height
\dotfill~\pageref{weight-depth-height}.}

{\bf \ref{dualities}. Dualities
\dotfill~\pageref{dualities}.}

{\bf \ref{detachable-notations}. Detachable and Flexible Notations
\dotfill~\pageref{detachable-notations}.}

{\bf \ref{total-ordering-polyzetas}. Total Ordering $\prec$ 
Between Polyzetas of Fixed Weight
\dotfill~\pageref{total-ordering-polyzetas}.}

{\bf \ref{intermediate-countings}. Intermediate Countings
\dotfill~\pageref{intermediate-countings}.}

{\bf \ref{shuffle-stuffle}. Shuffle Minus Stuffle
\dotfill~\pageref{shuffle-stuffle}.}

{\bf \ref{1-double-melange}. 
Computing $-(1)\ast\big(a_1,1^{b_1}, \dots,a_{\sf h},1^{b_{\sf h}}\big)
+ 1 \shuffle 0^{ a_1-1} 1\, 1^{b_1} \cdots 0^{ a_{\sf h} -1}
1\, 1^{b_{\sf h}}$
\dotfill~\pageref{1-double-melange}.}

{\bf \ref{2-double-melange}. 
Computing $-(2)\ast\big(a_1,1^{b_1}, \dots,a_{\sf h},1^{b_{\sf h}}\big)
+ 01 \shuffle 0^{ a_1-1} 1\, 1^{b_1} \cdots 0^{ a_{\sf h} -1}
1\, 1^{b_{\sf h}}$
\dotfill~\pageref{2-double-melange}.}

{\bf \ref{3-double-melange}. 
Computing $-(3)\ast\big(a_1,1^{b_1}, \dots,a_{\sf h},1^{b_{\sf h}}\big)
+ 001 \shuffle 0^{ a_1-1} 1\, 1^{b_1} \cdots 0^{ a_{\sf h} -1}
1\, 1^{b_{\sf h}}$
\dotfill~\pageref{3-double-melange}.}

{\bf \ref{21-double-melange}. 
Computing $-(2,1)\ast\big(a_1,1^{b_1}, \dots,a_{\sf h},1^{b_{\sf h}}\big)
+ 011 \shuffle 0^{ a_1-1} 1\, 1^{b_1} \cdots 0^{ a_{\sf h} -1}
1\, 1^{b_{\sf h}}$
\dotfill~\pageref{21-double-melange}.}

}\end{minipage}
\end{center}

\bigskip

\hfill
\begin{minipage}[t]{11.25cm}
{\footnotesize\sf\em
The Holy Grail in the field of MZV's {\em [Multiple Zeta Values]} is
an algorithm to express each MZV into an unique basis in a
constructive way. That way we would have a (hopefully) small procedure
rather than giant tables. Thus far this has not been found.}\hfill
J.A.M. {\sc Vermaseren}.
\end{minipage}

\section{Introduction}
\label{Introduction}
\HEAD{Introduction}{
Jo\"el {\sc Merker}, Universit\'e Paris-Sud Orsay, France}

\bigskip

As is known (\cite{ 
Zagier-1994, Hoffman-1997, JMOP-1999, Cartier-2001}, convergent polyzetas
are infinite numerical sums:
\[
\zeta\big(s_1,s_2,\dots,s_{\sf d}\big)
=
\sum_{n_1>n_2>\cdots>n_{\sf d}>0}\,
\frac{1}{(n_1)^{s_1}\,(n_2)^{s_2}\,\cdots\,
(n_{\sf d})^{s_{\sf d}}},
\]
for certain integers $s_1, s_2, \dots, s_{\sf d} \geqslant 1$, all
positive, with the additional requirement that $s_1 \geqslant 2$ in
order to insure (absolute) convergence of this multiple series in
which $n_1 > n_2 > \cdots > n_{\sf d} > 0$ are integers.  The number
${\sf d}$ of the appearing 
entries $s_i$ is called the {\sl depth} of the polyzeta, while
the sum of entries:
\[
s_1+s_2+\cdots+s_{\sf d}
=:
{\sf w}
\]
is said to be its {\sl weight}. 

It is expected\,\,---\,\,or conjectured\,\,---\,\,that all the 
so-called {\sl double shuffle relations}:
\[
0
=
-\,
\zeta\big(s_1,s_2,\dots,s_{\sf d}\big)
\ast
\zeta\big(t_1,t_2,\dots,t_{\sf e}\big)
+
\zeta\big(s_1,s_2,\dots,s_{\sf d}\big)
\shuffle
\zeta\big(t_1,t_2,\dots,t_{\sf e}\big),
\]
including the ones for which $s_1 = 1$\,\,---\,\,a single
non-convergent polyzeta in the two products is then erased by the
subtraction in this case\,\,---\,\,provide all $\Q$-linear relations
between the convergent polyzetas. Since double shuffle relations are
{\em $\Q$-linear} between polyzetas of any fixed weight, such a
conjecture amounts, metaphysically speaking, to believe that
easy-to-discover structures exhaust {\em all} 
polyzeta relations.

Intentionally, we will not produce here any `sexy' basic
presentation of conjectures and structures
(fundamental references are listed in the bibliography),
because our goal is instead to set up `handy' formulas 
towards a (new?) {\sl Polyzeta Calculus}.

A first immediate aspect of the problems that are open in the field
would be to understand {\em in a closed way} all what is contained in
these double shuffle relations.

A second, much more delicate, aspect would be to hope for
irrationality and even for transcendence concerning polyzetas {\em
modulo} these double shuffle relations, but this question will not
at all be dealt with here, because, due to its complexity, a
satisfactory settlement of it could well require yet several decades
of intensive mathematical research (\cite{ Cartier-2011}).

Concerning the first aspect, the so-called {\sl Hoffman conjecture}
expects that the polyzetas whose entries $s_i$ are equal all to either
$2$ or $3$ make up a basis of the vector space of all polyzetas of any
fixed weight ${\sf w}$ modulo all existing double shuffle
($\Q$-linear) relations within a same fixed weight.

The correctness of this quite appealing conjecture has been verified
on powerful computer machines up to weight ${\sf w} = 22$
(\cite{MMK-2008, BBV-2010}) and a
striking proof that every polyzeta is a certain $\Q$-linear
combination of the $\zeta ( s_1, \dots, s_{\sf d})$ with the $s_i 
\in \{ 2, 3\}$
was obtained recently by Francis Brown (\cite{
Brown-2011}).

Still, experts agree that a (wide?) gap exists between, on one hand,
abstract, indirect or non computationally complete arguments, and, on
the other hand, what our sparkling computer machines really handle
within their most intimate core.

Therefore, our present goal here will be to offer a contribution that,
hopefully, might be new, in the sense that it would launch the
construction of some archway for a bridge between brain and machine
towards satisfactory constructiveness in the $(2, 3)$-conjecture.  A
strong result of Jean \'Ecalle (\cite{ Ecalle-2009}) already showed in
an effective way that all entries $s_i = 1$ may be eliminated.

In 2008, Masanobu Kaneko, Masayuki Noro and Ken'ichi Tsurumaki (\cite{
MMK-2008}) showed that up to weight ${\sf w} = 20$, it suffices in
fact to look only at all double shuffle relations whose first polyzeta
$\zeta( s_1, \dots, s_{\sf d})$ is either $\zeta ( 1)$, $\zeta ( 2)$,
$\zeta ( 3)$ or $\zeta ( 2, 1)$\,\,---\,\,the consideration of duality
relations is also advantageous.  This improves the $(2, 3)$-conjecture
in a highly natural manner, because the number of relations obtained
in such a way for polyzetas of weight ${\sf w}$:
\[
2^{{\sf w}-3}
+
2^{{\sf w}-4}
+
2^{{\sf w}-5}
+
2^{{\sf w}-5}
\]
is {\em exactly equal} to the number:
\[
2^{{\sf w}-2}
\]
of polyzetas of weight ${\sf w}$, hence its exceeds in the right way
the expected number $2^{{\sf w}-2} - \delta_{\sf w}$ of {\em
dependent} polyzetas, where $\delta_2 = \delta_3 = \delta_4 = 1$,
$\delta_{\sf w} = \delta_{{\sf w} - 2} + \delta_{{\sf w} - 3}$
is the number of polyzetas of weight ${\sf w}$ having
entries belonging to the set $\{ 2, 3\}$.

In this paper, we develop
techniques and we establish theorems in order to explicitly write down
these four families of {\em conjecturably relevant} double shuffle relations:
\[
\label{1-2-3-21}
\small
\aligned
0
&
=
-\,
\zeta(1)
\ast
\zeta\big(a_1,1^{b_1},\dots,a_{\sf h},1^{b_{\sf h}}\big)
+
\zeta(1)
\shuffle
\zeta\big(a_1,1^{b_1},\dots,a_{\sf h},1^{b_{\sf h}}\big)
\ \ \ \ \ \ \ \ \ \ \ \ \
\explain{{\sf w}-1},
\\
0
&
=
-\,
\zeta(2)
\ast
\zeta\big(a_1,1^{b_1},\dots,a_{\sf h},1^{b_{\sf h}}\big)
+
\zeta(2)
\shuffle
\zeta\big(a_1,1^{b_1},\dots,a_{\sf h},1^{b_{\sf h}}\big)
\ \ \ \ \ \ \ \ \ \ \ \ \
\explain{{\sf w}-2},
\\
0
&
=
-\,
\zeta(3)
\ast
\zeta\big(a_1,1^{b_1},\dots,a_{\sf h},1^{b_{\sf h}}\big)
+
\zeta(3)
\shuffle
\zeta\big(a_1,1^{b_1},\dots,a_{\sf h},1^{b_{\sf h}}\big)
\ \ \ \ \ \ \ \ \ \ \ \ \
\explain{{\sf w}-3},
\\
0
&
=
-\,
\zeta(2,1)
\ast
\zeta\big(a_1,1^{b_1},\dots,a_{\sf h},1^{b_{\sf h}}\big)
+
\zeta(2,1)
\shuffle
\zeta\big(a_1,1^{b_1},\dots,a_{\sf h},1^{b_{\sf h}}\big)
\ \ \ \ \ \
\explain{{\sf w}-3}.
\endaligned
\]

By definition, the {\sl height} ${\sf h}$ of a polyzeta is the number
of its entries $s_i$ that are $\geqslant 2$.  It follows that a
general, arbitrary convergent polyzeta can always be written under the
specific form:
\[
\zeta\big(a_1,1^{b_1},a_2,1^{b_2},\dots,
a_{\sf h},1^{b_{\sf h}}\big),
\]
where the integer ${\sf h} \geqslant 1$ is its
height, where all the $a$-integers are $\geqslant 2$:
\[
a_1\geqslant 2,
\ \ \ \ \ \ \ \ 
a_2\geqslant 2,
\ \ \ \ \ \ \ \
\dots\dots,
\ \ \ \ \ \ \ \
a_{\sf h}\geqslant 2,
\]
and where all the $b$-integers are $\geqslant 0$:
\[
b_1\geqslant 0,
\ \ \ \ \ \ \ \
b_2\geqslant 0,
\ \ \ \ \ \ \ \
\dots\dots,
\ \ \ \ \ \ \ \
b_{\sf h}\geqslant 0.
\]

\begin{Theorem}
For every height ${\sf h} \geqslant 1$,
every entries $a_1 \geqslant 2$, \dots, $a_{\sf h} \geqslant 2$,
every entries $b_1 \geqslant 0$, \dots, $b_{\sf h} \geqslant 0$,
if one sets:
\[
\aligned
a_1+b_1+\cdots+a_{\sf h}+b_{\sf h}
&
=:
{\sf w}-1,
\\
1+b_1+\cdots+1+b_{\sf h}
&
=:
{\sf d},
\endaligned
\]
then the so-called {\em regularized} double shuffle relation: 
\[
0
=
-\,(1)
\ast
\big(a_1,1^{b_1},\dots,a_{\sf h},1^{b_{\sf h}}\big)
+
1
\shuffle
0^{a_1-1}1\,1^{b_1}\,\cdots\,0^{a_{\sf h}-1}1\,1^{b_{\sf h}}
\]
between $\zeta ( 1)$ and any $\zeta_{[{\sf w}-1,{\sf d},{\sf h}]}
\big( a_1, 1^{ b_1}, \dots, a_{\sf h}, 1^{ b_{\sf h}} \big)$ of weight
${\sf w} - 1 \geqslant 2$, of depth ${\sf d}$ and of height ${\sf h}$ writes out
after complete finalization:
\[
\boxed{
\aligned
0
&
=
-\,\sum_{1\leqslant i\leqslant{\sf h}}\,
\zeta_{[{\sf w},{\sf d},{\sf h}]}
\big(\bfdots,a_i+1,\bfdots\big)
-
\\
&
\ \ \ \ \
-
\sum_{1\leqslant j\leqslant{\sf h}
\atop b_j\geqslant 1}
\Bigg(
\sum_{b_j'+b_j''=b_j-1}\,
\zeta_{[{\sf w},{\sf d},{\sf h}+1]}
\big(\bfdots,1^{b_j'},2,1^{b_j''},\bfdots\big)
\Bigg)
+
\\
&
\ \ \ \ \
+
\sum_{1\leqslant j\leqslant{\sf h}}\,
\zeta_{[{\sf w},{\sf d}+1,{\sf h}]}
\big(\bfdots,1^{b_j+1},\bfdots\big)
+
\\
&
\ \ \ \ \
+
\sum_{1\leqslant i\leqslant{\sf h}
\atop a_i\geqslant 3}
\Bigg(
\sum_{a_i'+a_i''=a_i+1
\atop
a_i'\geqslant 2,\,\,a_i''\geqslant 2}\,
\zeta_{[{\sf w},{\sf d}+1,{\sf h}+1]}
\big(\bfdots,a_i',a_i'',\bfdots\big)
\Bigg).
\endaligned}
\]
\end{Theorem}

Here by convention, only the terms of the
initial polyzeta $\zeta_{[ {\sf w}-1, {\sf d}, {\sf h}]}$ that are
{\em changed} are written down, so that the symbol:
\[
\bfdots
\]
means that all other entries are unchanged, as one
might have guessed at once.

\begin{Theorem}
For every height ${\sf h} \geqslant 1$,
every entries $a_1 \geqslant 2$, \dots, $a_{\sf h} \geqslant 2$,
every entries $b_1 \geqslant 0$, \dots, $b_{\sf h} \geqslant 0$,
if one sets:
\[
\aligned
a_1+b_1+\cdots+a_{\sf h}+b_{\sf h}
&
=:
{\sf w}-2,
\\
1+b_1+\cdots+1+b_{\sf h}
&
=:
{\sf d},
\endaligned
\]
then the so-called double shuffle relation: 
\[
0
=
-\,(2)
\ast
\big(a_1,1^{b_1},\dots,a_{\sf h},1^{b_{\sf h}}\big)
+
01
\shuffle
0^{a_1-1}1\,1^{b_1}\,\cdots\,0^{a_{\sf h}-1}1\,1^{b_{\sf h}}
\]
between $\zeta (2)$ and any $\zeta_{[{\sf w}-2,{\sf d},{\sf h}]}
\big( a_1, 1^{ b_1}, \dots, a_{\sf h}, 1^{ b_{\sf h}} \big)$ of weight
${\sf w} - 2 \geqslant 2$, of depth ${\sf d}$ 
and of height ${\sf h}$ writes out after complete
finalization:
\[
\boxed{
\small
\aligned
0
&
=
-\,\sum_{1\leqslant i\leqslant{\sf h}}\,
\zeta_{[{\sf w},{\sf d},{\sf h}]}
\big(\bfdots,a_i+2,\bfdots\big)
-
\\
&
\ \ \ \ \
-\,
\sum_{1\leqslant j\leqslant{\sf h}\atop b_j\geqslant 1}
\Bigg(
\sum_{b_j'+b_j''=b_j-1
\atop
b_j'\geqslant 0,\,\,b_j''\geqslant 0}\,
\zeta_{[{\sf w},{\sf d},{\sf h}+1]}
\big(\bfdots,1^{b_j'},3,1^{b_j''},\bfdots\big)
\Bigg)
+
\\
&
\ \ \ \ \
+
\sum_{1\leqslant i\leqslant j\leqslant{\sf h}}\,
a_i\,(b_j+2)\cdot
\zeta_{[{\sf w},{\sf d}+1,{\sf h}]}
\big(\bfdots,a_i+1,\bfdots,1^{b_j+1},\bfdots\big)
+
\\
&
\ \ \ \ \
+
\sum_{1\leqslant j_1<j_2\leqslant{\sf h}\atop b_{j_1}\geqslant 1}
\Bigg(
\sum_{b_{j_1}'+b_{j_1}''=b_{j_1}-1
\atop
b_{j_1}'\geqslant 0,\,\,b_{j_1}''\geqslant 0}\,
(b_{j_2}+2)\cdot
\zeta_{[{\sf w},{\sf d}+1,{\sf h}+1]}
\big(\bfdots,1^{b_{j_1}'},2,1^{b_{j_1}''},\bfdots,1^{b_{j_2}+1},\bfdots\big)
\Bigg)
+
\\
&
\ \ \ \ \
+
\sum_{1\leqslant j\leqslant{\sf h}}
\Bigg(
\sum_{b_j'+b_j''=b_j
\atop
b_j'\geqslant 0,\,\,b_j''\geqslant 0}\,
b_j''\cdot
\zeta_{[{\sf w},{\sf d}+1,{\sf h}+1]}
\big(\bfdots,1^{b_j'},2,1^{b_j''},\bfdots\big)
\Bigg)
+
\\
&
\ \ \ \ \
+
\sum_{1\leqslant i\leqslant{\sf h}\atop a_i\geqslant 5}
\Bigg(
\sum_{a_i'+a_i''=a_i+2
\atop
a_i'\geqslant 3,\,\,a_i''\geqslant 2}\,
(a_i'-1)\cdot
\zeta_{[{\sf w},{\sf d}+1,{\sf h}+1]}
\big(\bfdots,a_i',1^0,a_i'',\bfdots\big)
\Bigg)
+
\\
&
\ \ \ \ \
+
\sum_{1\leqslant i_1<i_2\leqslant{\sf h}\atop a_{i_2}\geqslant 3}\,
\Bigg(
\sum_{a_{i_2}'+a_{i_2}''=a_{i_2}+1
\atop
a_{i_2}'\geqslant 2,\,\,a_{i_2''}\geqslant 2}\,
a_{i_1}
\cdot
\zeta_{[{\sf w},{\sf d}+1,{\sf h}+1]}
\big(\bfdots,a_{i_1}+1,\bfdots,a_{i_2}',1^0,a_{i_2}'',\bfdots\big)
\Bigg)
+
\\
&
\ \ \ \ \
+
\sum_{1\leqslant j<i\leqslant{\sf h}
\atop
b_j\geqslant 1,\,\,a_i\geqslant 3}
\Bigg(
\sum_{b_j'+b_j''=b_j-1
\atop
b_j'\geqslant 0,\,\,b_j''\geqslant 0}\,
\sum_{a_i'+a_i''=a_i+1
\atop
a_i'\geqslant 2,\,\,a_i''\geqslant 2}\,
\zeta_{[{\sf w},{\sf d}+1,{\sf h}+2]}
\big(\bfdots,1^{b_j'},2,1^{b_j''},\bfdots,a_i',1^0,a_i'',\bfdots\big)
\Bigg).
\endaligned}
\]
\end{Theorem}

\begin{Theorem}
For every height ${\sf h} \geqslant 1$,
every entries $a_1 \geqslant 2$, \dots, $a_{\sf h} \geqslant 2$,
every entries $b_1 \geqslant 0$, \dots, $b_{\sf h} \geqslant 0$,
if one sets:
\[
\aligned
a_1+b_1+\cdots+a_{\sf h}+b_{\sf h}
&
=:
{\sf w}-3,
\\
1+b_1+\cdots+1+b_{\sf h}
&
=:
{\sf d},
\endaligned
\]
then the so-called $\ast$-stuffle product between $\zeta ( 3)$ and
any $\zeta_{[{\sf w}-3, {\sf d}, {\sf h}]} \big( a_1, 1^{ b_1}, \dots,
a_{\sf h}, 1^{ b_{\sf h}} \big)$ of weight ${\sf w} - 3
\geqslant 2$, of depth
${\sf d}$ and of height ${\sf h}$ writes out after complete
finalization:
\[
\aligned
\zeta(3)
&
\ast
\zeta_{[{\sf w}-3,{\sf d},{\sf h}]}
\big(a_1,1^{b_1},\dots,a_{\sf h},1^{b_{\sf h}}\big)
=
\\
&
=
\sum_{1\leqslant i\leqslant{\sf h}}\,
\zeta_{[{\sf w},{\sf d},{\sf h}]}
\big(\bfdots,a_i+3,\bfdots\big)
+
\\
&
\ \ \ \ \
+
\sum_{1\leqslant j\leqslant{\sf h}
\atop b_j\geqslant 1}
\Bigg(
\sum_{b_j'+b_j''=b_j-1
\atop
b_j'\geqslant 0,\,\,b_j''\geqslant 0}\,
\zeta_{[{\sf w},{\sf d},{\sf h}+1]}
\big(\bfdots,1^{b_j'},4,1^{b_j''},\bfdots\big)
\Bigg)
+
\\
&
\ \ \ \ \
+
\zeta_{[{\sf w},{\sf d}+1,{\sf h}+1]}
\big(3,a_1,1^{b_1},\dots,a_{\sf h},1^{b_{\sf h}}\big)
+
\\
&
\ \ \ \ \
\sum_{1\leqslant j\leqslant{\sf h}}
\Bigg(
\sum_{b_j'+b_j''=b_j
\atop
b_j'\geqslant 0,\,b_j''\geqslant 0}\,
\zeta_{[{\sf w},{\sf d}+1,{\sf h}+1]}
\big(\bfdots,1^{b_j'},3,1^{b_j''},\bfdots\big)
\Bigg).
\endaligned
\]
\end{Theorem}

\begin{Theorem}
For every height ${\sf h} \geqslant 1$,
every entries $a_1 \geqslant 2$, \dots, $a_{\sf h} \geqslant 2$,
every entries $b_1 \geqslant 0$, \dots, $b_{\sf h} \geqslant 0$,
if one sets:
\[
\aligned
a_1+b_1+\cdots+a_{\sf h}+b_{\sf h}
&
=:
{\sf w}-3,
\\
1+b_1+\cdots+1+b_{\sf h}
&
=:
{\sf d},
\endaligned
\]
then the so-called $\shuffle$-stuffle product between $\zeta ( 3)$ and
any $\zeta_{[{\sf w}-3, {\sf d}, {\sf h}]} \big( a_1, 1^{ b_1}, \dots,
a_{\sf h}, 1^{ b_{\sf h}} \big)$ of weight ${\sf w} - 3 \geqslant 2$,
of depth ${\sf d}$ and of height ${\sf h}$ writes in full:
\[
\footnotesize
\aligned
001
&
\shuffle
0^{a_1-1}1\,1^{b_1}\cdots
0^{a_{\sf h}-1}1\,1^{b_{\sf h}}
=
\\
&
\ \ \ \ \
=
\zeta_{[{\sf w},{\sf d}+1,{\sf h}+1]}
\big(3,a_1,1^{b_1},\dots,a_k,1^{b_k}\big)
+
\\
&
\ \ \ \ \ \ \ \ \ \
+
\sum_{1\leqslant i\leqslant{\sf h}\atop a_i\geqslant 3}
\Bigg(
\sum_{a_i'+a_i''=a_i+3
\atop
a_i'\geqslant 4,\,\,a_i''\geqslant 2}\,
\frac{(a_i'-1)(a_i'-2)}{2}
\cdot
\zeta_{[{\sf w},{\sf d}+1,{\sf h}+1]}
\big(
\bfdots,a_i',a_i'',\bfdots
\big)
\Bigg)
+
\endaligned
\]
\[
\footnotesize
\aligned
&
+
\sum_{1\leqslant j\leqslant{\sf h}}
\Bigg(
\sum_{b_j'+b_j''=b_j
\atop
b_j'\geqslant 0,\,\,b_j''\geqslant 0}\,
(b_j''+1)
\cdot
\zeta_{[{\sf w},{\sf d}+1,{\sf h}+1]}
\big(
\bfdots,1^{b_j'},3,1^{b_j''},\bfdots
\big)
\Bigg)
+
\\
&
\ \ \ \ \
+
\sum_{1\leqslant j\leqslant{\sf h}\atop b_j\geqslant 1}
\Bigg(
\sum_{b_j'+b_j''+b_j'''=b_j-1
\atop
b_j'\geqslant 0,\,\,b_j''\geqslant 0,\,\,b_j'''\geqslant 0}\,
(b_j'''+1)
\cdot
\zeta_{[{\sf w},{\sf d}+1,{\sf h}+1]}
\big(
\bfdots,1^{b_j'},2,1^{b_j''},2,1^{b_j'''},\bfdots
\big)
\Bigg)
+
\endaligned
\]
\[
\footnotesize
\aligned
&
+
\sum_{1\leqslant i_1<i_2\leqslant{\sf h}
\atop
a_{i_2}\geqslant 3}
\Bigg(
\sum_{a_{i_2}'+a_{i_2}''=a_{i_2}+1
\atop
a_{i_2}'\geqslant 2,\,\,a_{i_2}''\geqslant 2}\,
\frac{(a_{i_1}+1)\,a_{i_1}}{2}
\cdot
\zeta_{[{\sf w},{\sf d}+1,{\sf h}+1]}
\big(
\bfdots,a_{i_1}+2,\bfdots,a_{i_2}',a_{i_2}'',\bfdots
\big)
\Bigg)
+
\\
&
+
\sum_{1\leqslant i\leqslant j\leqslant{\sf h}}
\Bigg(
\frac{(a_i+1)\,a_i}{2}\,(b_j+2)
\cdot
\zeta_{[{\sf w},{\sf d}+1,{\sf h}]}
\big(
\bfdots,a_i+2,\bfdots,1^{b_j+1},\bfdots
\big)
\Bigg)
+
\endaligned
\]
\[
\footnotesize
\aligned
&
+
\sum_{1\leqslant j<i\leqslant{\sf h}
\atop
b_j\geqslant 1,\,\,a_i\geqslant 3}
\Bigg(
\sum_{b_j'+b_j''=b_j-1
\atop
b_j'\geqslant 0,\,\,b_j''\geqslant 0}
\sum_{a_i'+a_i''=a_i+1
\atop
a_i'\geqslant 2,\,\,a_i''\geqslant 2}\,
\zeta_{[{\sf w},{\sf d}+1,{\sf h}+2]}
\big(
\bfdots,1^{b_j'},3,1^{b_j''},\bfdots,a_i',a_i'',\bfdots
\big)
\Bigg)
+
\\
&
+
\sum_{1\leqslant j<i\leqslant{\sf h}
\atop
b_j\geqslant 2,\,\,a_i\geqslant 3}
\Bigg(
\sum_{b_j'+b_j''+b_j'''=b_j-2
\atop
b_j'\geqslant,\,\,b_j''\geqslant 0,\,\,b_j'''\geqslant 0}
\sum_{a_i'+a_i''=a_i+1
\atop
a_i'\geqslant 2,\,\,a_i''\geqslant 2}
\zeta_{[{\sf w},{\sf d}+1,{\sf h}+3]}
\big(
\bfdots,1^{b_j'},2,1^{b_j''},2,1^{b_j'''},\bfdots,
a_i',a_i'',\bfdots
\big)
\Bigg)
+
\endaligned
\]
\[
\footnotesize
\aligned
&
+
\sum_{1\leqslant j_1<j_2\leqslant{\sf h}
\atop
b_{j_1}\geqslant 1}
\Bigg(
\sum_{b_{j_1}'+b_{j_1}''=b_{j_1}-1
\atop
b_{j_1}'\geqslant 0,\,\,b_{j_1}''\geqslant 0}\,
(b_{j_2}+1)
\cdot
\zeta_{[{\sf w},{\sf d}+1,{\sf h}+1]}
\big(
\bfdots,1^{b{j_1}'},3,1^{b_{j_1}''},\bfdots,1^{b_{j_2}+1},\bfdots
\big)
\Bigg)
+
\\
&
+
\sum_{1\leqslant j_1<j_2\leqslant{\sf h}
\atop
b_{j_1}\geqslant 2}
\Bigg(
\sum_{b_{j_1}'+b_{j_1}''+b_{j_1}'''=b_{j_1}-2
\atop
b_{j_1}'\geqslant 0,\,\,b_{j_1}''\geqslant 0,\,\,b_{j_1}'''\geqslant 0}
(b_{j_2}+2)
\cdot
\zeta_{[{\sf w},{\sf d}+1,{\sf h}+2]}
\big(
\bfdots,1^{b_{j_1}'},2,1^{b_{j_1}''},2,1^{b_{j_1}'''},
\bfdots,1^{b_{j_2}+1},\bfdots
\big)
\Bigg)
+
\endaligned
\]
\[
\footnotesize
\aligned
&
+
\sum_{1\leqslant i_1<i_2\leqslant{\sf h}\atop a_{i_2}\geqslant 3}
\Bigg(
\sum_{a_{i_2}'+a_{i_2}''=a_{i_2}+2
\atop
a{i_2}'\geqslant 3,\,\,a_{i_2}''\geqslant 2}
a_{i_1}\,(a_{i_2}'-1)
\cdot
\zeta_{[{\sf w},{\sf d}+1,{\sf h}+1]}
\big(
\bfdots,a_{i_1}+1,\bfdots,a_{i_2}',a_{i_2}'',\bfdots
\big)
\Bigg)
+
\\
&
+
\sum_{1\leqslant i\leqslant j\leqslant{\sf h}}
\Bigg(
\sum_{b_j'+b_j''=b_j
\atop
b_j'\geqslant 0,\,\,b_j''\geqslant 0}\,
a_i\,(b_j''+1)
\cdot
\zeta_{[{\sf w},{\sf d}+1,{\sf h}+1]}
\big(
\bfdots,a_i+1,\bfdots,1^{b_j'},2,1^{b_j''},\bfdots
\big)
\Bigg)
+
\endaligned
\]
\[
\footnotesize
\aligned
&
+
\sum_{1\leqslant j<i\leqslant{\sf h}
\atop
b_j\geqslant 1,\,\,a_i\geqslant 3}
\Bigg(
\sum_{b_j'+b_j''=b_j-1
\atop
b_j'\geqslant 0,\,\,b_j''\geqslant 0}
\sum_{a_i'+a_i''=a_i+2
\atop
a_i'\geqslant 3,\,\,a_i''\geqslant 2}
(a_i'-1)
\cdot
\zeta_{[{\sf w},{\sf d}+1,{\sf h}+2]}
\big(
\bfdots,1^{b_j'},2,1^{b_j''},\bfdots,a_i',a_i'',\bfdots
\big)
\Bigg)
+
\\
&
+
\sum_{1\leqslant j_1<j_2\leqslant{\sf h}
\atop
b_{j_1}\geqslant 1}
\Bigg(
\sum_{b_{j_1}'+b_{j_1}''=b_{j_1}-1
\atop
b_{j_1}'\geqslant 0,\,\,b_{j_1}''\geqslant 0}
\sum_{b_{j_2}'+b_{j_2}''=b_{j_2}
\atop
b_{j_2}'\geqslant 0,\,\,b_{j_2}''\geqslant 0}
(b_{j_2}''+1)
\cdot
\zeta_{[{\sf w},{\sf d}+1,{\sf h}+2]}
\big(
\bfdots,1^{b_{j_1}'},2,1^{b_{j_1}''},\bfdots,
1^{b_{j_2}'},2,1^{b_{j_2}''},\bfdots
\big)
\Bigg)
+
\endaligned
\]
\[
\footnotesize
\aligned
&
+
\sum_{1\leqslant i_1<i_2<i_3\leqslant{\sf h}
\atop
a_{i_3}\geqslant 3}
\Bigg(
\sum_{a_{i_3}'+a_{i_3}''=a_{i_3}+1
\atop
a_{i_3}'\geqslant 2,\,\,a_{i_3}''\geqslant 2}
a_{i_1}a_{i_2}
\cdot
\zeta_{[{\sf w},{\sf d}+1,{\sf h}+1]}
\big(
\bfdots,a_{i_1}+1,\bfdots,a_{i_2}+1,\bfdots,a_{i_3}',a_{i_3}'',\bfdots
\big)
\Bigg)
+
\\
&
+
\sum_{1\leqslant i_1<i_2\leqslant j\leqslant{\sf h}}\,
a_{i_1}a_{i_2}(b_j+2)
\cdot
\zeta_{[{\sf w},{\sf d}+1,{\sf h}]}
\big(
\bfdots,a_{i_1}+1,\bfdots,a_{i_2}+1,\bfdots,1^{b_j+1},\bfdots
\big)
+
\endaligned
\]
\[
\footnotesize
\aligned
&
+
\sum_{1\leqslant i_1\leqslant j<i_2\leqslant{\sf h}
\atop
b_j\geqslant 1,\,\,a_{i_2}\geqslant 3}
\Bigg(
\sum_{b_j'+b_j''=b_j-1
\atop
b_j'\geqslant 0,\,\,b_j''\geqslant 0}
\sum_{a_{i_2}'+a_{i_2}''=a_{i_2}+1
\atop
a_{i_2}'\geqslant 2,\,\,a_{i_2}''\geqslant 2}
a_{i_1}
\cdot
\zeta_{[{\sf w},{\sf d}+1,{\sf h}+2]}
\big(
\bfdots,a_{i_1}+1,\bfdots,1^{b_j'},2,1^{b_j''},\bfdots,
a_{i_2}',a_{i_2}'',\bfdots
\big)
\Bigg)
+
\\
&
+
\sum_{1\leqslant i\leqslant j_1<j_2\leqslant{\sf h}
\atop
b_{j_1}\geqslant 1}
\Bigg(
\sum_{b_{j_1}'+b_{j_1}''=b_{j_1}-1
\atop
b_{j_1}'\geqslant 0,\,\,b_{j_1}''\geqslant 0}
a_i(b_{j_2}+2)
\cdot
\zeta_{[{\sf w},{\sf d}+1,{\sf h}+1]}
\big(
\bfdots,a_i+1,\bfdots,1^{b_j'},2,1^{b_{j_1}''},\bfdots,1^{b_{j_2}+1},
\bfdots
\big)
\Bigg)
+
\endaligned
\]
\[
\footnotesize
\aligned
&
+
\sum_{1\leqslant j<i_1<i_2\leqslant{\sf h}
\atop
b_j\geqslant 1,\,\,a_{i_2}\geqslant 3}
\Bigg(
\sum_{b_j'+b_j''=b_j-1
\atop
b_j'\geqslant 0,\,\,b_j''\geqslant 0}
\sum_{a_{i_2}'+a_{i_2}''=a_{i_2}+1
\atop
a_{i_2}'\geqslant 2,\,\,a_{i_2}''\geqslant 2}
\zeta_{[{\sf w},{\sf d}+1,{\sf h}+2]}
\big(
\bfdots,1^{b_j'},2,1^{b_j''},\bfdots,
a_{i_1}+1,\bfdots,a_{i_2}',a_{i_2}'',\bfdots
\big)
\Bigg)
+
\\
&
+
\sum_{1\leqslant j_1<i\leqslant j_2\leqslant{\sf h}
\atop
b_{j_1}\geqslant 1}
\Bigg(
\sum_{b_{j_1}'+b_{j_1}''=b_{j_1}-1
\atop
b_{j_1'}\geqslant 0,\,\,b_{j_1}''\geqslant 0}
a_i(b_{j_2}+2)
\cdot
\zeta_{[{\sf w},{\sf d}+1,{\sf h}+1]}
\big(
\bfdots,1^{b_{j_1}'},2,1^{b_{j_1}''},\bfdots,
a_i+1,\bfdots,1^{b_{j_2}+1},\bfdots
\big)
\Bigg)
+
\endaligned
\]
\[
\footnotesize
\aligned
&
+
\sum_{1\leqslant j_1<j_2<i\leqslant{\sf h}
\atop
b_{j_1}\geqslant 1,\,b_{j_2}\geqslant 1,\,a_i\geqslant 3}\!\!\!
\Bigg(
\sum_{b_{j_1}'+b_{j_1}''=b_{j_1}-1
\atop
b_{j_1}'\geqslant 0,\,\,b_{j_1}''\geqslant 0}
\sum_{b_{j_2}'+b_{j_2}''=b_{j_2}-1
\atop
b_{j_2}'\geqslant 0,\,\,b_{j_2}''\geqslant 0}
\sum_{a_i'+a_i''=a_i+1
\atop
a_i'\geqslant 2,\,\,a_i''\geqslant 2}
\\
&
\ \ \ \ \ \ \ \ \ \ \ \ \ \ \ \ \ \ \ \ \ \ \ \ \ \ \ \ \ \ \ \ \ \ \ 
\ \ \ \ \ \ \ \ \ \ \ \ \ \ \ \ \ \ \ \ \ \ \ \ \ \ \ \ 
\zeta_{[{\sf w},{\sf d}+1,{\sf h}+3]}
\big(
\bfdots,1^{b_{j_1}'},2,1^{b_{j_1}''},\bfdots,
1^{b_{j_2}'},2,1^{b_{j_2}''},\bfdots,
a_i',a_i'',\bfdots
\big)
\Bigg)
+
\endaligned
\]
\[
\footnotesize
\aligned
&
+
\sum_{1\leqslant j_1<j_2<j_3\leqslant{\sf h}
\atop
b_{j_1}\geqslant 1,\,\,b_{j_2}\geqslant 1}
\Bigg(
\sum_{b_{j_1}'+b_{j_1}''=b_{j_1}-1
\atop
b_{j_1}'\geqslant 0,\,\,b_{j_1}''\geqslant 0}
\sum_{b_{j_2}'+b_{j_2}''=b_{j_2}-1
\atop
b_{j_2}'\geqslant 0,\,\,b_{j_2}''\geqslant 0}
\\
&
\ \ \ \ \ \ \ \ \ \ \ \ \ \ \ \ \ \ \ \ \ \ \ \ \ \ \ \ \ \ \ \ \ \ \ 
\ \ \ \ \ \ \ \ \ \ 
(b_{j_3}+2)
\cdot
\zeta_{[{\sf w},{\sf d}+1,{\sf h}+2]}
\big(
\bfdots,1^{b_{j_1}'},2,1^{b_{j_1}''},\bfdots,
1^{b_{j_2}'},2,1^{b_{j_2}''},\bfdots,
1^{b_{j_3}+1},\bfdots
\big)
\Bigg)
+
\endaligned
\]
\end{Theorem}

\begin{Theorem}
For every height ${\sf h} \geqslant 1$,
every entries $a_1 \geqslant 2$, \dots, $a_{\sf h} \geqslant 2$,
every entries $b_1 \geqslant 0$, \dots, $b_{\sf h} \geqslant 0$,
if one sets:
\[
\aligned
a_1+b_1+\cdots+a_{\sf h}+b_{\sf h}
&
=:
{\sf w}-3,
\\
1+b_1+\cdots+1+b_{\sf h}
&
=:
{\sf d},
\endaligned
\]
then the so-called $\ast$-stuffle product between $\zeta (2,1)$ and
any $\zeta_{[{\sf w}-3, {\sf d}, {\sf h}]} \big( a_1, 1^{ b_1}, \dots,
a_{\sf h}, 1^{ b_{\sf h}} \big)$ of weight ${\sf w} - 3
\geqslant 2$, of depth
${\sf d}$ and of height ${\sf h}$ writes in full:
\[
\aligned
\zeta(2,1)
&
\ast
\zeta_{[{\sf w}-3,{\sf d},{\sf h}]}
\big(a_1,1^{b_1},\dots,a_{\sf h},1^{b_{\sf h}}\big)
=
\\
&
=
\sum_{1\leqslant i_1< i_2\leqslant{\sf h}}
\zeta_{[{\sf w},{\sf d},{\sf h}]}
\big(
\bfdots,a_{i_1}+2,\bfdots,a_{i_2}+1,\bfdots
\big)
+
\endaligned
\]
\[
\footnotesize
\aligned
&
+
\sum_{1\leqslant i\leqslant j\leqslant{\sf h}
\atop
b_j\geqslant 1}
\Bigg(
\sum_{b_j'+b_j''=b_j-1
\atop
b_j'\geqslant 0,\,\,b_j''\geqslant 0}
\zeta_{[{\sf w},{\sf d},{\sf h}+1]}
\big(
\bfdots,a_i+2,\bfdots,1^{b_j'},2,1^{b_j''},\bfdots
\big)
\Bigg)
+
\\
&
+
\sum_{1\leqslant j<i\leqslant{\sf h}
\atop
b_j\geqslant 1}
\Bigg(
\sum_{b_j'+b_j''=b_j-1
\atop
b_j'\geqslant 0,\,\,b_j''\geqslant 0}
\zeta_{[{\sf w},{\sf d},{\sf h}+1]}
\big(
\bfdots,1^{b_j'},3,1^{b_j''},\bfdots,a_i+1,\bfdots
\big)
\Bigg)
+
\endaligned
\]
\[
\footnotesize
\aligned
&
+
\sum_{1\leqslant j_1<j_2\leqslant{\sf h}}
\Bigg(
\sum_{b_{j_1}'+b_{j_1}''=b_{j_1}-1
\atop
b_{j_1}'\geqslant 0,\,\,b_{j_1}''\geqslant 0}
\sum_{b_{j_2}'+b_{j_2}''=b_{j_2}-1
\atop
b_{j_2}'\geqslant 0,\,\,b_{j_2}''\geqslant 0}
\zeta_{[{\sf w},{\sf d},{\sf h}+2]}
\big(
\bfdots,1^{b_{j_1}'},3,1^{b_{j_1}''},\bfdots,
1^{b_{j_2}'},2,1^{b_{j_2}''},\bfdots
\big)
\Bigg)
+
\\
&
+
\sum_{1\leqslant i\leqslant j\leqslant{\sf h}}
(b_j+1)
\cdot
\zeta_{[{\sf w},{\sf d}+1,{\sf h}]}
\big(
\bfdots,a_i+2,\bfdots,1^{b_j+1},\bfdots
\big)
+
\endaligned
\]
\[
\footnotesize
\aligned
&
+
\sum_{1\leqslant j\leqslant{\sf h}\atop b_j\geqslant 1}
\Bigg(
\sum_{b_j'+b_j''=b_j-1
\atop
b_j'\geqslant 0,\,\,b_j''\geqslant 0}
(b_j''+1)
\cdot
\zeta_{[{\sf w},{\sf d}+1,{\sf h}+1]}
\big(
\bfdots,1^{b_j'},3,1^{b_j''},\bfdots
\big)
\Bigg)
+
\\
&
+
\sum_{1\leqslant j_1<j_2\leqslant{\sf h}\atop b_{j_1}\geqslant 1}
\Bigg(
\sum_{b_{j_1}'+b_{j_1}''=b_{j_1}-1
\atop
b_{j_1}'\geqslant 0,\,\,b_{j_1}''\geqslant 0}
(b_{j_2}+1)
\cdot
\zeta_{[{\sf w},{\sf d}+1,{\sf h}+1]}
\big(
\bfdots,1^{b_{j_1}'},3,1^{b_{j_1}''},\bfdots,1^{b_{j_2}+1},\bfdots
\big)
\Bigg)
+
\endaligned
\]
\[
\footnotesize
\aligned
&
+
\sum_{1\leqslant j\leqslant{\sf h}
\atop
b_j\geqslant 1}
\Bigg(
\sum_{b_j'+b_j''=b_j-1
\atop
b_j'\geqslant 0,\,\,b_j''\geqslant 0}
\zeta_{[{\sf w},{\sf d}+1,{\sf h}+2]}
\big(
2,\bfdots,1^{b_j'},2,1^{b_j''},\bfdots
\big)
\Bigg)
+
\\
&
+
\sum_{1\leqslant j\leqslant{\sf h}}
\Bigg(
\sum_{b_j'+b_j''+b_j'''=b_j
\atop
b_j'\geqslant 0,\,\,b_j''\geqslant 0,\,\,b_j'''\geqslant 0}
\zeta_{[{\sf w},{\sf d}+1,{\sf h}+2]}
\big(
\bfdots,1^{b_j'},2,1^{b_j''},2,1^{b_j'''},\bfdots
\big)
\Bigg)
+
\endaligned
\]
\[
\footnotesize
\aligned
&
+
\sum_{1\leqslant j_1<j_2\leqslant{\sf h}\atop b_{j_2}\geqslant 1}
\Bigg(
\sum_{b_{j_1}'+b_{j_1}''=b_{j_1}
\atop
b_{j_1}'\geqslant 0,\,\,b_{j_1}''\geqslant 0}
\sum_{b_{j_2}'+b_{j_2}''=b_{j_2}-1
\atop
b_{j_2}'\geqslant 0,\,\,b_{j_2}''\geqslant 0}
\zeta_{[{\sf w},{\sf d}+1,{\sf h}+2]}
\big(
\bfdots,1^{b_{j_1}'},2,1^{b_{j_1}''},\bfdots,
1^{b_{j_2}'},2,1^{b_{j_2}''},\bfdots
\big)
\Bigg)
+
\\
&
+
\sum_{1\leqslant i\leqslant{\sf h}}
\zeta_{[{\sf w},{\sf d}+1,{\sf h}+1]}
\big(
2,\bfdots,a_i+1,\bfdots
\big)
\Bigg)
+
\endaligned
\]
\[
\footnotesize
\aligned
&
+
\sum_{1\leqslant j<i\leqslant{\sf h}}
\Bigg(
\sum_{b_j'+b_j''=b_j
\atop
b_j'\geqslant 0,\,\,b_j''\geqslant 0}
\zeta_{[{\sf w},{\sf d}+1,{\sf h}+1]}
\big(
\bfdots,1^{b_j'},2,1^{b_j''},\bfdots,a_i+1,\bfdots
\big)
\Bigg)
+
\endaligned
\]
\[
\footnotesize
\aligned
&
+
\sum_{1\leqslant j\leqslant{\sf h}}
\Bigg(
\sum_{b_j'+b_j''+b_j'''=b_j
\atop
b_j'\geqslant 0,\,\,b_j''\geqslant 0,\,\,b_j'''\geqslant 0}
\zeta_{[{\sf w},{\sf d}+2,{\sf h}+1]}
\big(
\bfdots,1^{b_j'},2,1^{b_j''},1^{b_j'''+1},\bfdots
\big)
\Bigg)
+
\\
&
+
\sum_{1\leqslant j_1<j_2\leqslant{\sf h}}
\Bigg(
\sum_{b_{j_1}'+b_{j_1}''=b_{j_1}
\atop
b_{j_1}'\geqslant 0,\,\,b_{j_1}''\geqslant 0}
(b_{j_2}+1)
\cdot
\zeta_{[{\sf w},{\sf d}+2,{\sf h}+1]}
\big(
\bfdots,1^{b_{j_1}'},2,1^{b_{j_1}''},\bfdots,1^{b_{j_2}+1},\bfdots
\big)
\Bigg)
\\
&
+
\zeta_{[{\sf w},{\sf d}+2,{\sf h}+1]}
\big(2,1,\bfdots\big).
\endaligned
\]
\end{Theorem}

\begin{Theorem}
For every height ${\sf h} \geqslant 1$,
every entries $a_1 \geqslant 2$, \dots, $a_{\sf h} \geqslant 2$,
every entries $b_1 \geqslant 0$, \dots, $b_{\sf h} \geqslant 0$,
if one sets:
\[
\aligned
a_1+b_1+\cdots+a_{\sf h}+b_{\sf h}
&
=:
{\sf w}-3,
\\
1+b_1+\cdots+1+b_{\sf h}
&
=:
{\sf d},
\endaligned
\]
then the so-called $\shuffle$-shuffle product between $\zeta (2,1)$ and
any $\zeta_{[{\sf w}-3, {\sf d}, {\sf h}]} \big( a_1, 1^{ b_1}, \dots,
a_{\sf h}, 1^{ b_{\sf h}} \big)$ of weight ${\sf w} - 3
\geqslant 2$, of depth
${\sf d}$ and of height ${\sf h}$ writes out in full:
\[
\footnotesize
\aligned
011
&
\shuffle
0^{a_1-1}1\,1^{b_1}\cdots
0^{a_{\sf h}-1}1\,1^{b_{\sf h}}
=
\endaligned
\]
\[
\footnotesize
\aligned
&
=
\sum_{1\leqslant i\leqslant{\sf h}\atop a_i\geqslant 2}
\Bigg(
\sum_{a_i'+a_i''=a_i+2
\atop
a_i'\geqslant 2,\,\,a_i''\geqslant 2}\,
(a_i'-1)
\cdot
\zeta_{[{\sf w},{\sf d}+2,{\sf h}+1]}
\big(
\bfdots,a_i',1,a_i'',\bfdots
\big)
\Bigg)
+
\\
&
+
\sum_{1\leqslant i\leqslant{\sf h}\atop a_i\geqslant 3}
\Bigg(
\sum_{a_i'+a_i''+a_i'''=a_i+3
\atop
a_i'\geqslant 2,\,\,a_i''\geqslant 2,\,\,a_i'''\geqslant 2}\,
(a_i'-1)
\cdot
\zeta_{[{\sf w},{\sf d}+2,{\sf h}2]}
\big(
\bfdots,a_i',a_i'',a_i''',\bfdots
\big)
\Bigg)
+
\endaligned
\]
\[
\footnotesize
\aligned
&
+
\sum_{1\leqslant j\leqslant{\sf h}\atop b_j\geqslant 1}
\Bigg(
\sum_{b_j'+b_j''=b_j-1
\atop
b_j'\geqslant 0,\,\,b_j''\geqslant 0}\,
\frac{(b_j''+3)(b_j''+2)}{2}
\cdot
\zeta_{[{\sf w},{\sf d}+2,{\sf h}+1]}
\big(
\bfdots,1^{b_j'},2,1^{b_j''+2},\bfdots
\big)
\Bigg)
+
\\
&
+
\sum_{1\leqslant i_1<i_2\leqslant{\sf h}\atop a_{i_2}\geqslant 3}
\Bigg(
\sum_{a_{i_1}'+a_{i_1}''=a_{i_1}+2}\,
(a_{i_1}'-1)
\cdot
\zeta_{[{\sf w},{\sf d}+2,{\sf h}+2]}
\big(
\bfdots,a_{i_1}',a_{i_1}'',\bfdots,a_{i_2}',a_{i_2}'',\bfdots
\big)
\Bigg)
+
\endaligned
\]
\[
\footnotesize
\aligned
&
+
\sum_{1\leqslant i\leqslant j\leqslant{\sf h}}
\Bigg(
\sum_{a_i'+a_i''=a_i+2
\atop
a_i'\geqslant 2,\,\,a_i''\geqslant 2}\,
(a_i'-1)(b_j+2)
\cdot
\zeta_{[{\sf w},{\sf d}+2,{\sf h}+1]}
\big(
\bfdots,a_i',a_i'',\bfdots,1^{b_{j_2}+1},\bfdots
\big)
\Bigg)
+
\\
&
+
\sum_{1\leqslant j<i\leqslant{\sf h}
\atop
b_j\geqslant 1,\,\,a_i\geqslant 3}
\Bigg(
\sum_{b_j'+b_j''=b_j-1
\atop
b_j'\geqslant 0,\,\,b_j''\geqslant 0}
\sum_{a_i'+a_i''=a_i+1
\atop
a_i'\geqslant 2,\,\,a_i''\geqslant 2}\,
(b_j''+1)
\cdot
\zeta_{[{\sf w},{\sf d}+2,{\sf h}+2]}
\big(
\bfdots,1^{b_j'},2,1^{b_j''+1},\bfdots,a_i',a_i'',\bfdots
\big)
\Bigg)
+
\endaligned
\]
\[
\footnotesize
\aligned
&
+
\sum_{1\leqslant j_1<j_2\leqslant{\sf h}
\atop
b_{j_1}\geqslant 1}
\Bigg(
\sum_{b_{j_1}'+b_{j_1}''=b_{j_1}-1
\atop
b_{j_1}'\geqslant 0,\,\,b_{j_1}''\geqslant 0}\,
(b_{j_1}''+2)(b_{j_2}+2)
\cdot
\zeta_{[{\sf w},{\sf d}+2,{\sf h}+1]}
\big(
\bfdots,1^{b_{j_1}'},2,1^{b_{j_1}''+1},
\bfdots,1^{b_{j_2}+1},\bfdots
\big)
\Bigg)
+
\\
&
+
\sum_{1\leqslant i_1<i_2\leqslant{\sf h}
\atop
a_{i_2}\geqslant 3}
\Bigg(
\sum_{a_{i_2}'+a_{i_2}''=a_{i_2}+1
\atop
a_{i_2}'\geqslant 2,\,\,a_{i_2}''\geqslant 2}\,
a_{i_1}
\cdot
\zeta_{[{\sf w},{\sf d}+2,{\sf h}+1]}
\big(
\bfdots,a_{i_1}+1,\bfdots,a_{i_2}',1,a_{i_2}'',\bfdots
\big)
\Bigg)
+
\endaligned
\]
\[
\footnotesize
\aligned
&
+
\sum_{1\leqslant i_1<i_2\leqslant{\sf h}
\atop
a_{i_2}\geqslant4}
\Bigg(
\sum_{a_{i_2}'+a_{i_2}''+a_{i_2}'''=a_{i_2}+2
\atop
a_{i_2}'\geqslant 2,\,\,a_{i_2}''\geqslant 2,\,\,a_{i_2}'''\geqslant 2}\,
a_{i_1}
\cdot
\zeta_{[{\sf w},{\sf d}+2,{\sf h}+2]}
\big(
\bfdots,a_{i_1}+1,\bfdots,a_{i_2}',a_{i_2}'',a_{i_2}''',\bfdots
\big)
\Bigg)
+
\\
&
+
\sum_{1\leqslant i\leqslant j\leqslant{\sf h}}\,
a_i\,
\frac{(b_j+3)(b_j+2)}{2}
\cdot
\zeta_{[{\sf w},{\sf d}+2,{\sf h}]}
\big(
\bfdots,a_i+1,\bfdots,1^{b_j+2},\bfdots
\big)
+
\endaligned
\]
\[
\footnotesize
\aligned
&
+
\sum_{1\leqslant j<i\leqslant{\sf h}
\atop
b_j\geqslant 1,\,\,a_i\geqslant 3}
\Bigg(
\sum_{b_j'+b_j''=b_j-1
\atop
b_j'\geqslant 0,\,\,b_j''\geqslant 0}
\sum_{a_i'+a_i''=a_i+1
\atop
a_i'\geqslant 2,\,\,a_i''\geqslant 2}
\zeta_{[{\sf w},{\sf d}+2,{\sf h}+2]}
\big(
\bfdots,1^{b_j'},2,1^{b_j''},\bfdots,
a_i',1,a_i'',\bfdots
\big)
\Bigg)
+
\\
&
+
\sum_{1\leqslant j<i\leqslant{\sf h}
\atop
b_j\geqslant 1,\,\,a_i\geqslant 4}
\Bigg(
\sum_{b_j'+b_j''=b_j-1
\atop
b_j'\geqslant 0,\,\,b_j''\geqslant 0}
\sum_{a_i'+a_i''+a_i'''=a_i+2
\atop
a_i'\geqslant 2,\,\,a_i''\geqslant 2,\,\,a_i'''\geqslant 2}
\zeta_{[{\sf w},{\sf d}+2,{\sf h}+3]}
\big(
\bfdots,1^{b_j'},2,1^{b_j''},\bfdots,a_i',a_i'',a_i''',\bfdots
\big)
\Bigg)
+
\endaligned
\]
\[
\footnotesize
\aligned
&
+
\sum_{1\leqslant j_1<j_2\leqslant{\sf h}
\atop
b_{j_1}\geqslant 1}
\Bigg(
\sum_{b_{j_1}'+b_{j_1}''=b_{j_1}-1
\atop
b_{j_1}'\geqslant 0,\,\,b_{j_1}''\geqslant 0}\,
\frac{(b_{j_2}+3)(b_{j_2}+2)}{2}
\cdot
\zeta_{[{\sf w},{\sf d}+2,{\sf h}+1]}
\big(
\bfdots,1^{b_{j_1}'},2,1^{b_{j_1}''},\bfdots,
1^{b_{j_2}+2},\bfdots
\big)
\Bigg)
+
\\
&
+
\sum_{1\leqslant i_1<i_2<i_3\leqslant{\sf h}
\atop
a_{i_2}\geqslant 3,\,\,a_{i_3}\geqslant 3}
\Bigg(
\sum_{a_{i_2}'+a_{i_2}''=a_{i_2}+1
\atop
a_{i_2}'\geqslant 2,\,\,a_{i_2}''\geqslant 2}
\sum_{a_{i_3}'+a_{i_3}''=a_{i_3}+1
\atop
a_{i_3}'\geqslant 2,\,\,a_{i_3}''\geqslant 2}\,
\\
&
\ \ \ \ \ \ \ \ \ \ \ \ \ \ \ \ \ \ \ \ \ \ \ \ \ \ \ \ \ \ \ \ \ \
\ \ \ \ \ \ \ \ \ \ \ \ \ \ \ \ \ \ \ \ \ \ \ \ 
a_{i_1}
\cdot
\zeta_{[{\sf w},{\sf d}+2,{\sf h}+2]}
\big(
\bfdots,a_{i_1}+1,\bfdots,a_{i_2}',a_{i_2}'',\bfdots,
a_{i_3}',a_{i_3}'',\bfdots
\big)
\Bigg)
+
\endaligned
\]
\[
\footnotesize
\aligned
&
+
\sum_{1\leqslant i_1<i_2\leqslant j\leqslant{\sf h}
\atop
a_{i_2}\geqslant 3}
\Bigg(
\sum_{a_{i_2'}+a_{i_2}''=a_{i_2}+1
\atop
a_{i_2}'\geqslant 2,\,\,a_{i_2}''\geqslant 2}\,
a_{i_1}(b_j+2)
\cdot
\zeta_{[{\sf w},{\sf d}+2,{\sf h}+1]}
\big(
\bfdots,a_{i_1}+1,\bfdots,a_{i_2}',a_{i_2}'',\bfdots,
1^{b_j+1},\bfdots
\big)
\Bigg)
+
\\
&
+
\sum_{1\leqslant i_1\leqslant j<i_2\leqslant{\sf h}}
\Bigg(
\sum_{a_{i_2}'+a_{i_2}''=a_{i_2}+1
\atop
a_{i_2}'\geqslant 2,\,\,a_{i_2}''\geqslant 2}\,
a_{i_1}(b_j+2)
\cdot
\zeta_{[{\sf w},{\sf d}+2,{\sf h}+1]}
\big(
\bfdots,a_{i_1}+1,\bfdots,1^{b_j+1},\bfdots,
a_{i_2}',a_{i_2}'',\bfdots
\big)
\Bigg)
+
\endaligned
\]
\[
\footnotesize
\aligned
&
+
\sum_{1\leqslant i\leqslant j_1<j_2\leqslant{\sf h}}\,
a_i(b_{j_1}+2)(b_{j_2}+2)
\cdot
\zeta_{[{\sf w},{\sf d}+2,{\sf h}]}
\big(
\bfdots,a_i+1,\bfdots,1^{b_{j_1}+1},\bfdots,
1^{b_{j_2}+1},\bfdots
\big)
+
\\
&
+
\sum_{1\leqslant j<i_1<i_2\leqslant{\sf h}
\atop a_{i_1}\geqslant 3,\,\,a_{i_2}\geqslant 3}
\Bigg(
\sum_{b_j'+b_j''=b_j-1
\atop
b_j'\geqslant 0,\,\,b_j''\geqslant 0}
\sum_{a_{i_1}'+a{i_1}''=a_{i_1}+1
\atop
a_{i_1}'\geqslant 2,\,\,a_{i_1}''\geqslant 2}
\sum_{a_{i_2}'+a{i_2}''=a_{i_2}+1
\atop
a_{i_2}'\geqslant 2,\,\,a_{i_2}''\geqslant 2}
\\
&
\ \ \ \ \ \ \ \ \ \ \ \ \ \ \ \ \ \ \ \ \ \ \ \ \ \ \ \ \ \ \ \ \ \
\ \ \ \ \ \ \ \ \ \ \ \ \ \ \ \ \ \ \ \ \ \ \ \ 
\zeta_{[{\sf w},{\sf d}+2,{\sf h}+3]}
\big(
\bfdots,1^{b_j'},2,1^{b_j''},\bfdots,
a_{i_1}',a_{i_1}'',\bfdots,a_{i_2}',a_{i_2}'',\bfdots
\big)
\Bigg)
+
\endaligned
\]
\[
\footnotesize
\aligned
&
+
\sum_{1\leqslant j_1<i\leqslant j_2\leqslant{\sf h}}
\Bigg(
\sum_{b_{j_1}'+b_{j_1}''=b_{j_1}-1
\atop
b_{j_1}'\geqslant 0,\,\,b_{j_1}''\geqslant 0}
\sum_{a_i'+a_i''=a_i+1
\atop
a_i'\geqslant 2,\,\,a_i''\geqslant 2}\,
\\
&
\ \ \ \ \ \ \ \ \ \ \ \ \ \ \ \ \ \ \ \ \ \ \ \ \ \ \ \ \ \ \ \ \ \
\ \ \ \ \ \ \ 
(b_{j_2}+2)
\cdot
\zeta_{[{\sf w},{\sf d}+2,{\sf h}+2]}
\big(
\bfdots,1^{b_{j_1}'},2,1^{b_{j_1}''},\bfdots,
a_i',a_i'',\bfdots,1^{b_{j_2}+1},\bfdots
\big)
\Bigg)
+
\\
&
+
\sum_{1\leqslant j_1<j_2<i\leqslant{\sf h}}
\Bigg(
\sum_{b_{j_1}'+b_{j_1}''=b_{j_1}-1
\atop
b_{j_1}'\geqslant 0,\,\,b_{j_1}''\geqslant 0}
\sum_{a_i'+a_i''=a_i+1
\atop
a_i'\geqslant 2,\,\,a_i''\geqslant 2}\,
\\
&
\ \ \ \ \ \ \ \ \ \ \ \ \ \ \ \ \ \ \ \ \ \ \ \ \ \ \ \ \ \ \ \ \ \
\ \ \ \ \ \ \ \ \ \
(b_{j_2}+2)
\cdot
\zeta_{[{\sf w},{\sf d}+2,{\sf h}+2]}
\big(
\bfdots,1^{b_{j_1}'},2,1^{b_{j_1}''},\bfdots,
1^{b_{j_2}+1},\bfdots,a_i',a_i'',\bfdots
\big)
\Bigg)
+
\endaligned
\]
\[
\footnotesize
\aligned
&
+
\sum_{1\leqslant j_1<j_2<j_3\leqslant{\sf h}}
\Bigg(
\sum_{b_{j_1}'+b_{j_1}''=b_{j_1}-1
\atop
b_{j_1}'\geqslant 0,\,\,b_{j_1}''\geqslant 0}\,
\\
&
\ \ \ \ \ \ \ \ \ \ \ \ \ \ \ \ \ \ \ \ \ \ \ \ \ \ \ \ \ \ \ 
(b_{j_2}+2)(b_{j_3}+2)
\cdot
\zeta_{[{\sf w},{\sf d}+2,{\sf h}+1]}
\big(
\bfdots,1^{b_{j_1}'},2,1^{b_{j_1}''},\bfdots,
1^{b_{j_2}+1},\bfdots,1^{b_{j_3}+1},\bfdots
\big)
\Bigg)
+
\\
&
+
\zeta_{[{\sf w},{\sf d}+2,{\sf h}+1]}
\big(\bfdots,2,1).
\endaligned
\]
\end{Theorem}

\subsection{Acknowledgments and research `trail'}
Once in January 2011, at the {\em S\'eminaire de Philosophie des
Math\'ematiques} (\'Ecole Normale Sup\'erieure, Paris), Pierre
Cartier (\cite{ Cartier-2011}) expressed the natural desire that a
constructive approach existed towards Hoffman's $(2,3)$-conjecture,
and this motivated us to invest energy in order to set up usable
general formulas for double shuffle relations. Quickly, we realized
that presumably, the four families on p.~\pageref{1-2-3-21}
of double shuffle relations
could suffice, we checked this manually up to weight ${\sf w} = 7$,
and we also spent enough time in Spring 2011 to finalize all the
formulas presented here. However, after a while of further experiments
and explorations (the content of which will appear in some forthcoming
{\footnotesize\sf arxiv.org} prepublications), a prudent mathematical
wise fear appeared in our thoughts that possibly, some hidden
unpleasant complexities could probably contradict the sufficiency 
of the four families of 
relations on p.~\pageref{1-2-3-21}.
It was only in July 2012 that the {\footnotesize\sf arxiv.org}
preprint~\cite{ Combariza-2012} of German Combariza made us aware that
the improved $(2, 3)$-conjecture we devised intuitively without
conjecturable certainty had already been verified by Masanobu Kaneko,
Masayuki Noro and Ken'ichi Tsurumaki to be correct on powerful
computer machines up to weight ${\sf w} = 20$! This wonderful confirmation
decided us to typeset our manuscript.

Special thanks are addressed to Professor Masanobu Kaneko for sending
us a pdf file of the publication~\cite{ MMK-2008} and for transmitting
data that are unreachable by hand.  Further interactions between mind
and machine are expectable, hopefully towards a full constructive
resolution of various polyzeta conjectures in arbitrary weight.

In July 2012, Olivier Bouillot provided us with precious Maple files
which automatize the production of double-shuffle relations, a tool
we are very grateful for.

Of course, the (extensive) formulas provided in this article are only
a preliminary\,\,---\,\,in a sense easy and elementary\,\,---\,\,step
towards constructiveness in the $(2, 3)$-conjecture, for experts know
well that double shuffle computations just lie on the threshold of
higher level linear algebra computations with matrix minors whose
size, unfortunately, increases exponentially with the weights of the
polyzetas.

\section{Weight, depth, height}
\label{weight-depth-height}
\HEAD{Weight, Depth, Height}{
Jo\"el {\sc Merker}, Universit\'e Paris-Sud Orsay, France}

\subsection{Principles of $\zeta$-notations}
Thus, any convergent polyzeta writes out:
\[
\zeta\big(s_1,s_2,\dots,s_{\sf d}\big)
=
\sum_{n_1>n_2>\cdots>n_{\sf d}>0}\,
\frac{1}{(n_1)^{s_1}\,(n_2)^{s_2}\,\cdots\,
(n_{\sf d})^{s_{\sf d}}},
\]
for certain integers $s_1, s_2, \dots, s_{\sf d} \geqslant 1$,
all positive, with
the additional requirement that $s_1 \geqslant 2$ in order to insure
(absolute) convergence of this multiple series in which
$n_1 > n_2 > \cdots > n_{\sf d} > 0$ are integers.

\begin{Definition}
Classically, the integer ${\sf d} \geqslant 1$ is called the {\sl
depth} ${\sf d} ( \zeta)$ 
of the polyzeta $\zeta\big( s_1, s_2, \dots, s_{\sf d}
\big)$, while the total sum of the $s_i$:
\[
{\sf w}
:=
s_1+s_2+\cdots+s_{\sf d}
\]
is said to be its {\sl weight} ${\sf w}(\zeta)$.
\end{Definition}

Here are a few simple instances:
\[
\zeta(2,1),
\ \ \ \ \ \
\zeta(3,2),
\ \ \ \ \ \
\zeta(6,2),
\ \ \ \ \ \
\zeta(2,1,2,1,1),
\ \ \ \ \ \
\zeta(5,3,1),
\]
the weights and the depths of which may be at once computed mentally.

Although it appears in the
literature to be a less widespread feature, the concept of {\sl height} 
is also intrinsically necessary to deal with when one enters
a {\em more-in-depth}, general, systematic {\sl Polyzeta Calculus}.

\begin{Definition}
The {\sl height} of a convergent polyzeta $\zeta ( s_1, \dots, s_{\sf d})$ 
is the number of entries $s_i$ that are $\geqslant 2$:
\[
{\sf h}
(\zeta)
:=
\mathmotsf{Card}\,
\big\{
i\in\{1,2,\dots,{\sf d}(\zeta)\}\,
\colon\
s_i\geqslant 2
\big\}.
\]
\end{Definition}

Of course, the height is always $\geqslant 1$.
Equivalently, one could just count the number of entries 
that are equal to $1$.
It is appropriate, then, to adapt the notation in order to 
specify how many entries $1$ are repeated one after the other.
Here are a few examples:
\[
\aligned
\zeta(3,1,2,1,1)
&
=
\zeta(3,1,2,1^2),
\\
\zeta(4,1,1,1)
&
=
\zeta(4,1^3),
\\
\zeta(5,3,2)
&
=
\zeta(5,1^0,3,1^0,2,1^0),
\endaligned
\] 
and such {\em coincidences of notations} will be admitted throughout
the present paper. With the intention of lightening the writing, we
want to denote the number of $1$ that are present in some continuous
succession, say $b \geqslant 0$ times, using just by a plain exponent,
generally as:
\[
1^b.
\] 
However, we will not use any further specific symbol like
$1^{\# b}$, $1^{*2}$, $1^{\sim 1}$. With such a convention, it must
then be clear that $1^0$ means that no $1$ is present at all.
Furthermore, between any two successive entries that are $\geqslant
2$, as {\em e.g.} in $\zeta ( 6, 2)$, an $1^0$ is always implicitly
present, and a last $1^0$ is also present when the very last $s_{\sf
d}$ happens to be $\geqslant 2$, as for instance in:
\[
\zeta(6,2)
=
\zeta(6,1^0,2,1^0).
\]

\begin{Notation}
A general, arbitrary convergent polyzeta will always be written under the
specific form:
\[
\boxed{
\zeta\big(a_1,1^{b_1},a_2,1^{b_2},\dots,
a_{\sf h},1^{b_{\sf h}}\big)}\,,
\]
where the integer ${\sf h} \geqslant 1$ is its
height, where all the $a$-integers are $\geqslant 2$:
\[
a_1\geqslant 2,
\ \ \ \ \ \ \ \ 
a_2\geqslant 2,
\ \ \ \ \ \ \ \
\dots\dots,
\ \ \ \ \ \ \ \
a_{\sf h}\geqslant 2,
\]
and where all the $b$-integers are $\geqslant 0$:
\[
b_1\geqslant 0,
\ \ \ \ \ \ \ \
b_2\geqslant 0,
\ \ \ \ \ \ \ \
\dots\dots,
\ \ \ \ \ \ \ \
b_{\sf h}\geqslant 0.
\]
\end{Notation}

In summary, {\em three fundamental numerical quantities are attached
to any convergent polyzeta $\zeta$}: 

\medskip$\square$\,
its weight ${\sf w} ( \zeta) = a_1 + b_1 + \cdots + a_{\sf h} + b_{\sf h}$;

\medskip$\square$\,
its depth ${\sf d} ( \zeta) = 1 + b_1 + \cdots + 1 + b_{\sf h}$;
 
\medskip$\square$\,
and its height ${\sf h} ( \zeta)$.

\medskip\noindent
Also, it will be admitted throughout the present text that {\em any
polyzeta can be written either with or without pointing out the
implicit $1^0$ between any two subsequent $a_i \geqslant 2$, 
$a_{ i+1} \geqslant 2$ in it}.

\section{Dualities}
\label{dualities}
\HEAD{Dualities}{
Jo\"el {\sc Merker}, Universit\'e Paris-Sud Orsay, France}

Before proceeding further, it is appropriate, at this point, 
to speak of the important and very useful relations of duality
(to be studied and exploited in a forthcoming publication).
Here are a few examples:
\[
\aligned
\zeta(3)
&
=
\zeta(2,1)
\ \ \ \ \ \ \ \ \ \ \ \ \ \ \ \ \ \ \ \ \ \
\text{\rm that is to say:}
\ \ \ \ \ \ \ \ \ \ \ \ \ \ \ 
\zeta(3,1^0)
=
\zeta(2,1),
\\
\zeta(6,2)
&
=
\zeta(2,2,1,1,1,1)
\ \ \ \ \ \ \ \
\text{\rm that is to say:}
\ \ \ \ \
\zeta(6,1^0,2,1^0)
=
\zeta(2,1^4).
\endaligned
\]
Sometimes, the {\sl dual} to a polyzeta $\zeta$ will be denoted using
the sign $\circ$: 
\[
\zeta^\circ
:=
\mathmotsf{dual of the polyzeta}\,\,\zeta,
\]
because
in the existing theory, the sign $*$ is already
used in order to denote the stuffle product ({\em see} below), and
because we want to reserve the prime symbol
$'$ for other technical purposes. 

In order to find the {\sl dual} of
any polyzeta:
\[
\zeta\big(a_1,1^{b_1},\dots,a_{{\sf h}-1},1^{b_{{\sf h}-1}},
a_{\sf h},1^{b_{\sf h}}\big),
\]
the recipe is simple as soon as one adopts the general notation 
of the previous section:

\medskip$\square$
read the entries backwards, starting from the last entry;

\medskip$\square$
replace any incoming $1^{b_j}$ by a new single
entry $b_j+2$ which is $\geqslant 2$;

\medskip$\square$
replace any incoming $a_i \geqslant 2$ by a new $1^{a_i-2}$ with 
a number of repetitions of $1$'s which is 
$\geqslant 0$.

\begin{DualityTheorem}
Any convergent polyzeta is equal to its 
{\small\sf\em dual} polyzeta:
\[
\boxed{
\small
\aligned
\zeta
\big(a_1,1^{b_1},\dots,a_{{\sf h}-1},1^{b_{{\sf h}-1}},
a_{\sf h},1^{b_{\sf h}}
\big)
&
=
\zeta
\big(
b_{\sf h}+2,1^{a_{\sf h}-2},
b_{{\sf h}-1}+2,1^{a_{{\sf h}-1}-2},
\dots\dots,
b_1+2,1^{a_1-2}
\big)
\\
&
=:
\Big[\zeta\big(a_1,1^{b_1},\dots,a_{\sf h},1^{b_{\sf h}}\big)\Big]^\circ,
\endaligned}
\]
for every height ${\sf h} \geqslant 1$, every $a_1, \dots, a_{\sf h}
\geqslant 2$ and every $b_1, \dots, b_{\sf h} \geqslant 0$.
Furthermore, duality is an involution:
\[
({\zeta^\circ})^\circ
=
\zeta
\]
between convergent polyzetas.
\end{DualityTheorem}
 
The reader is referred to~\cite{ Zudilin-2003}, p.11, for an 
elementary
proof of these relations of duality which will not be copied here; it
comes from a simple and natural change of variable in a Chen integral
representing the polyzeta.

In weight ${\sf w} = 2$, there is only one duality: $\zeta ( 2) =
\zeta ( 2)$, but it gives no information. Polyzetas equal to their
dual will be said to be {\sl self-dual}, they exist only when the
weight is even, and they can well be counted, as a forthcoming work
will show.

In weight ${\sf w} = 3$, there is only one duality:
\[
\zeta(3)^\circ
=
\zeta(2,1),
\]
namely Euler's relation. 

In weight ${\sf w} = 4$, there are two dualities:
\[
\zeta(4)^\circ
=
\zeta(2,1,1)
\ \ \ \ \ \ \ \ \ \ \ \
\text{\rm and}
\ \ \ \ \ \ \ \ \ \ \ \
\zeta(2,2)^\circ
=
\zeta(2,2),
\]
the latter being self-dual.

In weight ${\sf w} = 5$, there are four dualities:
\[
\zeta(5)^\circ
=
\zeta(2,1,1,1),
\ \ \ \ \ \ \ \
\zeta(4,1)^\circ
=
\zeta(3,1,1),
\ \ \ \ \ \ \ \
\zeta(3,2)^\circ
=
\zeta(2,2,1),
\ \ \ \ \ \ \ \
\zeta(2,3)^\circ
=
\zeta(2,1,2).
\]

In weight ${\sf w} = 6$, there are six dualities:
\[
\aligned
\zeta(6)^\circ
&
=
\zeta(2,1,1,1,1),
\ \ \ \ \ \ \ \ \ \ \ 
\zeta(5,1)^\circ
=
\zeta(3,1,1,1),
\ \ \ \ \ \ \ \ \ \ \ \ \ \ \ 
\zeta(4,2)^\circ
=
\zeta(2,2,1,1),
\\
\zeta(3,3)^\circ
&
=
\zeta(2,1,2,1),
\ \ \ \ \ \ \ \ \ \ \ \ \ \ \,
\zeta(2,4)^\circ
=
\zeta(2,1,1,2),
\ \ \ \ \ \ \ \ \ \ \
\zeta(3,1,2)^\circ
=
\zeta(2,3,1).
\endaligned
\]

\begin{Lemma}
Weight and height remain unchanged through duality:
\[
{\sf w}\big(\zeta^\circ)
=
{\sf w}(\zeta)
\ \ \ \ \ \ \ \ \ \ \ 
\text{\rm and}
\ \ \ \ \ \ \ \ \ \ \ 
{\sf h}\big(\zeta^\circ\big)
=
{\sf h}(\zeta),
\]
while depth is symmetrized across the medium weight $\frac{\sf w}{2}$:
\[
{\sf d}\big(\zeta^\circ\big)
=
{\sf w}(\zeta)
-
{\sf d}(\zeta).
\]
\end{Lemma}

\proof
We use the general writing $\zeta \big( a_1, 1^{ b_1}, \dots,
a_{\sf h}, 1^{b_{\sf h}}\big)$. Firstly, invariance
of weight amounts to the trivial arithmetical identity:
\[
a_1+b_1+\cdots+a_{\sf h}+b_{\sf h}
=
b_{\sf h}+2+a_{\sf h}-2
+\cdots+
b_1+2+a_1-2.
\]
Secondly, invariance of height is immediately visible
in the duality theorem. 
Thirdly and lastly, 
one has by definition:
\[
\aligned
{\sf d}(\zeta)
=
1+b_1+\cdots+1+b_{\sf h}
\ \ \ \ \ \ \ \ \ \
\text{\rm and}
\ \ \ \ \ \ \ \ \ \
{\sf w}(\zeta)
=
a_1+b_1+\cdots+a_{\sf h}+b_{\sf h},
\endaligned
\]
while on the dual side:
\[
\aligned
d\big(\zeta^\circ\big)
=
1+a_{{\sf h}-2}
+\cdots+
1+a_1-2
&
=
a_1+b_1-1-b_1
+\cdots+
a_{\sf h}+b_{\sf h}-1-b_{\sf h}
\\
&
=
{\sf w}(\zeta)-{\sf d}(\zeta),
\endaligned
\]
which completes this quite elementary proof.
\endproof

\section{Detachable and Flexible Notations}
\label{detachable-notations}
\HEAD{Detachable and Flexible Notations}{
Jo\"el {\sc Merker}, Universit\'e Paris-Sud Orsay, France}

The previous section showed that the notation:
\[
\zeta\big(a_1,1^{b_1},\dots,a_{\sf h},1^{{\sf b}_h}\big)
\]
is, in particular, well adapted to the expression of the dualities.
But in several circumstances, it is advisable, and useful, to also
make visible the weight ${\sf w}$, the depth ${\sf d}$ and the height
${\sf h}$ of any written polyzeta, since only the height is visible in
such a writing. This is why we shall regularly employ the more
precise notation:
\[
\zeta_{[{\sf w},{\sf d},{\sf h}]}
\big(a_1,1^{b_1},\dots,a_{\sf h},1^{{\sf b}_h}\big),
\]
in which ${\sf w}$, ${\sf d}$, ${\sf h}$ appear as lower case indices,
just before the entries of the polyzeta. Sometimes, when it is
understood from the context what the weight is, the mention of ${\sf
w}$ will be allowed to be dropped:
\[
\zeta_{[{\sf d},{\sf h}]}
\big(a_1,1^{b_1},\dots,a_{\sf h},1^{{\sf b}_h}\big).
\]
What matters most is that the notations be {\em flexible}, with the lower
case indices ${}_{[{\sf w}, {\sf d}, {\sf h}]}$ or ${}_{[ {\sf d},
{\sf h}]}$ being detachable, either present or absent.

In concrete examples for instance, it is rarely advisable to specify
weights, depths and heights, because they may be reconstituted by a
glance. In fact, writing $\zeta_{ [ 8, 2, 2]} ( 6, 2)$ for $\zeta ( 6,
2)$ just makes the reading harder.

Beyond this, we shall also even adopt the convention of {\em sometimes not
writing the letter $\zeta$ at all}, writing for instance Euler's
relation simply under the form:
\[
(3)
=
(2,1).
\]
Since only polyzetas will be dealt with in this
paper, no risk of ambiguity, of confusion, or of misunderstanding exists.
We would like to mention that in
all our hand manuscripts, we found it more economical and expeditious
not to write the $\zeta$ letters, exactly as one would do in computer
programming.

\begin{Summaryaboutnotations}
{\em In various places, four flexible, detachable, interchangeable
notations will be employed to denote a general, arbitrary polyzeta:}
\[
\boxed{
\aligned
&
\big(a_1,1^{b_1},\dots,a_{\sf h},1^{b_{\sf h}}\big),
\\
&
\zeta\big(a_1,1^{b_1},\dots,a_{\sf h},1^{b_{\sf h}}\big),
\\
&
\zeta_{[{\sf d},{\sf h}]}
\big(a_1,1^{b_1},\dots,a_{\sf h},1^{b_{\sf h}}\big),
\\
&
\zeta_{[{\sf w},{\sf d},{\sf h}]}
\big(a_1,1^{b_1},\dots,a_{\sf h},1^{b_{\sf h}}\big),
\endaligned}
\]
{\em and these four notations will be considered
to be completely equivalent.}
\end{Summaryaboutnotations}

\begin{Lemma}
For any weight ${\sf w} \geqslant 2$, the total number of convergent
polyzetas having weight ${\sf w}$ is equal to:
\[
{\bf n}({\sf w})
:=
2^{{\sf w}-2}.
\]
\end{Lemma}

\proof
As we saw implicitly above when listing the dualities in weights ${\sf
w} = 2, 3, 4, 5, 6$, this is already known to be true for small
${\sf w}$.

By induction on ${\sf w}$, suppose now that the number of convergent
polyzetas $\zeta( s_1, \dots, s_{\sf d})$ of weight ${\sf w}$ is
indeed equal to $2^{{\sf w} - 2}$. Then we claim that any convergent
polyzeta of weight ${\sf w} + 1$ may be obtained from:
\[
\zeta(s_1,\dots,s_{\sf d})
\]
in exactly two non-overlapping ways:

\medskip\noindent$\square$
either by adding $+1$ to the last entry, getting:
\[
\zeta(s_1,\dots,s_{\sf d}+1);
\]

\medskip\noindent$\square$
or by concatenating $1$ after the last entry, getting:
\[
\zeta(s_1,\dots,s_{\sf d},1).
\]

\medskip\noindent
One easily convinces oneself of this fact by observing that any
polyzeta of weight ${\sf w} + 1$ either has $1$ as its last entry, or
has a last entry which is $\geqslant 2$. Consequently, the number of
convergent polyzetas of weight ${\sf w} + 1$ equals twice that of
weight ${\sf w}$, namely $2\cdot 2^{{\sf w} - 2} = 2^{{\sf w} + 1 -
2}$, completing the induction.
\endproof

Since the number of polyzetas grows exponentially with
their weights, one must set up a fine ordering between all of them.

\section{Total Ordering $\prec$ Between Polyzetas of Fixed Weight}
\label{total-ordering-polyzetas}
\HEAD{Total Ordering $\prec$ Between Polyzetas of Fixed Weight}{
Jo\"el {\sc Merker}, Universit\'e Paris-Sud Orsay, France}

To launch the main computations of the paper, the next initial goal is
to set up an appropriate total ordering between all polyzetas having
equal weight. Several mental speculations that are uneasy to reproduce
imposed the specific choice presented here.

In accordance to experimental and speculative conjectures about the
diophantine properties of polyzetas which conducted
experts to believe that no
algebraic relation exists between polyzetas of different weights, only
polyzetas of equal weight will be compared.

Thus, let us take any two {\em distinct} polyzetas of equal weight
${\sf w}$:
\[
\aligned
&
\zeta_{[{\sf w},{\sf d'},{\sf h'}]}
\big(
a_1',1^{b_1'},\dots,a_{{\sf h'}}',1^{b_{\sf h'}'}\big)
\\
\text{\rm and:}
\ \ \ \ \ \ \ \ \ \ \
&
\zeta_{[{\sf w},{\sf d''},{\sf h''}]}
\big(
a_1'',1^{b_1''},\dots,a_{{\sf h''}}'',1^{b_{\sf h''}''}\big),
\endaligned
\]
but of possibly unequal depths ${\sf d}'$, ${\sf d''}$ and heights
${\sf h'}$, ${\sf h''}$, with entries that are arbitrary.

First of all, increasing depth will be the dominant criterion of
ordering, and we declare that
\[
\zeta_{[{\sf w},{\sf d'},{\sf h'}]}
\big(
a_1',1^{b_1'},\dots,a_{{\sf h'}}',1^{b_{\sf h'}'}\big)
\prec
\zeta_{[{\sf w},{\sf d''},{\sf h''}]}
\big(
a_1'',1^{b_1''},\dots,a_{{\sf h''}}'',1^{b_{\sf h''}''}\big),
\]
(strictly)
for any entries in both sides whenever:
\[
{\sf d}'
<
{\sf d}''.
\]

Next, when the two depths are equal while heights are unequal, namely when:
\[
{\sf d'}
=
{\sf d''}
=:
{\sf d}
\ \ \ \ \ \ \ \ \ \ \ \ \ \ \
\text{\rm but}
\ \ \ \ \ \ \ \ \ \ \ \ \ \ \
{\sf h'}<{\sf h''},
\]
the height will be the sub-dominant criterion, and 
we also declare in this case that:
\[
\zeta_{[{\sf w},{\sf d},{\sf h'}]}
\big(
a_1',1^{b_1'},\dots,a_{{\sf h'}}',1^{b_{\sf h'}'}\big)
\prec
\zeta_{[{\sf w},{\sf d},{\sf h''}]}
\big(
a_1'',1^{b_1''},\dots,a_{{\sf h''}}'',1^{b_{\sf h''}''}\big),
\]
for any entries in both sides.

When the two depths and the two heights are equal: 
\[
{\sf d'}
=
{\sf d''}
=:
{\sf d}
\ \ \ \ \ \ \ \ \ \ \ \ \ \ \
\text{\rm and}
\ \ \ \ \ \ \ \ \ \ \ \ \ \ \
{\sf h'}
=
{\sf h''}
=:
{\sf h},
\]
the third criterion for ordering will be the place
where the $1$ lie. Thus, consider two distinct
polyzetas of same weight:
\[
\zeta_{[{\sf w},{\sf d},{\sf h}]}
\big(a_1',1^{b_1'},\dots,a_{\sf h}',1^{b_{\sf h}'}\big)
\ \ \ \ \ \ \ \ \ \ \ \ \ \ \
\text{\rm and}
\ \ \ \ \ \ \ \ \ \ \ \ \ \ \
\zeta_{[{\sf w},{\sf d},{\sf h}]}
\big(a_1'',1^{b_1''},\dots,a_{\sf h}'',1^{b_{\sf h}''}\big)
\]
having equal height ${\sf h}$ and equal depth:
\[
a_1'+b_1'+\cdots+a_{\sf h}'+b_{\sf h}'
=
a_1''+b_1''+\cdots+a_{\sf h}''+b_{\sf h}''
=:
{\sf d}.
\]
It follows that the number of $1$ present in each of them is the same,
for this number is, as known, equal to the difference between the
(common) height and the (common) depth:
\[
b_1'+\cdots+b_{\sf h}'
=
b_1''+\cdots+b_{\sf h}''
=
{\sf d}-{\sf h}.
\]
Consequently, in any such two polyzetas
to be compared and ordered, the places of the $1$ are encoded
by two multiindices:
\[
\N^{\sf h}
\ni
\big(b_1',\dots,b_{\sf h}'\big)
\ \ \ \ \ \ \ \ \ \ \ \ \ \ \
\text{\rm and}
\ \ \ \ \ \ \ \ \ \ \ \ \ \ \
\big(b_1'',\dots,b_{\sf h}''\big)
\in 
\N^{\sf h}
\]
of equal length:
\[
\vert b'\vert
=
b_1'+\cdots+b_{\sf h}'
=
b_1''+\cdots+b_{\sf h}''
=
\vert b_{\sf h}''\vert.
\]
Then by definition, whenever these two multiindices are
distinct:
\[
b'\neq b'',
\]
we declare that:
\[
\zeta_{[{\sf w},{\sf d},{\sf h}]}
\big(a_1',1^{b_1'},\dots,a_{\sf h}',1^{b_{\sf h}'}\big)
\prec
\zeta_{[{\sf w},{\sf d},{\sf h}]}
\big(a_1'',1^{b_1''},\dots,a_{\sf h}'',1^{{\sf h}''}\big)
\]
if the first multiindex encoding the $1$ is {\em larger}
than the second one, with respect to the
{\em reverse lexicographic ordering}:
\[
\big(b_1',\dots,b_{\sf h}'\big)
>_{\smallmathmotsf{revlex}}
\big(b_1'',\dots,b_{\sf h}''\big).
\]
We recall that the so-called {\sl reverse lexicographic ordering} 
is the standard lexicographic ordering, though read backwards
from the last term, namely for any two distinct
multiindices $b' \neq b''$ in $\N^h$, one defines:
\[
\big(b_1',\dots,b_{\sf h}'\big)
>_{\smallmathmotsf{revlex}}
\big(b_1'',\dots,b_{\sf h}''\big)
\]
in all of the following mutually exclusive and exhaustive circumstances:

\medskip\noindent$\square$
when $b_{\sf h}'>b_{\sf h}''$;

\medskip\noindent$\square$
or when $b_{\sf h'} = b_{\sf h}''$ but $b_{{\sf h}-1}' > b_{{\sf h}-1}''$;

\medskip\noindent$\square$
$\cdots\cdots\cdots\cdots\cdots\cdots\cdots\cdots\cdots\cdots\cdots\cdots$

\medskip\noindent$\square$
or when $b_{\sf h'} = b_{\sf h}''$,
$b_{{\sf h}-1}' > b_{{\sf h}-1}''$, \dots, $b_2' = b_2''$ but
$b_1' > b_1''$.

\medskip
For instance, convergent 
polyzetas of depth ${\sf d} = 5$ and height ${\sf h} = 3$
are ordered in the following schematic way
(entries are drawn vertically without parentheses):
\[
\begin{array}{c}
a
\\
a
\\
a
\\
1
\\
1
\end{array}
\prec
\begin{array}{c}
a
\\
a
\\
1
\\
a
\\
1
\end{array}
\prec
\begin{array}{c}
a
\\
1
\\
a
\\
a
\\
1
\end{array}
\prec
\begin{array}{c}
a
\\
a
\\
1
\\
1
\\
a
\end{array}
\prec
\begin{array}{c}
a
\\
1
\\
a
\\
1
\\
a
\end{array}
\prec
\begin{array}{c}
a
\\
1
\\
1
\\
a
\\
a,
\end{array}
\]
whatever the entries $a \geqslant 2$ are.

It therefore only remains to compare and to order any two
distinct convergent polyzetas having same weight (as agreed), 
same depth, same height and in which the places
of the $1$ are exactly the same, so that:
\[
b_1'=b_1''=:b_1,\ \ \ \ \
b_2'=b_2''=:b_2,\ \ \ \ \
\cdots\cdots,\ \ \ \ \
b_{\sf h}'=b_{\sf h}''=:b_{\sf h},
\]
that is to say, to compare any two:
\[
\zeta_{[{\sf w},{\sf d},{\sf h}]}
\big(a_1',1^{b_1},\dots,a_{\sf h}',1^{b_{\sf h}}\big)
\ \ \ \ \ \ \ \ \ \ \ \ \ \ \
\text{\rm and}
\ \ \ \ \ \ \ \ \ \ \ \ \ \ \
\zeta_{[{\sf w},{\sf d},{\sf h}]}
\big(a_1'',1^{b_1},\dots,a_{\sf h}'',1^{b_{\sf h}}\big).
\]
In such a circumstance, one observes that equality
of weights:
\[
a_1'+b_1+\cdots+a_{\sf h}'+b_{\sf h}
=
a_1''+b_1+\cdots+a_{\sf h}''+b_{\sf h}
\]
and equality of the number of $1$ immediately gives that
the two multiindices:
\[
\N^{\sf h}
\ni
\big(a_1',\dots,a_{\sf h}'\big)
\ \ \ \ \ \ \ \ \ \ \ \ \ \ \
\text{\rm and}
\ \ \ \ \ \ \ \ \ \ \ \ \ \ \
\big(a_1'',\dots,a_{\sf h}''\big)
\in
\N^{\sf h}
\]
which encode all the entries that are $\geqslant 2$ in the
two polyzetas have equal length:
\[
\vert a'\vert
=
a_1'+\cdots+a_{\sf h}'
=
a_1''+\cdots+a_{\sf h}''
=
\vert a''\vert.
\]
Most simply then, to compare such two distinct
multiindices $a'$ and $a''$ with $\vert a' \vert = 
\vert a'' \vert$, 
we will declare that:
\[
\zeta_{[{\sf w},{\sf d},{\sf h}]}
\big(a_1',1^{b_1},\dots,a_{\sf h}',1^{b_{\sf h}}\big)
\prec
\zeta_{[{\sf w},{\sf d},{\sf h}]}
\big(a_1'',1^{b_1},\dots,a_{\sf h}'',1^{b_{\sf h}}\big),
\]
whenever the first multiindex is {\em larger} than the
second one:
\[
\big(a_1',\dots,a_{\sf h}'\big)
\prec_{\smallmathmotsf{lex}}
\big(a_1'',\dots,a_{\sf h}''\big),
\]
with respect to the standard lexicographic ordering.

For instance, the ordering of all polyzetas of weight $10$, of depth
$6$, of height $3$, and of $1$-type $\zeta \big( a, 1, a,
1, 1, a \big)$ is the following\,\,---\,\,again, 
entries are drawn vertically without parentheses\,\,---\,:
\[
\begin{array}{c}
4
\\
1
\\
2
\\
1
\\
1
\\
2
\end{array}
\prec
\begin{array}{c}
3
\\
1
\\
3
\\
1
\\
1
\\
2
\end{array}
\prec
\begin{array}{c}
3
\\
1
\\
2
\\
1
\\
1
\\
3
\end{array}
\prec
\begin{array}{c}
2
\\
1
\\
4
\\
1
\\
1
\\
2
\end{array}
\prec
\begin{array}{c}
2
\\
1
\\
3
\\
1
\\
1
\\
3
\end{array}
\prec
\begin{array}{c}
2
\\
1
\\
2
\\
1
\\
1
\\
4.
\end{array}
\]

\subsection{Total ordering in small weights ${\sf w} = 3, 4, 5, 6$}
To illustrate this, let us list convergent polyzetas in accordance
to this ordering (entries are
drawn vertically without parentheses):
\[
\footnotesize
\aligned
\begin{array}{ccc}
& & {\bf 2} 
\\
\underline{{\sf w}=3:}\ \
& {\bf 3} & {\bf 1}
\end{array}
\ \ \ \ \ \ \ \ \ \ \ \ \ \ \ \ \ \ \ \ \ \ 
\begin{array}{ccccc}
\\
& & & & {\bf 2}
\\
&  & {\bf 3} & {\bf 2} & {\bf 1}
\\
\underline{{\sf w}=4:}\ \
& {\bf 4} & {\bf 1} & {\bf 2} & {\bf 1}
\end{array}
\endaligned
\]
\[
\footnotesize
\aligned
\begin{array}{ccccccccc}
& & & &
& & & & {\bf 2}
\\
& & & &
& {\bf 3} & {\bf 2} & {\bf 2} & {\bf 1}
\\
& & {\bf 4} & {\bf 3} & {\bf 2}
& {\bf 1} & {\bf 2} & {\bf 1} & {\bf 1}
\\
\underline{{\sf w}=5:}\ \
& {\bf 5} & {\bf 1} & {\bf 2} & {\bf 3}
& {\bf 1} & {\bf 1} & {\bf 2} & {\bf 1}
\end{array}
\endaligned
\]
\[
\footnotesize
\aligned
\begin{array}{ccccccccccccccccc}
& & & & 
& & & &
& & & &
& & & & {\bf 2}
\\
& & & & 
& & & &
& & & & {\bf 2}
& {\bf 2} & {\bf 2} & {\bf 2} & {\bf 1}
\\
& & & & 
& & {\bf 4} & {\bf 3} & {\bf 2}
& {\bf 3} & {\bf 2} & {\bf 2} & {\bf 1}
& {\bf 2} & {\bf 1} & {\bf 1} & {\bf 1}
\\
& & {\bf 5} & {\bf 4} & {\bf 3}
& {\bf 2} & {\bf 1} & {\bf 2} & {\bf 3}
& {\bf 1} & {\bf 1} & {\bf 2} & {\bf 1}
& {\bf 1} & {\bf 2} & {\bf 1} & {\bf 1}
\\
\underline{{\sf w}=6:}\ \
& {\bf 6} & {\bf 1} & {\bf 2} & {\bf 3}
& {\bf 4} & {\bf 1} & {\bf 1} & {\bf 1}
& {\bf 2} & {\bf 3} & {\bf 2} & {\bf 1}
& {\bf 1} & {\bf 1} & {\bf 2} & {\bf 1}
\end{array}
\endaligned
\]

\subsection{Summarized presentation of the total ordering between
convergent polyzetas of fixed weight} 
{\em By
definition, two arbitrary distinct
convergent polyzetas of equal weight are
smaller one to another:}
\[
\zeta_{[{\sf w},{\sf d'},{\sf h'}]}
\big(a_1',1^{b_1'},\dots,a_{\sf h'}',1^{b_{\sf h'}'}\big)
\prec
\zeta_{[{\sf w},{\sf d'},{\sf h'}]}
\big(a_1'',1^{b_1''},\dots,a_{\sf h''}'',1^{b_{\sf h''}''}\big)
\]
{\em whenever firstly, either the depth of the first is smaller
than that of the second:}
\[
1+b_1'+\cdots+1+b_{\sf h'}'
=
{\sf d'}
<
{\sf d''}
=
1+b_1''+\cdots+1+b_{\sf h''}'',
\]
{\em or secondly, whenever their depths coincide, but the height
of the first is smaller than that of the second:}
\[
{\sf d'}
=
{\sf d''}
\ \ \ \ \ \ \ \ \ \ \ \ \ \ 
\text{\rm but}
\ \ \ \ \ \ \ \ \ \ \ \ \ \ 
{\sf h'}
<
{\sf h''},
\]
{\em or thirdly, whenever
both their depths and their heights ${\sf h'} = {\sf h''}
=: {\sf h}$ coincide, but their $b$-integers encoding
their places of the $1$ satisfy:}
\[
\big(b_1',\dots,b_{\sf h}'\big)
>_{\smallmathmotsf{revlex}}
\big(b_1'',\dots,b_{\sf h}''\big),
\]
{\em or fourthly and lastly, whenever depths, heights and places
of the $1$ coincide, but:}
\[
\big(a_1',\dots,a_{\sf h}'\big)
>_{\smallmathmotsf{lex}}
\big(a_1'',\dots,a_{\sf h}''\big).
\]

\section{Intermediate Countings}
\label{intermediate-countings}
\HEAD{Intermediate Countings}{
Jo\"el {\sc Merker}, Universit\'e Paris-Sud Orsay, France}

We remember that the total number of convergent
polyzetas of weight ${\sf w}$ is equal to: 
\[
{\bf n}({\sf w})
:=
2^{{\sf w}-2}.
\]
More precise intermediate countings can now be provided.

\begin{Lemma}
The total number of convergent polyzetas: 
\[
\zeta_{[{\sf w},{\sf d},{\sf h}]}
\big(a_1,1^{b_1},\dots,a_{\sf h},1^{b_{\sf h}}\big)
\]
of weight ${\sf w}$, of depth ${\sf d}$, of height ${\sf h}$ and whose
entries $1$ are placed according to a fixed multiindex $(b_1, \dots,
b_{\sf h}) \in \N^{\sf h}$ is equal to:
\[
{\bf n}\big({\sf w},{\sf d},{\sf h},(b_1,\dots,b_{\sf h})\big)
:=
\binom{{\sf w}-{\sf d}-1}{{\sf h}-1}.
\]
\end{Lemma}

Observe that this number does not depend on the assignment of places
to the $1$, {\em i.e.} ${\bf n} \big( {\sf w}, {\sf d}, {\sf h},
(b_1, \dots, b_{\sf h}) \big)$ 
is\,\,---\,\,visibly\,\,---\,\,independent of $(b_1,
\dots, b_{\sf h}) \in \N^{\sf h}$.

\proof
Indeed, in such convergent polyzetas, only the integers $a_i \geqslant
2$ can still vary, but with the only constraint that their sum equals:
\[
\aligned
a_1+\cdots+a_{\sf h}
&
=
{\sf w}-b_1-\cdots-b_{\sf h}
\\
&
=
{\sf w}+{\sf h}-{\sf d}.
\endaligned
\]
But it is either elementary to check or already known that:
\[
\mathmotsf{Card}
\big\{
\big(\widetilde{a}_1,\dots,\widetilde{a}_{\sf h}\big)
\in
\N^{\sf h}\colon\
\widetilde{a}_1+\cdots+\widetilde{a}_{\sf h}
=
{\sf m}
\big\}
=
\binom{{\sf m}+{\sf h}-1}{{\sf h}-1}.
\]
Since the $a_i$ must be $\geqslant 2$, in order to come
to some $\widetilde{ a_i}$ that are $\geqslant 0$, one has
to subtract $2$ each time:
\[
\underbrace{a_1-2}_{=:\,\widetilde{a}_1}
+\cdots+
\underbrace{a_{\sf h}-2}_{=:\,\widetilde{a}_{\sf h}}
=
{\sf w}+{\sf h}-{\sf d}-2\,{\sf h}
=
\underbrace{{\sf w}-{\sf d}-{\sf h}}_{=:\,{\sf m}}
\]
so that substituting ${\sf m} := {\sf w} + {\sf h} - {\sf d}$ in the
binomial completes
the proof.
\endproof

\begin{Lemma}
The total number of convergent polyzetas
$\zeta_{[ {\sf w}, {\sf d}, {\sf h}]}$ of weight ${\sf w}$,
of depth ${\sf d}$ and of height ${\sf h}$ is equal to:
\[
{\bf n}({\sf w},{\sf d},{\sf h})
:=
\binom{{\sf d}-1}{{\sf d}-{\sf h}}\,
\binom{{\sf w}-{\sf d}-1}{{\sf h}-1}.
\]
\end{Lemma}

\proof 
Of course, the first entry of a convergent polyzeta
$\zeta_{[ {\sf w}, {\sf d}, {\sf h}]}$ must be $\geqslant
2$, namely it must always be an $a_1 \geqslant 2$. But then,
the ${\sf d} - {\sf h}$ existing $1$ can take any place in the ${\sf
d} - 1$ entries of a convergent $\zeta_{[ {\sf w}, {\sf d}, {\sf h}]}$
which lie after $a_1$. As is known, the number of places that the
$1$ can be so given in a convergent $\zeta_{[ {\sf w}, {\sf d}, {\sf h}]}$
is equal to the plain binomial:
\[
\binom{{\sf d}-1}{{\sf d}-{\sf h}}.
\] 
Once the places of the $1$ are fixed, the number $\binom{{\sf w} -
{\sf d} - 1}{{\sf h} - 1}$ of variations of the ${\sf a}_i$
was counted just above, and lastly, the total number of possibilities is
multiplicative.
\endproof

\begin{Corollary}
The total number of convergent polyzetas of weight ${\sf w}$ and of
depth ${\sf d}$ is equal to:
\[
{\bf n}({\sf w},{\sf d})
:=
\sum_{{\sf h}=1}^{{\sf h}=2}\,
\binom{{\sf d}-1}{{\sf d}-{\sf h}}\,
\binom{{\sf w}-{\sf d}-1}{{\sf h}-1}.
\qed
\]
\end{Corollary}

For fun and for later use, let us list these dimensions ${\bf n} (
{\sf w}, {\sf d})$ in weights $6$, $7$, $8$, $9$, $10$.
\[
\underline{{\sf w}=6:}
\ \ \ \ \ \ \ \ \ \ \ \ \ 
\left.
\aligned
1
&
=
{\bf 1}
\\
1\cdot 1+1\cdot 3
&
=
{\bf 4}
\\
1\cdot 1+2\cdot 2+1\cdot 1
&
=
{\bf 6}
\\
1\cdot 1+1\cdot 3
&
=
{\bf 4}
\\
1\cdot 1
&
=
{\bf 1}
\endaligned
\ \ \ \ 
\right]
\ \ \ 
{\bf 16}.
\]
\[
\underline{{\sf w}=7:}
\ \ \ \ \ \ \ \ \ \ \ \ \ 
\left.
\aligned
1
&
=
{\bf 1}
\\
1\cdot 1+1\cdot 4
&
=
{\bf 5}
\\
1\cdot 1+2\cdot 3+1\cdot 3
&
=
{\bf 10}
\\
1\cdot 1+2\cdot 3+1\cdot 3
&
=
{\bf 10}
\\
1\cdot 1+1\cdot 4
&
=
{\bf 5}
\\
1\cdot 1
&
=
{\bf 1}
\endaligned
\ \ \ \ 
\right]
\ \ \ 
{\bf 32}.
\]
\[
\underline{{\sf w}=8:}
\ \ \ \ \ \ \ \ \ \ \ \ \ 
\left.
\aligned
1
&
=
{\bf 1}
\\
1\cdot 1+1\cdot 5
&
=
{\bf 6}
\\
1\cdot 1+2\cdot 4+1\cdot 6
&
=
{\bf 15}
\\
1\cdot 1+3\cdot 3+3\cdot 3+1\cdot 1
&
=
{\bf 20}
\\
1\cdot 1+2\cdot 4+1\cdot 6
&
=
{\bf 15}
\\
1\cdot 1+1\cdot 5
&
=
{\bf 6}
\\
1
&
=
{\bf 1}
\endaligned
\ \ \ \ 
\right]
\ \ \ 
{\bf 64}.
\]
\[
\underline{{\sf w}=9:}
\ \ \ \ \ \ \ \ \ \ \ \ \ 
\left.
\aligned
1
&
=
{\bf 1}
\\
1\cdot 1+1\cdot 6
&
=
{\bf 7}
\\
1\cdot 1+2\cdot 5+1\cdot 10
&
=
{\bf 21}
\\
1\cdot 1+3\cdot 4+3\cdot 6+1\cdot 4
&
=
{\bf 35}
\\
1\cdot 1+3\cdot 4+3\cdot 6+1\cdot 4
&
=
{\bf 35}
\\
1\cdot 1+2\cdot 5+1\cdot 10
&
=
{\bf 21}
\\
1\cdot 1+1\cdot 6
&
=
{\bf 7}
\\
1
&
=
{\bf 1}
\endaligned
\ \ \ \ 
\right]
\ \ \ 
{\bf 128}.
\]
\[
\underline{{\sf w}=10:}
\ \ \ \ \ \ \ \ \ \ \ \ \ 
\left.
\aligned
1
&
=
{\bf 1}
\\
1\cdot 1+1\cdot 7
&
=
{\bf 8}
\\
1\cdot 1+2\cdot 6+1\cdot 15
&
=
{\bf 28}
\\
1\cdot 1+3\cdot 5+3\cdot 10+1\cdot 10
&
=
{\bf 56}
\\
1\cdot 1+4\cdot 4+6\cdot 6+4\cdot 4+1\cdot 1
&
=
{\bf 70}
\\
1\cdot 1+3\cdot 5+3\cdot 10+1\cdot 10
&
=
{\bf 56}
\\
1\cdot 1+2\cdot 6+1\cdot 15
&
=
{\bf 28}
\\
1\cdot 1+1\cdot 7
&
=
{\bf 8}
\\
1
&
=
{\bf 1}
\endaligned
\ \ \ \ 
\right]
\ \ \ 
{\bf 256}.
\]

\section{Shuffle Minus Stuffle}
\label{shuffle-stuffle}
\HEAD{Shuffle Minus Stuffle}{
Jo\"el {\sc Merker}, Universit\'e Paris-Sud Orsay, France}

The result of the usual multiplication between any two polyzetas is
known to be re-expressible in {\em two} ways, either as a polyzeta
whose entries come from the combinatorial $\ast$-stuffle product
between the two concerned entries, or as a polyzeta whose entries come
from the combinatorial $\shuffle$-shuffle product between the two
concernend entries:
\[
\aligned
&
\underbrace{\zeta
\big(\alpha_1,1^{\beta_1},\dots,\alpha_{\sf g},1^{b_{\sf g}}\big)
\cdot
\zeta
\big(a_1,1^{b_1},\dots,a_{\sf h},1^{b_{\sf h}}\big)}_{
\smallmathmotsf{usual multiplication}}
=
\\
&
\ \ \ \ \ \ \ \ \ \ \ \ \ \ \ \ \ \ \ \ \ \ \ \ \ \ \ \ \ \ \ \ \ \
\ \ \ \ \ \ \ \ \ \ \ \ \ \ \ \ \ \ \ 
=
\zeta
\Big(
\underbrace{
\big(\alpha_1,1^{\beta_1},\dots,\alpha_{\sf g},1^{b_{\sf g}}\big)
\ast
\big(a_1,1^{b_1},\dots,a_{\sf h},1^{b_{\sf h}}\big)}_{
\smallmathmotsf{combinatorial s\underline{\bf t}uffle product between the two
entries}}
\Big)
\\
&
\ \ \ \ \ \ \ \ \ \ \ \ \ \ \ \ \ \ \ \ \ \ \ \ \ \ \ \ \ \ \ \ \ \
\ \ \ \ \ \ \ \ \ \ \ \ \ \ \ \ \ \ \ 
=
\zeta
\Big(
\underbrace{
\big(
0^{\alpha_1-1}1\,1^{\beta_1}\cdots 0^{\alpha_{\sf g}-1}1\,
1^{\beta_{\sf g}}
\big)
\shuffle
\big(
0^{a_1-1}1\,1^{b_1}\cdots 0^{a_{\sf h}-1}1\,1^{b_{\sf h}}
\big)
}_{
\smallmathmotsf{combinatorial s\underline{\bf h}uffle product between the two
entries}}
\Big).
\endaligned
\]
Here, we use the simplest binary symbols 
$0$ and $1$\,\,---\,\,instead of $x_0$ and $x_1$, or
instead of $x$ and $y$ as is usually done\,\,---\,\,to translate
any polyzeta in its second encoding: 
\[
\boxed{
\underbrace{
\big(
a_1,1^{b_1},\dots,a_{\sf h},1^{b_{\sf h}}
\big)}_{
\smallmathmotsf{encoding used}
\atop
\smallmathmotsf{to compute $\ast$-products}}
\overset{\longleftrightarrow}{=}
\underbrace{
\big(
0^{a_1-1}1\,1^{b_1}\cdots 0^{a_{\sf h}-1}1\,1^{b_{\sf h}}
\big)}_{
\smallmathmotsf{encoding used}
\atop
\smallmathmotsf{to compute $\shuffle$-products}}.
}
\]
Since the introductory literature is rich enough, we will not make any
further reminder here. Observably, we will not copy the classical
inductive rules for the computation of stuffle and shuffle products,
{\em because the non-closed character of these inductive rules is a
landmark of a hidden defective knowledge concerning what they really
give}.

Just let us restrict ourselves with expressing in full generality the
known-to-hold double shuffle relations that we want, in the four
special cases mentioned in the introduction, to close up, namely 
in the particular circumstances when:
\[
\big(
\alpha_1,1^{\beta_1},\dots,\alpha_{\sf g},1^{\beta_{\sf h}}
\big)
=
\left\{
\aligned
&
(1);
\\
&
(2);
\\
&
(3);
\\
&
(2,1).
\endaligned\right.
\]

\begin{Theorem}
For every two heights ${\sf h} \geqslant 1$, ${\sf g} \geqslant 1$,
every two collections of entries $\alpha_1 \geqslant 2$, \dots,
$\alpha_{\sf g} \geqslant 2$, $a_1 \geqslant 2$, \dots, $a_{\sf h}
\geqslant 2$ and for every two collections of entries $\beta_1
\geqslant 0$, \dots, $\beta_{\sf g} \geqslant 0$, $b_1 \geqslant 0$,
\dots, $b_{\sf h} \geqslant 0$, one has the so-called {\sl
double-shuffle relations}:
\[
\aligned
0
&
=
-\,
\big(
\alpha_1,1^{\beta_1},\dots,\alpha_{\sf g},1^{b_{\sf g}}
\big)
\ast
\big(
a_1,1^{b_1},\dots,a_{\sf h},1^{b_{\sf h}}
\big)
+
\\
&
\ \ \ \ \
+
\big(
0^{\alpha_1-1}1\,1^{\beta_1}\cdots 0^{\alpha_{\sf g}-1}1\,
1^{\beta_{\sf g}}
\big)
\shuffle
\big(
0^{a_1-1}1\,1^{b_1}\cdots 0^{a_{\sf h}-1}1\,1^{b_{\sf h}}
\big).
\endaligned
\]
By automatic cancelling out of two nonconvergent polyzetas both equal
to $(1, a_1, 1^{b_1}, \dots, a_{\sf h}, 1^{ b_{\sf h}} \big)$
through the subtraction, 
these double shuffle relations also hold true when:
\[
\big(\alpha_1,1^{\beta_1},\dots,\alpha_{\sf g},1^{\beta_{\sf g}}\big)
=
(1).
\]
\end{Theorem}

\section{Computing $-(1)\ast\big(a_1,1^{b_1}, \dots,a_{\sf h},
1^{b_{\sf h}}\big) + 1 \shuffle 0^{ a_1-1} 1\, 1^{b_1} \cdots 0^{
a_{\sf h}-1} 1\, 1^{b_{\sf h}}$}
\label{1-double-melange}
\HEAD{Computing $-(1)\ast\big(a_1,1^{b_1}, \dots,a_{\sf h},
1^{b_{\sf h}}\big) + 1 \shuffle 0^{ a_1-1} 1\, 1^{b_1} \cdots 0^{
a_{\sf h}-1} 1\, 1^{b_{\sf h}}$}{
Jo\"el {\sc Merker}, Universit\'e Paris-Sud Orsay, France}

\subsection{Computing firstly the stuffle product
$(1)\ast\big(a_1,1^{b_1}, \dots,a_{\sf h},1^{b_{\sf h}}\big)$}
Granted that we want a result of weight ${\sf w}$, 
we aim to compute in general the
stuffle product of $\zeta ( 1)$ with an arbitrary
convergent polyzeta of weight ${\sf w}-1$:
\[
\zeta(1)
\ast
\zeta_{[{\sf w}-1,{\sf d},{\sf h}]}
\big(a_1,1^{b_1},\dots,a_{\sf h},1^{b_{\sf h}}\big),
\]
that it to say in greater length, we aim to compute:
\[
(1)
\ast
\big(a_1,
\overbrace{1,\dots,1}^{b_1},a_2,\overbrace{1,\dots,1}^{b_2},
\bfdots,
a_{\sf h},
\overbrace{1,\dots,1}^{b_{\sf h}}\big).
\]
According to the general rule, one must either add $1$ to each one of
the entries of $\zeta_{[{\sf w}-1,{\sf d},{\sf h}]}$, or insert $1$
between each two of its entries, including the extremes. For
instance (remember that we allow dropping the letter
$\zeta$):
\[
\aligned
(\green{1})\ast(4,1,1)
&
=
(\green{5},1,1)+(4,\green{2},1)+(4,1,\green{2})
+
\ \ \ \ \ \ \ \ \ \ \ \ \ \ \ \ \ \ \ \ \ \ \ \ \ \ \ \ \ \ \ \ \ \ \ \ \ \ \
\explain{additions}
\\
&
\ \ \ \ \
+
(\green{1},4,1,1)+(4,\green{1},1,1)+(4,1,\green{1},1)+(4,1,1,\green{1})
\ \ \ \ \ \ \
\explain{insertions},
\endaligned
\]
and this gives, after gathering equal terms:
\[
(1)\ast(4,1,1)
=
(5,1,1)+(4,2,1)+(4,1,2)+(1,4,1,1)+3\,(4,1,1,1).
\]
Drawing on this example, let us perform the general computation
in the simpler case where the height ${\sf h}$ is equal to $1$:
\[
\aligned
(\green{1})
\ast
\big(a_1,
\overbrace{1,\dots,1}^{b_1}
\big)
&
=
\big(a_1+\green{1},1,\dots,1\big)
+
\big(a_1,\green{2},1,\dots,1\big)
+
\big(a_1,1,\green{2},1,\dots,1\big)
+
\\
&
\ \ \ \ \
+\cdots+
\big(a_1,1,\dots,\green{2},1\big)
+
\big(a_1,1,\dots,1,\green{2}\big)
+
\ \ \ \ \ \ \ \ \ \ \ \ \ \ \ \ 
\explain{additions}
\\
&
\ \ \ \ \ 
+
\big(\green{1},a_1,1,\dots,1\big)
+
\big(a_1,\green{1},1,\dots,1\big)
+
\big(a_1,1,\green{1},1,\dots,1\big)
+
\\
&
\ \ \ \ \
+\cdots+
\big(a_1,1,\dots,1,\green{1},1)
+
\big(a_1,1,\dots,1,\green{1}\big)
\ \ \ \ \ \ \ \ \ \ \ \ \ \ \ \ \ \ \ \ \ 
\explain{insertions}.
\endaligned
\]
In the last two lines, all terms except the first one are
equal so that after gathering:
\[
\aligned
\mathmotsf{insertions}
&
=
\big(1,a_1,
\overbrace{1,\dots,1}^{b_1}\big)
+
(b_1+1)\big(a_1,
\overbrace{1,1,\dots,1}^{b_1+1}
\big)
\\
&
=
\zeta\big(\red{1},a_1,1^{b_1}\big)
+
(b_1+1)\,\zeta\big(a_1,1^{b_1+1}\big).
\endaligned
\]
Of course, the first $\zeta\big( \red{ 1}, a_1, 1^{ b_1} \big)$ is
non-convergent, but as is known, it will simply disappear in the
double shuffle-stuffle subtraction:
\[
-\,(1)\ast\big(a_1,1^{b_1}\big)
+
(1)\shuffle\big(a_1,1^{b_1}\big),
\]
since it will also appear in $1 \shuffle \big( a_1, 1^{ b_1} \big)$.
Here and everywhere below, we intentionally write the minus sign
{\em in the first position}, 
because the terms involving
the stuffle $\ast$ are of smallest depth, so that
they should come {\em first} if one respects the chosen total order between
polyzetas of fixed weight.

Let us also rewrite the result of additions under a concise sum-like form: 
\[
\mathmotsf{additions}
=
\big(a_1+1,1^{b_1}\big)
+
\sum_{b_1'+b_1''=b_1-1
\atop
b_1'\geqslant 0,\,\,b_1''\geqslant 0}\,
\big(
a_1,\dots,1^{b_1'},2,1^{b_1''},\dots
\big).
\]
An inspection of depths and heights of appearing terms enables
us to summarize the result in the special case ${\sf h} = 1$.

\begin{Lemma}
For every $a_1 \geqslant 2$ and every $b_1 \geqslant 0$, if one
sets:
\[
a_1+b_1
=:
{\sf w}-1
\ \ \ \ \ \ \
\text{\rm and}
\ \ \ \ \ \ \
1+b_1
=:
{\sf d},
\]
then the $*$-stuffle product of $\zeta(1)$ with any $\zeta_{[ {\sf
w}-1, {\sf d}, 1]} \big( a_1, 1^{ b_1} \big)$ of weight ${\sf w}
- 1 \geqslant 2$, of height ${\sf d}$ 
and of height $1$ writes out, after complete finalization:
\[
\aligned
\zeta(1)
\ast
\zeta_{[{\sf w}-1,{\sf d},1]}
\big(a_1,1^{b_1}\big)
&
=
\zeta_{[{\sf w},{\sf d},1]}
\big(a_1+1,1^{b_1}\big)
+
\\
&
\ \ \ \ \
+
\sum_{b_1'+b_1''=b_1-1
\atop
b_1'\geqslant 0,\,\,b_1''\geqslant 0}\,
\zeta_{[{\sf w},{\sf d},2]}
\big(a_1,1^{b_1'},2,1^{b_1''}\big)
+
\\
&
\ \ \ \ \
+
\underbrace{\zeta_{[{\sf w},{\sf d}+1,1]}
\big(1,a_1,1^{b_1}\big)}_{
\smallmathmotsf{non-convergent}
\atop
\smallmathmotsf{but will disappear}}
+
(b_1+1)\,\zeta_{[{\sf w},{\sf d}+1,1]}
\big(a_1,1^{b_1+1}\big).
\endaligned
\]
\end{Lemma}

Notice that, except for the single appearing non-convergent 
polyzeta, the terms are ordered in accordance with the
total (sub-)ordering:
\[
\zeta_{[{\sf w},{\sf d'},{\sf h'}]}
\prec
\zeta_{[{\sf w},{\sf d''},{\sf h''}]}
\ \ \ \ \
\text{\rm if}\ \ \ \ \
\left\{
\aligned
&
\text{\rm either}\ \
{\sf d'}
<
{\sf d''}
\\
&
\text{\rm or}\ \ \ \ \ \ \ \
{\sf d'}
=
{\sf d''}
\ \ \ \ \
\text{\rm but}\ \ 
{\sf h'}
<
{\sf h''}.
\endaligned
\right.
\]
Now, in the general case ${\sf h} \geqslant 1$, the procedure will be
quite similar, with no obstacle of combinatorics understanding.

One one hand, additions of $1$ will concern firstly the $a_i$, $i = 1,
\dots, {\sf h}$, and secondly, the $1$ of each $1^{ b_j}$, $j = 1,
\dots, {\sf h}$, hence we will have similar concise sum-like
expressions for addition terms.

One the second hand, insertions of the additional $1$ will be achieved
in such a way that, if $1$ is inserted {\em just after} an $a_i$, one
will consider that {\em equivalently}, the $1$ is inserted just before
the first $1$ of $1^{ b_i}$ which sits just after the $a_i$ in
question, so that one gets a group of $(b_i+1)$ entries $1$, namely
$1^{ b_i+1}$ (the $1^{ b_i}$ `{\em absorbs}' the inserted $1$).
Completely similarly, if the additional $1$ is inserted {\em just
before} an $a_i$, and if $i$ is $\geqslant 2$, one will consider
that the inserted $1$ is {\em absorbed} by the group of $1$ of the
$1^{ b_{ i-1}}$ which sits just before the $a_i$ in question, so that
one gets an $1^{ b_{ i-1} + 1}$ in this way. The only exception is
when $1$ is inserted just before $a_1$, and this produces the only
non-convergent polyzeta $\big( 1, a_1, 1^{ b_1}, \dots, a_{\sf h}, 1^{
b_{\sf h}} \big)$, which will anyway disappear at the end.

Granted these explanations, we may state the general lemma without
further proof, but with a few comments just after. Later, subsequent
statements that are less direct and more complex will be established
with all details.

\begin{Lemma}
For every height ${\sf h} \geqslant 1$, every entries $a_1 \geqslant 2$,
\dots, $a_{\sf h} \geqslant 2$, every entries $b_1 \geqslant 0$, \dots,
$b_{\sf h} \geqslant 0$, if one sets:
\[
\aligned
a_1+b_1+\cdots+a_{\sf h}+b_{\sf h}
&
=:
{\sf w}-1,
\\
1+b_1+\cdots+1+b_{\sf h}
&
=:
{\sf d},
\endaligned
\]
then the $*$-stuffle product of $\zeta (1)$ with any $\zeta_{[ {\sf
w}-1, {\sf d}, {\sf h}]} \big( a_1, 1^{ b_1}, \dots, a_{\sf h},
1^{ b_{\sf h}} \big)$ of weight ${\sf w} - 1 \geqslant 2$, 
of depth ${\sf d}$ and of height ${\sf h}$
writes out, after complete finalization:
\[
\aligned
\zeta(1)
\ast
\zeta_{[{\sf w}-1,{\sf d},{\sf h}]}
\big(a_1,1^{b_1},\dots,a_{\sf h},1^{b_{\sf h}}\big)
&
=
\sum_{i=1}^{\sf h}\,
\zeta_{[{\sf w},{\sf d},{\sf h}]}
\big(\bfdots,a_i+1,\bfdots\big)
+
\\
&
\ \ \ \ \
+
\sum_{j=1\atop b_j\geqslant 1}^{\sf h}\,
\Bigg(
\sum_{b_j'+b_j''=b_j-1
\atop
b_j'\geqslant 0,\,\,b_j''\geqslant 0}\,
\zeta_{[{\sf w},{\sf d},{\sf h}+1]}
\big(
\bfdots,1^{b_j'},2,1^{b_j''},\bfdots
\big)
\Bigg)
+
\endaligned
\]
\[
\aligned
&
+
\underbrace{\zeta_{[{\sf w},{\sf d}+1,{\sf h}]}
\big(1,a_1,1^{b_1},\dots,a_{\sf h},1^{b_{\sf h}}\big)}_{
\smallmathmotsf{non-convergent}
\atop
\smallmathmotsf{but will disappear}}
+
\\
& 
+
\sum_{j=1}^{\sf h}\,
(b_j+1)\,\zeta_{[{\sf w},{\sf d}+1,{\sf h}]}
\big(
\bfdots,1^{b_j+1},\bfdots
\big).
\endaligned
\]
\end{Lemma}

Here by convention, in the right-hand side, only the terms of the
initial polyzeta $\zeta_{[ {\sf w}-1, {\sf d}, {\sf h}]}$ that are
{\em changed} are written down, so that the symbol:
\[
\bfdots
\]
means that all other entries are unchanged.

With constancy and regularity, the letter '$i$' will always be used in
relation with the entries $a_i \geqslant 2$, especially in summation
symbols. Later below, $i'$, $i''$, $i'''$ and $i_1, i_2, i_3$ will
also be used.

Similarly, the letter $j$ will always be used in link with the entries
$b_j \geqslant 0$, specially in summation symbols.

\subsection{Computing secondly the shuffle product
$( 1) \shuffle 0^{ a_1-1}1\, 1^{b_1} \cdots 0^{ a_{\sf h}-1}1\,
 1^{b_{\sf h}}$} To begin with, let us study the case of smallest
height ${\sf h} = 1$, so that one has to insert $1$ in any place
between two successive entries of the polyzeta $\big( a_1, 1^{ b_1}
\big)$, after re-coding it in terms of $0$ and $1$:

\bigskip
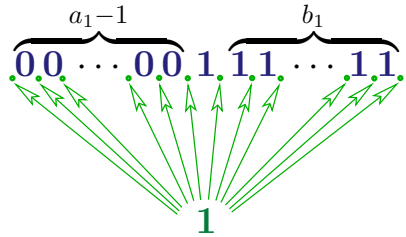
\begin{center}
\input 1-all-of.pstex_t
\end{center}

\noindent
Since its depth equals $a_1 + b_1$, 
there visibly are $a_1 + b_1 + 1$ possible insertions of $1$.
Let us set apart the case where $1$ is inserted in the
front place:
\[
\gbf{1}\,
\zsqov{a_1-1}\,
\unsqov{b_1},
\]
since it corresponds to the only obtained non-convergent polyzeta
$\zeta \big( 1, a_1, 1^{ b_1} \big)$.
Suppose now that the $1$ is inserted {\em strictly inside}
the group of $0$:

\bigskip
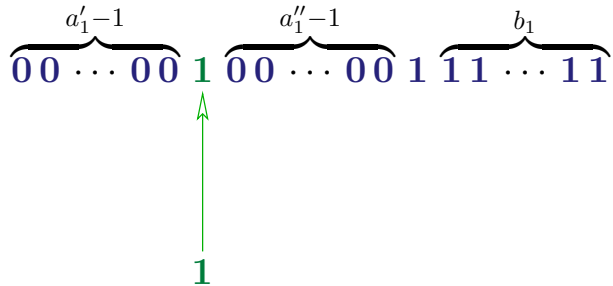
\begin{center}
\input strictly-inside-0.pstex_t
\end{center}

\noindent
with $a_1 ' - 1 \geqslant 1$ in the diagram, with
$a_1 '' - 1 \geqslant 1$ and with:
\[
a_1'-1+a_1''-1
=
a_1-1
=
\mathmotsf{unchanged total number of}\,\,0,
\]
so that one has:
\[
a_1'+a_1''=a_1+1.
\]
In other words, we have required that the 
inserted $1$ does not enter in contact neither
with the $1$ which terminates the series
$0 0 \cdots 0 0 1$ which encodes $a_1$,
nor with the $1$ which constitute $1^{ b_1}$,
and one obtains in sum:
\[
\sum_{a_1'+a_1''=a_1+1
\atop
a_1'\geqslant 2,\,\,
a_1''\geqslant 2}\,
\zeta
\big(
a_1',a_1'',1^{b_1}
\big).
\]

\medskip\noindent{\bf Principle of agglutination for shuffle-insertions 
of $1$.} {\em Any additional entry:
\[
\gbf{1}
\]
which is shuffle-inserted in some general polyzeta:
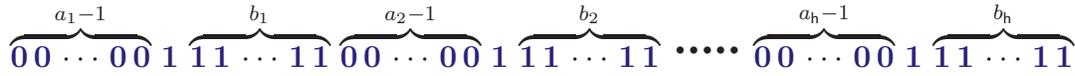
\[
\zsqov{a_1-1}\,
\unsqov{b_1}\,
\zsqov{a_2-1}\,
\unsqov{b_2}\,
\,\upbfdots\,
\zsqov{a_{\sf h}-1}\,
\unsqov{b_{\sf h}}\,
\]
at some place lying in direct contact with one of the $1$ of some
block $1\, 1^{ b_j}$: 

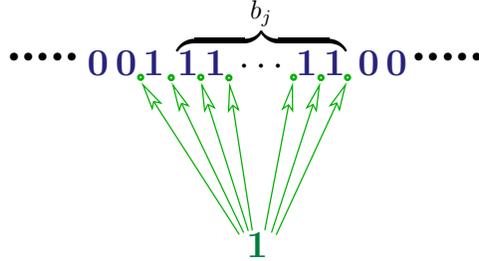
\bigskip
\begin{center}
\input 1-agglutination-1.pstex_t
\end{center}

\noindent
should be thought of as being automatically {\small\sf\em
agglutinated} to the group of $1$ of the $1\,1^{ b_j}$ in question, and
therefore, other insertions of such an additional $1$ in the $0$ of
some block $0^{ a_i - 1}$:

\bigskip
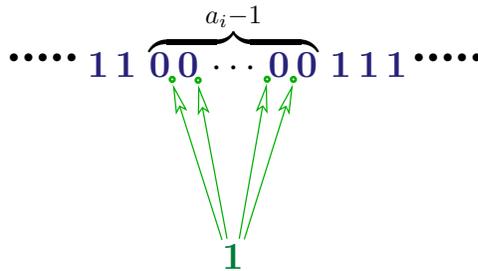
\begin{center}
\input re-strictly-inside-0.pstex_t
\end{center}

\noindent
should always be done
{\small\sf\em strictly inside} these $0$,
namely never in any two extremities of a block $0^{ a_i - 1}$.}

\medskip
In accordance with such a rule and coming back to our sample $(1)
\shuffle \big( a_1, 1^{ b_1} \big)$ studied in the simple case of
height ${\sf h} = 1$, it therefore only remains to look at all
insertions of $1$ which enter in contact with the $1^{ b_1}$:

\bigskip
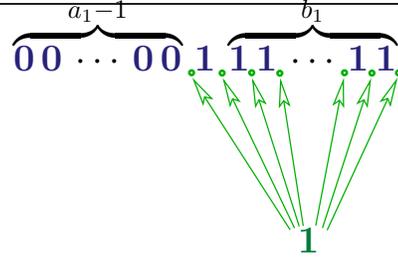
\begin{center}
\input 1-a1-b1.pstex_t
\end{center}

\noindent
But evidently, in each one of these $b_1 + 2$ circumstances, one
gets:
\[
\zsqov{a_1-1}\,\bbf{1}\,
\overbrace{\bbf{1}\,\bbf{1}\,
\cdots\,\gbf{1}\,\cdots\,\bbf{1}\,\bbf{1}}^{b_1+1}
=
\big(a_1,1^{b_1+1}\big),
\]
so that one obtains in sum:
\[
(b_1+2)\,
\zeta\big(a_1,1^{b_1+1}\big).
\]

\begin{Lemma}
For every $a_1 \geqslant 2$ and every $b_1 \geqslant 0$,
if one sets:
\[
a_1+b_1=:{\sf w}-1
\ \ \ \ \ \ \ \ \ \ \ \ \ 
\text{\rm and}
\ \ \ \ \ \ \ \ \ \ \ \ \ 
1+b_1=:{\sf d},
\]
then the $\shuffle$-shuffle product of $\zeta ( 1)$
with any $\zeta_{[ {\sf w}-1, {\sf d}, 1]}
\big( a_1, 1^{ b_1} \big)$ of weight ${\sf w}-1 \geqslant 2$, 
of depth ${\sf d}$ and of height
$1$ writes out, after complete finalization:
\[
\aligned
\zeta(1)
\shuffle
\zeta_{[{\sf w}-1,{\sf d},1]}
\big(a_1,1^{b_1}\big)
&
=
(b_1+2)\,\zeta_{[{\sf w},{\sf d}+1,1]}
\big(a_1,1^{b_1+1}\big)
+
\underbrace{\zeta_{[{\sf w},{\sf d}+1,1]}
\big(1,a_1,1^{b_1}\big)}_{\smallmathmotsf{non-convergent}}
+
\\
&
\ \ \ \ \
+
\sum_{a_1'+a_1''=a_1+1
\atop
a_1'\geqslant 2,\,\,
a_1''\geqslant 2}\,
\zeta_{[{\sf w},{\sf d}+1,2]}
\big(a_1',a_1'',1^{b_1}\big).
\qed
\endaligned
\]
\end{Lemma}

In the general case of arbitrary height ${\sf h} \geqslant 1$: 
\[
\gbf{1}
\shuffle
\zsqov{a_1-1}\,
\unsqov{b_1}\,
\upbfdots
\zsqov{a_{\sf h}-1}\,
\unsqov{b_{\sf h}},
\]
the reasonings are quite similar. Setting apart the insertion of $1$
in the front place which provides the only non-convergent polyzeta, one
firstly considers all possible insertions of $1$ strictly inside some
group of $0$ of some $0^{ a_i - 1}$:

\bigskip
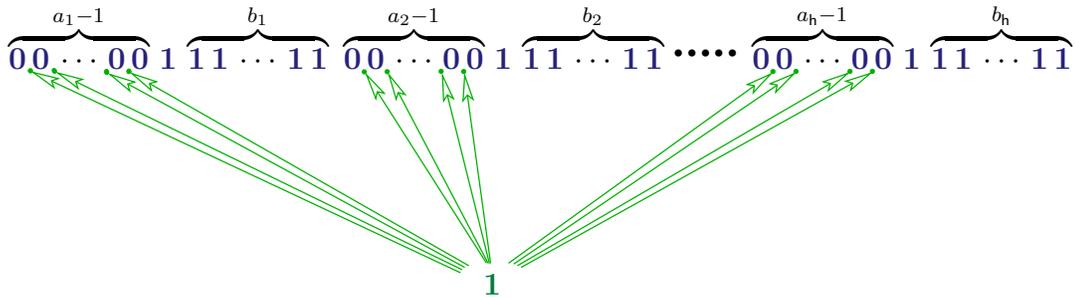
\begin{center}
\input 1-ins-1-in-00-00.pstex_t
\end{center}

\noindent
and this gives without delay:
\[
\sum_{i=1\atop a_i\geqslant 3}^{\sf h}\,
\Bigg(
\sum_{
a_i'+a_i''=a_i+1
\atop
a_i'\geqslant 2,\,\,a_i''\geqslant 2}\,
\zeta_{[{\sf w},{\sf d}+1,{\sf h}+1]}
\big(
\bfdots,
a_i',a_i'',
\bfdots
\big)
\Bigg).
\]
(Let us remind that implicitly, there is an $1^0$ between $a_i'$ and
$a_i''$ just here.) Secondly, in accordance with the principle of
agglutination stated above, one considers all possible insertions of $1$ which
happen to lie in direct contact with some $1$ of some
block $1\, 1^{ b_j}$:

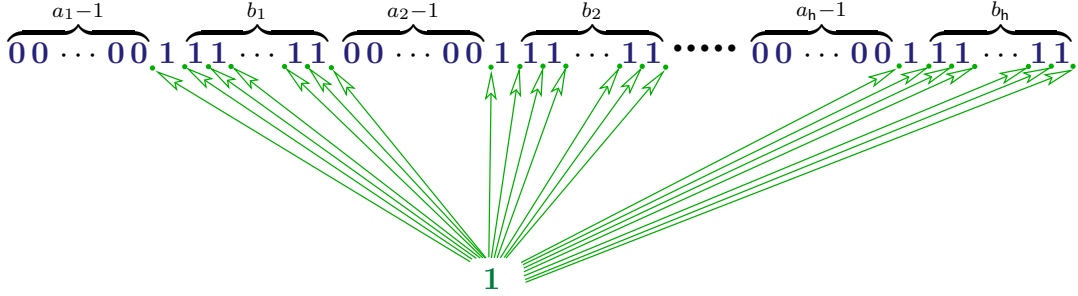
\bigskip
\begin{center}
\input 1-ins-1-in-111-11.pstex_t
\end{center}
 
\noindent
and naturally, one obtains a sum, for $j = 1, \dots, {\sf h}$, of
terms like the first one of the preceding lemma:
\[
\sum_{j=1}^{\sf h}\,
(b_j+2)\,
\zeta_{[{\sf w},{\sf d}+1,{\sf h}]}
\big(
\bfdots,1^{b_j+1},\bfdots
\big).
\]
All the explanations provided so far offer us the following result, 
fully proved now. 

\begin{Lemma}
For every height ${\sf h} \geqslant 1$,
every entries $a_1 \geqslant 2$, \dots, $a_{\sf h} \geqslant 2$,
every entries $b_1 \geqslant 0$, \dots, $b_{\sf h} \geqslant 0$,
if one sets:
\[
\aligned
a_1+b_1+\cdots+a_{\sf h}+b_{\sf h}
&
=:
{\sf w}-1,
\\
1+b_1+\cdots+1+b_{\sf h}
&
=:
{\sf d},
\endaligned
\]
then the $\shuffle$-shuffle product of $\zeta ( 1)$ with
any $\zeta_{[{\sf w}-1,{\sf d},{\sf h}]}
\big( a_1, 1^{ b_1}, \dots, a_{\sf h}, 1^{ b_{\sf h}} \big)$
of weight ${\sf w} - 1 \geqslant 2$, of depth ${\sf d}$ and of height ${\sf h}$
writes out, after complete finalization:
\[
\aligned
\zeta(1)
\shuffle
\zeta_{[{\sf w}-1,{\sf d},{\sf h}]}
\big(
a_1,1^{b_1},\dots,a_{\sf h},1^{b_{\sf h}}
\big)
&
=
\sum_{j=1}^{\sf h}\,
(b_j+2)\,
\zeta_{[{\sf w},{\sf d}+1,{\sf h}]}
\big(
\bfdots,1^{b_j+1},\bfdots
\big)
+
\\
&
\ \ \ \ \
+
\underbrace{\zeta_{[{\sf w},{\sf d}+1,{\sf h}]}
\big(1,a_1,1^{b_1},\dots,a_{\sf h},1^{b_{\sf h}}\big)}_{
\smallmathmotsf{non-convergent}
\atop
\smallmathmotsf{but will disappear}}
+
\\
&
\ \ \ \ \ 
+
\sum_{i=1\atop a_i\geqslant 3}^{\sf h}\,
\Bigg(
\sum_{
a_i'+a_i''=a_i+1
\atop
a_i'\geqslant 2,\,\,a_i''\geqslant 2}\,
\zeta_{[{\sf w},{\sf d}+1,{\sf h}+1]}
\big(
\bfdots,
a_i',a_i'',
\bfdots
\big)
\Bigg).
\endaligned
\]
\end{Lemma}

A concrete example\,\,---\,\,respecting the above
three-lines writing\,\,---\,\,is:
\[
\aligned
(1)\shuffle\big(3,1,4,1\big)
&
=
3\,\big(3,1,1,4,1\big)
+
3\,\big(3,1,4,1,1\big)
+
\\
&
\ \ \ \ \
+
\big(1,3,1,4,1\big)
+
\\
&
\ \ \ \ \
+
\big(2,2,1,4,1\big)
+
\big(3,1,3,2,1\big)
+
\big(3,1,2,3,1\big).
\endaligned
\]

\subsection{The subtraction}
Now that we have computed both the stuffle
and the shuffle production of $\zeta( 1)$ with any
$\zeta \big( a_1, 1^{ b_1}, \dots, a_{\sf h}, 
1^{ b_{\sf h}} \big)$, it only remains to perform the
final subtraction of the two results provided
by the lemmas above:
\[
\footnotesize
\aligned
0
&
=
-\,(1)
\ast
\big(a_1,1^{b_1},\dots,a_{\sf h},1^{b_{\sf h}}\big)
+
1
\shuffle
\big(a_1,1^{b_1},\dots,a_{\sf h},1^{b_{\sf h}}\big)
\\
&
=
-\,\sum_{i=1}^{\sf h}\,
\big(\bfdots,a_i+1,\bfdots\big)
-
\\
&
\ \ \ \ \
-
\sum_{j=1\atop b_j\geqslant 1}^{\sf h}
\Bigg(
\sum_{b_j'+b_j''=b_j-1
\atop
b_j'\geqslant 0,\,\,b_j''\geqslant 0}\,
\big(\bfdots,1^{b_j'},2,1^{b_j''},\bfdots\big)
\Bigg)
-
\\
&
\ \ \ \ \
-
\zero{\big(1,a_1,1^{b_1},\dots,a_{\sf h},1^{b_{\sf h}}\big)}
-
\sum_{j=1}^{\sf h}\,
(b_j+1)\,
\big(\bfdots,1^{b_j+1},\bfdots\big)
+
\\
&
\ \ \ \ \
+
\sum_{j=1}^{\sf h}\,(b_j+2)\,
\big(\bfdots,1^{b_j+1},\bfdots\big)
+
\\
&
\ \ \ \ \
+
\zero{\big(1,a_1,1^{b_1},\dots,a_{\sf h},1^{b_{\sf h}}\big)}
+
\\
&
\ \ \ \ \
+
\sum_{i=1\atop a_i\geqslant 3}^{\sf h}
\Bigg(
\sum_{a_i'+a_i''=a_i+1
\atop
a_i'\geqslant 2,\,\,a_i''\geqslant 2}\,
\big(\bfdots,a_i',a_i'',\bfdots\big)
\Bigg).
\endaligned
\]
As is known and as was expected, the two (underlined) non-convergent
polyzetas annihilate.  Also, the terms involving the multiplicities
$-(b_j+1)$ in the fourth line and the ones involving the
multiplicities $(b_j+2)$ in the fifth line immediately collapse, while
no other terms simplify.

With the mentions of weights, of depths and of heights, we therefore
have gained the first fundamental theorem of the present article,
classically attributed to Hoffman, but exhibited here in a more
detailed way, including rigorous explicitation of quantifiers.

\begin{Theorem}
For every height ${\sf h} \geqslant 1$,
every entries $a_1 \geqslant 2$, \dots, $a_{\sf h} \geqslant 2$,
every entries $b_1 \geqslant 0$, \dots, $b_{\sf h} \geqslant 0$,
if one sets:
\[
\aligned
a_1+b_1+\cdots+a_{\sf h}+b_{\sf h}
&
=:
{\sf w}-1,
\\
1+b_1+\cdots+1+b_{\sf h}
&
=:
{\sf d},
\endaligned
\]
then the so-called double shuffle relation: 
\[
0
=
-\,(1)
\ast
\big(a_1,1^{b_1},\dots,a_{\sf h},1^{b_{\sf h}}\big)
+
1
\shuffle
0^{a_1-1}1\,1^{b_1}\,\cdots\,0^{a_{\sf h}-1}1\,1^{b_{\sf h}}
\]
between $\zeta ( 1)$ and any $\zeta_{[{\sf w}-1,{\sf d},{\sf h}]}
\big( a_1, 1^{ b_1}, \dots, a_{\sf h}, 1^{ b_{\sf h}} \big)$ of weight
${\sf w} - 1 \geqslant 2$, of depth ${\sf d}$ 
and of height ${\sf h}$ writes out after complete
finalization:
\[
\boxed{
\aligned
0
&
=
-\,\sum_{i=1}^{\sf h}\,
\zeta_{[{\sf w},{\sf d},{\sf h}]}
\big(\bfdots,a_i+1,\bfdots\big)
-
\\
&
\ \ \ \ \
-
\sum_{j=1\atop b_j\geqslant 1}^{\sf h}
\Bigg(
\sum_{b_j'+b_j''=b_j-1}\,
\zeta_{[{\sf w},{\sf d},{\sf h}+1]}
\big(\bfdots,1^{b_j'},2,1^{b_j''},\bfdots\big)
\Bigg)
+
\\
&
\ \ \ \ \
+
\sum_{j=1}^{\sf h}\,
\zeta_{[{\sf w},{\sf d}+1,{\sf h}]}
\big(\bfdots,1^{b_j+1},\bfdots\big)
+
\\
&
\ \ \ \ \
+
\sum_{i=1\atop a_i\geqslant 3}^{\sf h}
\Bigg(
\sum_{a_i'+a_i''=a_i+1
\atop
a_i'\geqslant 2,\,\,a_i''\geqslant 2}\,
\zeta_{[{\sf w},{\sf d}+1,{\sf h}+1]}
\big(\bfdots,a_i',a_i'',\bfdots\big)
\Bigg).
\endaligned}
\]
\end{Theorem}

In fact, the standard statement of Hoffman's relations ({\em cf.}
{\em e.g.} \cite{Zudilin-2003}, p.~3):
\[
\aligned
0
&
=
-
\sum_{k=1}^{\sf d}\,
\zeta\big(s_1,\dots,s_{k-1},s_k+1,s_{k+1},\dots,s_{\sf d}\big)
+
\\
&
\ \ \ \ \
+
\sum_{k=1\atop s_k\geqslant 2}^{\sf d}\,
\sum_{j=0}^{s_k-2}\,
\zeta\big(s_1,\dots,s_{k-1},s_k-j,j+1,s_{k+1},\dots,s_{\sf d}\big),
\endaligned
\]
usually makes no mention of the variable depths and of the variable
heights of the appearing polyzetas, while the summations are not
organized according to any ordering between all polyzetas of fixed
weight.

\smallskip

So is our first\,\,---\,\,by far the easiest amongst
six\,\,---\,\,theorem.

\section{Computing $-(2) \ast \big(a_1,1^{b_1}, \dots,a_{\sf h},
1^{b_{\sf h}}\big) 
+ 01 \shuffle 0^{ a_1 -1}1\, 1^{b_1} \cdots
0^{a_{\sf h}-1}1\, 1^{b_{\sf h}}$}
\label{2-double-melange}
\HEAD{Computing $-(2) \ast \big(a_1,1^{b_1}, \dots,a_{\sf h},
1^{b_{\sf h}}\big) 
+ 01 \shuffle 0^{ a_1 -1}1\, 1^{b_1} \cdots
0^{a_{\sf h}-1}1\, 1^{b_{\sf h}}$}{
Jo\"el {\sc Merker}, Universit\'e Paris-Sud Orsay, France}

\subsection{Computing firstly the stuffle product 
$(2) \ast \big(a_1,1^{b_1}, \dots,a_{\sf h},
1^{b_{\sf h}}\big)$}
What next goal should be is clear: write out the general double
shuffle relations with $\zeta ( 2)$ as a first member.

We start by treating the $\ast$-stuffle product, in the (simpler)
sample case of height ${\sf h} = 1$. Exactly as in the case of $\zeta
( 1) \ast \big( a_1, 1^{ b_1} \big)$, the computation of $\zeta ( 2)
\ast \big( a_1, 1^{ b_1} \big)$ involves addition terms and insertion
terms:
\[
\aligned
(1)\ast\big(a_1,1^{b_1}\big)
&
=
\big(a_1+2,1^{b_1}\big)
+
\big(a_1,3,1^{b_1-1}\big)
+\cdots+
\big(a_1,1^{b_1-1},3\big)
+
\ \ \ \ \ \ \ \ \ \
\explain{additions}
\\
&
\ \ \ \ \
+
\big(2,a_1,1^{b_1}\big)
+
\big(a_1,2,1^{b_1}\big)
+\cdots+
\big(a_1,1^{b_1},2\big),
\ \ \ \ \ \ \ \ \ \ \ \ \ \ \ \ \ \ \ \,
\explain{insertions}
\endaligned
\]
the only difference being that no non-convergent polyzeta appears, now.
Using summation symbols, the result may then be expressed as:
\[
\aligned
(2)\ast\big(a_1,1^{b_1}\big)
&
=
\big(a_1+2,1^{b_1}\big)
+
\sum_{b_1'+b_1''=b_1-1
\atop
b_1'\geqslant 0,\,\,b_1''\geqslant 0}\,
\big(a_1,1^{b_1'},3,1^{b_1''}\big)
+
\\
&
\ \ \ \ \
+
\big(2,a_1,1^{b_1}\big)
+
\sum_{b_1'+b_1''=b_1
\atop
b_1'\geqslant 0,\,\,b_1''\geqslant 0}\,
\big(a_1,1^{b_1'},2,1^{b_1''}\big).
\endaligned
\]
The general case of arbitrary height ${\sf h} \geqslant 1$ is 
quite similar, hence we state it directly without additional words.

\begin{Lemma}
For every height ${\sf h} \geqslant 1$, every entries $a_1 \geqslant 2$,
\dots, $a_{\sf h} \geqslant 2$, every entries $b_1 \geqslant 0$, \dots,
$b_{\sf h} \geqslant 0$, if one sets:
\[
\aligned
a_1+b_1+\cdots+a_{\sf h}+b_{\sf h}
&
=:
{\sf w}-2,
\\
1+b_1+\cdots+1+b_{\sf h}
&
=:
{\sf d},
\endaligned
\]
then the $*$-stuffle product of $\zeta (2)$ with any $\zeta_{[ {\sf
w}-2, {\sf d}, {\sf h}]} \big( a_1, 1^{ b_1}, \dots, a_{\sf h},
1^{ b_{\sf h}} \big)$ of weight ${\sf w} - 2 \geqslant 2$,
of depth ${\sf d}$ and
of height ${\sf h}$ writes out, after complete finalization:
\[
\aligned
\zeta(2)
\ast
\zeta_{[{\sf w}-2,{\sf d},{\sf h}]}
\big(a_1,1^{b_1},\dots,a_{\sf h},1^{b_{\sf h}}\big)
&
=
\sum_{i=1}^{\sf h}\,
\zeta_{[{\sf w},{\sf d},{\sf h}]}
\big(\bfdots,a_i+2,\bfdots\big)
+
\\
&
\ \ \ \ \
+
\sum_{j=1\atop b_j\geqslant 1}^{\sf h}\,
\Bigg(
\sum_{b_j'+b_j''=b_j-1
\atop
b_j'\geqslant 0,\,\,b_j''\geqslant 0}\,
\zeta_{[{\sf w},{\sf d},{\sf h}+1]}
\big(\bfdots,1^{b_j'},3,1^{b_j''},\bfdots\big)
\Bigg)
+
\\
&
\ \ \ \ \
+
\zeta_{[{\sf w},{\sf d}+1,{\sf h}+1]}
\big(2,a_1,1^{b_1},\dots,a_{\sf h},1^{b_{\sf h}}\big)
+
\\
&
\ \ \ \ \
+
\sum_{j=1}^{\sf h}
\Bigg(
\sum_{b_j'+b_j''=b_j
\atop
b_j'\geqslant 0,\,\,b_j''\geqslant 0}\,
\zeta_{[{\sf w},{\sf d}+1,{\sf h}+1]}
\big(\bfdots,1^{b_j'},2,1^{b_j''},\bfdots\big)
\Bigg).
\endaligned
\]
\end{Lemma} 

Notice that in the right-hand side the polyzetas appear firstly
according to increasing depth, secondly according to increasing height.

\subsection{Computing secondly the shuffle product 
$01 \shuffle 0^{ a_1-1} 1\, 1^{b_1} \cdots 0^{ a_{\sf h}-1} 1\,
1^{b_{\sf h}}$} \label{shuffle-01} Now, the computational task
begins to be more substantial. One has to shuffle-insert $01$:
\[
\gbf{01}
\shuffle
\zsqov{a_1-1}\,\unsqov{b_1}\,
\upbfdots\,
\zsqov{a_{\sf h}-1}\,\unsqov{b_{\sf h}}
\]
in a general polyzeta of weight, say, equal to: 
\[
a_1+b_1+\cdots+a_{\sf h}+b_{\sf h}
=:
{\sf w} - 2, 
\]
of depth, say, equal to: 
\[
1+b_1+\cdots+1+b_{\sf h}
=:
{\sf d},
\]
and of height visibly equal to ${\sf h}$. Of course, the total number
of obtained polyzetas must be equal to the number of choices of two
(binary) digits among ${\sf w} - 2 + 2 = {\sf w}$ digits, namely:
\[
\aligned
\mathmotsf{number of terms}
&
=
\binom{\sf w}{2}
=
\frac{{\sf w}({\sf w}-1)}{2}
\\
&
=
\frac{
\big(a_1+b_1+\cdots+a_{\sf h}+b_{\sf h}+2\big)
\big(a_1+b_1+\cdots+a_{\sf h}+b_{\sf h}+1\big)}{2}.
\endaligned
\]
Another simple initial observation is:

\begin{Lemma}
All terms of $01 \shuffle 0^{ a_1 - 1} 1\, 1^{ b_1}
\cdots 0^{ a_{\sf h}-1} 1 \, 1^{ b_{\sf h}}$ are of depth ${\sf d} +1$.
\end{Lemma}

\proof
Remembering that the depth of a general polyzeta is:
\[
\aligned
{\sf d}
&
=
1+b_1+\cdots+1+b_{\sf h}
\\
&
=
\mathmotsf{number of}\,1\,\mathmotsf{in its encoding}\,
0^{a_1-1}1\,1^{b_1}
\cdots
0^{a_{\sf h}-1}1\,1^{b_{\sf h}},
\endaligned
\]
it is then clear that any shuffle-insertion of $01$
always increases the number of $1$ by one unit, exactly.
\endproof

How will we, then, insert the $0$ and the $1$ of $01$ inside a general
polyzeta $0^{ a_1 - 1} 1\, 1^{ b_1} \cdots 0^{ a_{\sf h}-1} 1 \, 1^{
b_{\sf h}}$?  The principle of agglutination stated in the preceding
section for insertions of $1$ also has a mirror-companion concerning
insertions of $0$.

\medskip\noindent{\bf Principle of agglutination 
for shuffle-insertions of $0$.} {\em Any additional entry: 
\[
\gbf{0}
\] 
which is
shuffle-inserted in some general polyzeta:
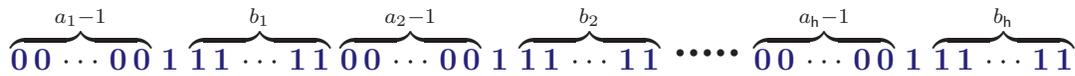
\[
\zsqov{a_1-1}\,
\unsqov{b_1}\,
\zsqov{a_2-1}\,
\unsqov{b_2}\,
\,\upbfdots\,
\zsqov{a_{\sf h}-1}\,
\unsqov{b_{\sf h}}\,
\]
at some place lying in direct contact with one of the $0$ of some
block $0^{ a_i - 1}$: 

\bigskip
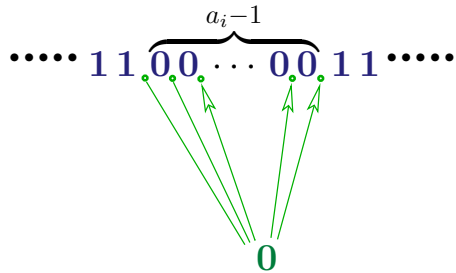
\begin{center}
\input 0-agglutination-0.pstex_t
\end{center}

\noindent
should be thought of as being automatically {\small\sf\em
agglutinated} to the group of $0$ of the block $0^{ a_i-1}$ in
question, and therefore, other insertions of such an additional $0$ in
the $1$ of some block $1\,1^{ b_j}$:

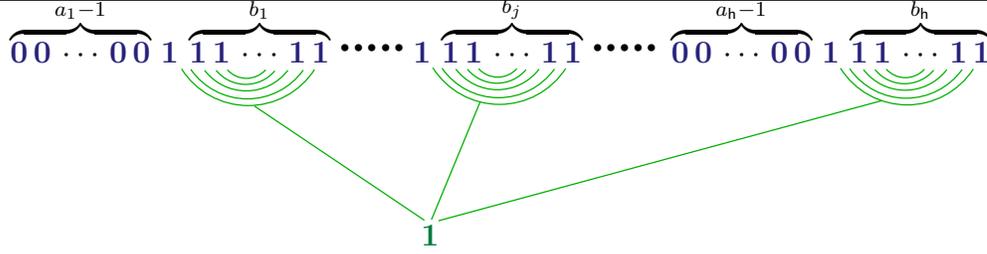
\bigskip
\begin{center}
\input re-strictly-inside-1.pstex_t
\end{center}

\noindent
should always be done {\small\sf\em strictly inside} these $1$, namely
never in any two extremities of a block $1\,1^{ b_j}$.}

\medskip
So we have to insert the $0$ and the $1$ of $01$
one after the other in a general polyzeta. The possibilities are:

\medskip$\square$
the $0$ goes into the $0$:

\bigskip
\begin{center}
\input 0-goes-in-the-00-00.pstex_t
\end{center}

\medskip$\square$
the $0$ goes into the $1$:

\bigskip
\begin{center}
\input 0-goes-in-the-111-11.pstex_t
\end{center}

\medskip$\square$
the $1$ goes into the $0$:

\bigskip
\begin{center}
\input 1-goes-in-the-00-00.pstex_t
\end{center}

\medskip$\square$
the $1$ goes into the $1$:

\bigskip
\begin{center}
\input 1-goes-in-the-111-11.pstex_t
\end{center}

Since one has to insert the $0$ and the $1$ of $01$ {\em
simultaneously one after the other} in a general polyzeta, four {\em
common} occurrences hold:

\medskip
$\boxed{i_1 \vert i_2}$\,:
the $0$ goes into the $0$ of some $0^{ a_{i_1} -1}$
and the $1$ also goes into the $0$ of some $0^{ a_{ i_2} - 1}$;

\medskip
$\boxed{i\vert j}$\,:
the $0$ goes into the $0$ of some $0^{ a_i -1}$
while the $1$ goes into the $1$ of some $1\, 1^{ b_j}$; 

\medskip
$\boxed{j\vert i}$\,:
the $0$ goes into the $1$ of some $1\,1^{ b_j}$ while the $1$ goes 
into the $0$ of some $0^{ a_i - 1}$; 

\medskip
$\boxed{j_1\vert j_2}$\,:
the $0$ goes into the $1$ of some
$1\, 1^{ b_{ j_1}}$ and the $1$ also goes into the $1$
of some $1\, 1^{ b_{ j_2}}$.

\medskip
Of course there also are the two somehow special subcases of the first
one and of the fourth one, respectively, namely when $i_1 = i_2$
and when $j_1 = j_2$: 

\medskip
$\boxed{ i\, i}$\,:
the $0$ goes into the $0$ of some $0^{ a_i -1}$
and the $1$ also goes into the $0$ of the {\em same} $0^{ a_i - 1}$;

\medskip
$\boxed{ j\, j}$\,:
the $0$ goes into the $1$ of some
$1\, 1^{ b_j}$ and the $1$ also goes into the $1$
of the {\em same} $1\, 1^{ b_j}$.

\medskip\noindent
These two last subcases will be treated separately,
hence precisely six cases will be dealt with, in fact.

\medskip
Although the depth of any polyzeta appearing in $01 \shuffle 0^{ a_1 -
1} 1 \, 1^{ b_1} \cdots 0^{ a_{\sf h}-1} 1 1^{ b_{\sf h}}$ is
always equal to ${\sf d} + 1$, the heights of such appearing polyzetas
can vary, as examples show. {\em A precise control of the heights in
question\,\,---\,\,to be made below\,\,---\,\,happens to be
available only when one applies the two principles agglutinations,
for the insertion of the $0$ of $01$ and as well, for the insertion
of the $1$ of $01$.} 
Complete explanations
being provided in just a while, let us list
in advance the obtained heights\,\,---\,\,we intentionally change the
order of appearance of the
six cases, so as to fit with the ordering we introduced in
Section~\ref{total-ordering-polyzetas}\,\,---\,\,:

\medskip$\square$\,
First family $\boxed{i\vert j}$ with $1 \leqslant i \leqslant j \leqslant
{\sf h}$\,: depth ${\sf d} + 1$; height ${\sf h}$.

\medskip$\square$\,
Second family $\boxed{j_1\vert j_2}$ with $1 \leqslant j_1 < j_2
\leqslant {\sf h}$: depth ${\sf d} + 1$; height ${\sf h} + 1$;

\medskip$\square$\,
Third family $\boxed{j\, j}$ with $1 \leqslant j \leqslant {\sf h}$:
depth ${\sf d} + 1$; height ${\sf h} + 1$;

\medskip$\square$\,
Fourth family $\boxed{i\, i}$ with $1 \leqslant i\leqslant {\sf h}$:
depth ${\sf d} + 1$; height ${\sf h} + 1$;

\medskip$\square$\,
Fifth family $\boxed{i_1\vert i_2}$ with $1 \leqslant i_1 < i_2
\leqslant {\sf h}$: depth ${\sf d} + 1$; height ${\sf h} + 1$;

\medskip$\square$\,
Sixth family $\boxed{ j \vert i}$ with $1 \leqslant j < i \leqslant 
{\sf h}$: depth ${\sf d} + 1$; height ${\sf h} + 2$.

\subsection{First family $\boxed{i \vert j}$\,}
Thus, assume that $0$ goes into a block of $0$ while $1$ goes into
a block of $1$: 

\bigskip
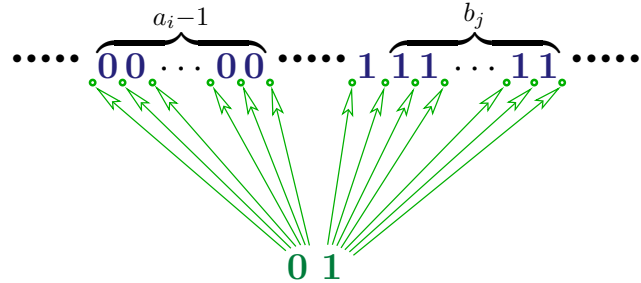
\begin{center}
\input 0-goes-00-1-goes-111.pstex_t
\end{center}

\noindent
Of course, one has $i \leqslant j$ here, for $0$ must always be inserted left
to the insertion of $1$ (shuffle rule).
For arbitrary fixed $i$ and $j$ with $1
\leqslant i \leqslant j \leqslant {\sf h}$, it is then
diagrammatically visible and rigorously clear that for every such
simultaneous insertion of $01$, one obtains the same
polyzeta, and the result is:
\[
a_i\,(b_j+2)\cdot
\zeta
\Big(
\bfdots\,
\overbrace{\bbf{0}\,\bbf{0}\,\cdots\,\gbf{0}\,
\cdots\,\bbf{0}\,\bbf{0}}^{a_i}\,
\bfdots
\bbf{1}\,\overbrace{\bbf{1}\,\bbf{1}\,\cdots\,
\gbf{1}\,\cdots\,\bbf{1}\,\bbf{1}}^{b_j+1}
\bfdots
\Big).
\]
One easily convinces oneself that the height is unchanged. As a result,
the polyzetas generated by this first {\sf f}amily\,\,---\,\,with the 
letter $\mathcal{ F}$\,\,---\,\,may be collected using
a triangular summation symbol: 
\[
\boxed{
\mathcal{F}_{i\vert j}
:=
\sum_{1\leqslant i\leqslant j\leqslant{\sf h}}\,
a_i\,(b_j+2)\,
\zeta_{[{\sf w},{\sf d+1},{\sf h}]}
\big(
\bfdots,a_i+1,\bfdots,1^{b_j+1},\bfdots
\big).}
\]
We recall that the bold dots mean that other entries of the polyzeta
are unchanged. Lastly, we observe {\em passim} that the total number
of terms in this family $\mathcal{ F}_{ i\vert j}$ equals:
\[
{\bf n}_{i\vert j}
:=
\sum_{1\leqslant i\leqslant j\leqslant{\sf h}}\,
a_i\,b_j
+
\sum_{1\leqslant i\leqslant j\leqslant{\sf h}}\,
2\,a_i.
\]
In fact, in Subsection~\ref{final-counting-checking} below,
we will check that the total sum of the numbers of
polyzetas appearing in each one of our six families:
\[
\footnotesize
\aligned
{\bf n}_{i\vert j}
+
{\bf n}_{j_1\vert j_2}
+
{\bf n}_{j\,j}
+
{\bf n}_{i\,i}
+
{\bf n}_{i_1\vert i_2}
+
{\bf n}_{j\vert i}
&
=
\frac{{\sf w}({\sf w}-1)}{2}
\\
&
=
\frac{\big(a_1+b_1+\cdots+a_{\sf h}+b_{\sf h}+2\big)
\big(a_1+b_1+\cdots+a_{\sf h}+b_{\sf h}+1\big)}{2}
\endaligned
\]
is indeed equal to the total expected number shown at the beginning of
Subsection~\ref{shuffle-01}.  In the future two
Sections~\ref{3-double-melange} and~\ref{21-double-melange} below, the
checkings will become a bit harder when dealing with $001 \shuffle 0^{
a_1 - 1} 1 \, 1^{ b_1}\, \cdots \, 0^{ a_{\sf h} - 1} 1 \, 1^{
b_{\sf h}}$ and with $011 \shuffle 0^{ a_1 - 1} 1 \, 1^{ b_1}\,
\cdots \, 0^{ a_{\sf h} - 1} 1 \, 1^{ b_{\sf h}}$ (respectively),
hence let us admit that the present computational level remains
accessible.

\subsection{Second family $\boxed{j_1 \vert j_2}$ 
with $1 \leqslant j_1 < j_2 \leqslant {\sf h}$\,}
Suppose now that the $0$ and the $1$ of $01$ both 
go into a block of $1$, the blocks being {\em distinct},
say into $1\, 1^{ b_{ j_1}}$ and into $1\, 1^{ b_{ j_2}}$
(respectively), for some $1 \leqslant j_1 < j_2 \leqslant {\sf h}$.

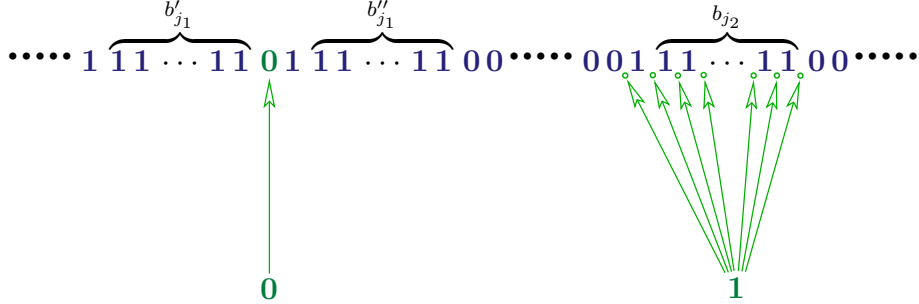
\bigskip
\begin{center}
\input 0-j1-1-j2.pstex_t
\end{center}

\noindent
According to the agglutination principles, the $1$ on the
right can take $b_{ j_2} + 2$ positions, while the $0$ on the
left must be inserted {\em strictly} inside $1\, 1^{ b_{ j_1}}$.
Consequently, there are two integers $b_{ j_1}' \geqslant 0$
and $b_{ j_1}'' \geqslant 0$ as in the diagram, 
which satisfy\,\,---\,\,notice that one more $1$ 
is kept just after the inserted $0$\,\,---\,:
\[
b_{j_1}'+b_{j_1}''=b_{j_1}-1.
\]
As a result, the original $1^{ b_{ j_1}}$ is replaced by $1^{
b_{j_1}'}, 2, 1^{ b_{ j_1}''}$, while the $1^{ b_{ j_2}}$ is
replaced by $1^{ b_{ j_2} + 1}$, with multiplicity $(b_{ j_2} +
2)$. Thus, the polyzetas generated by this second family may be
collected as follows using two summations symbols:
\[
\boxed{
\mathcal{F}_{j_1\vert j_2}
:=
\sum_{1\leqslant j_1<j_2\leqslant{\sf h}\atop b_{j_1}\geqslant 1}
\Bigg(
\sum_{b_{j_1}'+b_{j_1}''=b_{j_1-1}
\atop
b_{j_1}'\geqslant 0,\,\,b_{j_1}''\geqslant 0}\,
(b_{j_2}+2)
\cdot
\zeta_{[{\sf w},{\sf d}+1,{\sf h}+1]}
\big(
\bfdots,1^{b_{j_1'}},2,1^{b_{j_1}''},\bfdots,
1^{b_{j_2+1}},\bfdots
\big)
\Bigg).
}
\]
Lastly, let us observe that the total number of
terms in this second family $\mathcal{ F}_{ j_1 \vert j_2}$
equals:
\[
{\bf n}_{j_1\vert j_2}
:=
\sum_{1\leqslant j_1<j_2\leqslant{\sf h}}\,
b_{j_1}\,(b_{j_2}+2).
\]

\subsection{Third family $\boxed{j\,j}$ with 
$1 \leqslant j \leqslant {\sf h}$\,}
Next, assume that the $0$ and the $1$ of $01$ go together
into the {\em same} block of $1$, a case
which would correspond to the limit case
$j_1 = j_2 =: j$ in the second family that just precedes:

\bigskip
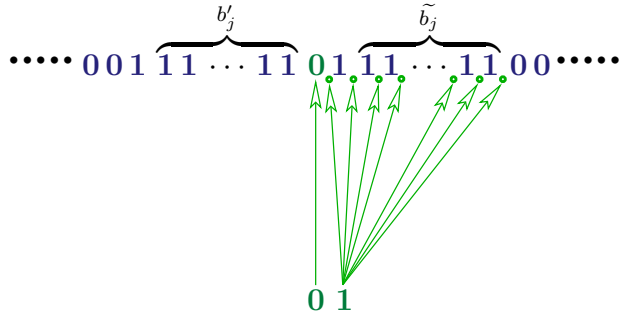
\begin{center}
\input 0-j-1-j.pstex_t
\end{center}

\noindent
Thus, there are two integers $b_j' \geqslant 0$ and
$\widetilde{ b}_j \geqslant 0$ as shown in the diagram 
satisfying:
\[
b_j'+\widetilde{b}_j
=
b_j-1.
\]
As shuffle rules dictate, the insertion of the $1$ must be
done after that of the $0$. Once this is performed, there is
one more $1$ in the right group:
\[
\bfdots\,
\bbf{0}\,\bbf{0}\,
\unsqov{b_j'}\,
\gbf{0}\,
\overbrace{\bbf{1}\,\bbf{1}\,\bbf{1}\,\cdots\,\gbf{1}\,\cdots\,
\bbf{1}\,\bbf{1}}^{b_j''}\,
\bbf{0}\,\bbf{0}\,\bfdots,
\]
namely we may set:
\[
b_j''
:=
\widetilde{b}_j+1,
\]
so that $b_j'$ and $b_j''$ now satisfy:
\[
b_j'+b_j''
=
b_j.
\]
As a result, the polyzetas generated by this third family
may be collected as follows using two summation symbols:
\[
\boxed{
\mathcal{F}_{j\,j}
:=
\sum_{1\leqslant j\leqslant{\sf h}}
\Bigg(
\sum_{b_j'+b_j''=b_j
\atop
b_j'\geqslant 0,\,\,b_j''\geqslant 0}\,
(b_j''+1)\cdot
\zeta_{[{\sf w},{\sf d}+1,{\sf h}+1]}
\big(\bfdots,1^{b_j'},2,1^{b_j''},\bfdots
\big)
\Bigg).
}
\]
Clearly, the number of terms present here equals:
\[
{\bf n}_{j\,j}
:=
\sum_{1\leqslant j\leqslant{\sf h}}\,
\frac{(b_j+2)(b_j+1)}{2}.
\]
But before passing to the fourth family, we must point out
in advance that, for each such fixed block $1\, 1^{ b_j}$, 
the last obtained polyzeta, namely the one
for $b_j' = b_j$ and $b_j'' = 0$, happens to be the following
special polyzeta which might appear elsewhere
because of some possible ambiguity:
\[
\bfdots\,
\bbf{0}\,\bbf{0}\,
\unsqov{b_j}\,
\gbf{0\,1}\,
\bbf{0}\,\bbf{0}\,
\bfdots.
\]
Indeed, there is no reason that the families we delineated above do
not have small overlaps, and in fact, there will be some (restricted)
overlaps (only) between the third and the fourth families\,\,---\,\,but 
we will
avoid them.

\subsection{Fourth family $\boxed{i\,i}$ with 
$1 \leqslant i \leqslant {\sf h}$\,}
Now, the $0$ and the $1$ of the $01$ are assumed
to go into the {\em same} group of $0$, say
into some block $0^{ a_i - 1}$:

\bigskip
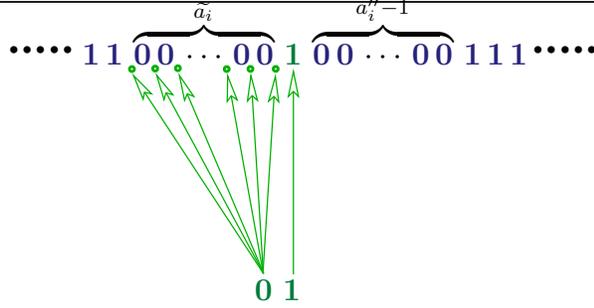
\begin{center}
\input 0-i-1-i.pstex_t
\end{center}

\noindent
Again, the principles of agglutination command to insert
firstly the $1$ strictly inside this block
$0^{ a_i - 1}$, so that
in the exhibited diagram the following (in)equalities
must hold:
\[
\widetilde{a_i}
\geqslant 1,
\ \ \ \ \ \ \ \
a_i''-1\geqslant 1,
\ \ \ \ \ \ \ \ \
\text{\rm and of course}
\ \ \ \ \ \ \ \ \
\widetilde{a_i}+a_i''-1
=
a_i-1.
\]
Then the $0$ is inserted at any place before the
insertion of the $1$, and there are $\widetilde{ a_i} + 1$
possibilities giving the same polyzeta:
\[
\bfdots\,
\bbf{1}\,\bbf{1}\,
\overbrace{\bbf{0}\,\bbf{0}\,\cdots\,\gbf{0}\,\cdots\,
\bbf{0}\,\bbf{0}}^{a_i'-1}\,
\gbf{1}\,
\zsqov{a_i''-1}\,\bbf{1}\,\bbf{1}\,\bbf{1}\,
\bfdots,
\]
in which it is natural to set:
\[
a_i'-1
:=
\widetilde{a}_i+1,
\ \ \ \ \ \ \ \ \ \ \ \ \
\text{\rm that is to say:}
\ \ \ \ \ \ \ \ \
a_i'
=
\widetilde{a_i}+2.
\]
Doing so, with $a_i' \geqslant 1$, one excludes the polyzeta:
\[
\bfdots\,\bbf{1}\,\bbf{1}\,
\gbf{01}\,\zsqov{a_i-1}\,\bbf{1}\,\bbf{1}\,\bbf{1}\,
\bfdots
\]
which appeared already in the third family, {\em except when $i = 1$,
so that one must add:}
\[
\gbf{0\,1}\,\zsqov{a_1-1}\,\unsqov{b_1}\,
\bfdots\,
\zsqov{a_{\sf h}-1}\,\unsqov{b_{\sf h}}
=
\big(2,a_1,1^{b_1},\dots,a_{\sf h},1^{b_{\sf h}}\big).
\]
As a result, without any overlapping with what precedes, 
the polyzetas of the fourth family are:
\[
\boxed{
\aligned
\mathcal{F}_{i\,i}
&
:=
\big(
2,a_1,1^{b_1},\dots,a_{\sf h},1^{b_{\sf h}}
\big)
+
\\
&
\ \ \ \ \
+
\sum_{1\leqslant i\leqslant{\sf h}\atop a_i\geqslant 5}
\Bigg(
\sum_{a_i'+a_i''=a_i+2
\atop
a_i'\geqslant 3,\,\,a_i''\geqslant 2}\,
(a_i'-1)\cdot
\zeta_{[{\sf w},{\sf d}+1,{\sf h}+1]}
\big(
\bfdots,a_i',1^0,a_i',\bfdots
\big)
\Bigg).
\endaligned
}
\]
Lastly, the total number of terms in this fourth family may readily be
checked to be equal to:
\[
{\bf n}_{i\,i}
:=
1
+
\sum_{1\leqslant i\leqslant{\sf h}}
\bigg(
\frac{a_i(a_i-1)}{2}-1
\bigg).
\]

\subsection{Fifth family $\boxed{i_1 \vert i_2}$ 
with $1 \leqslant i_1 < i_2 \leqslant {\sf h}$\,} Next, assume that
the $0$ and the $1$ of $01$ go into {\em distinct} blocks $0^{ a_{
i_1} - 1}$ and $0^{ a_{ i_2} - 1}$, for some $1 \leqslant i_1 <
i_2 \leqslant {\sf h}$:

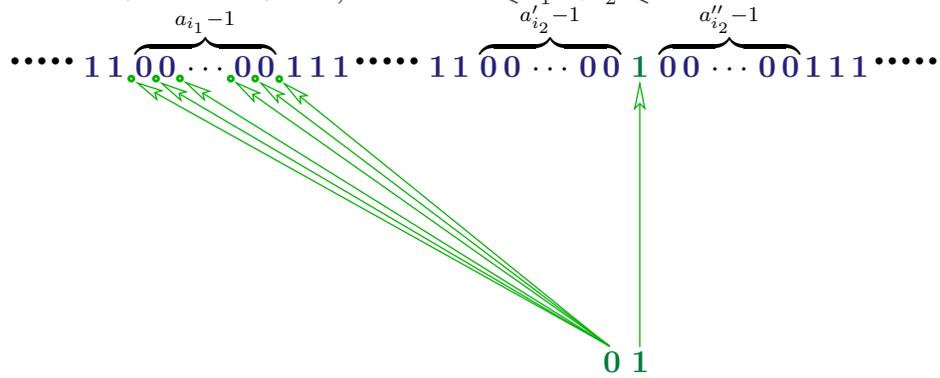
\bigskip
\begin{center}
\input 0-i1-1-i2.pstex_t
\end{center}

\noindent
The placement of the $0$ on the left always gives $a_{ i_1} + 1$ (with
multiplicity $a_{ i_1}$) instead of $a_{ i_1}$ in the original
polyzeta, while the placement of $1$ on the right imposes to replace
the $a_{ i_2}$ in the original polyzeta by $a_{ i_2}', a_{ i_2}''$
with:
\[
a_{i_2}'-1+a_{i_2}''-1
=
a_{i_2}-1,
\]
that is to say:
\[
a_{i_2}'+a_{i_2}''
=
a_{i_2}+1.
\]
As a result, the polyzetas obtained in this fifth family are:
\[
\boxed{
\mathcal{F}_{i_1\vert i_2}
:=
\sum_{1\leqslant i_1<i_2\leqslant{\sf h}\atop a_{i_2}\geqslant 3}
\Bigg(
\sum_{a_{i_2}'+a_{i_2}''=a_{i_2}+1
\atop
a_{i_2}'\geqslant 2,\,\,a_{i_2}''\geqslant 2}\,
a_{i_1}\cdot
\zeta_{[{\sf w},{\sf d}+1,{\sf h}+1]}
\big(
\bfdots,a_{i_1}+1,\bfdots,a_{i_2}',1^0,a_{i_2}'',\bfdots
\big)
\Bigg).}
\]
Visibly, their total number is:
\[
{\bf n}_{i_1\vert i_2}
:=
\sum_{1\leqslant i_2<i_2\leqslant{\sf h}}\,
a_{i_1}\,(a_{i_2}-2).
\]

\subsection{Sixth family $\boxed{j \vert i}$ 
with $1 \leqslant j < i \leqslant {\sf h}$\,} 
Finally, assume that the $0$ of $01$ goes into a block
$1\, 1^{ b_j}$ and that the $1$ of $01$ goes
into a block $0^{ a_i-1}$, with
$0 \leqslant j < i \leqslant {\sf h}$\,\,---\,\,the inequality
$j < i$ must indeed be strict\,\,---\,:

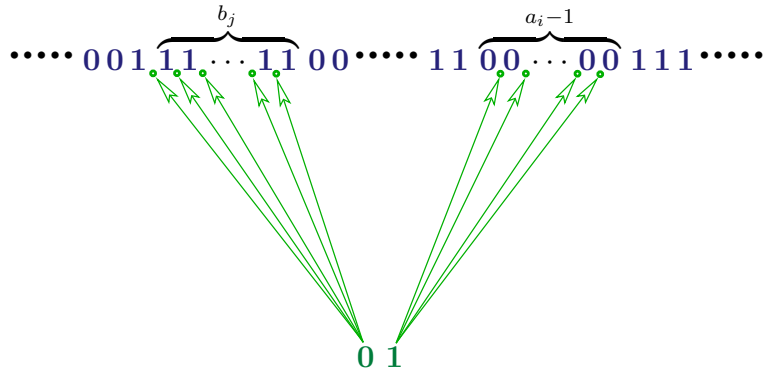
\bigskip
\begin{center}
\input 0-j-1-i.pstex_t
\end{center}

\noindent
The agglutination principles require that all insertions are strict.
No repetition appears, no multiplicity holds, all obtained polyzetas:
\[
\bfdots
\bbf{0}\,\bbf{0}\,
\unsqov{b_j'}\,
\gbf{0}\,
\unsqov{b_j''}\,
\bbf{0}\,\bbf{0}\,
\bfdots\,
\bbf{1}\,\bbf{1}\,
\zsqov{a_i'-1}\,
\gbf{1}\,
\zsqov{a_i''-1}\,
\bbf{1}\,\bbf{1}\,\bbf{1}\,
\bfdots.
\]
are pairwise distinct, and the appearing new integers must satisfy:
\[
b_j'+b_j''=b_j-1
\ \ \ \ \ \ \ \ \ \ \ \ \ \ \
\text{\rm and}
\ \ \ \ \ \ \ \ \ \ \ \ \ \ \
a_i'+a_i''
=
a_i+1.
\]
As a result, the polyzetas of this sixth and last family are
expressed by means of three summation symbols:
\[
\boxed{
\mathcal{F}_{j\vert i}
:=
\sum_{1\leqslant j<i\leqslant{\sf h}
\atop
b_j\geqslant 1,\,\,a_i\geqslant 3}
\Bigg(
\sum_{b_j'+b_j''=b_j-1
\atop
b_j'\geqslant 0,\,\,b_j''\geqslant 0}\,
\sum_{a_i'+a_i''=a_i+1
\atop
a_i'\geqslant 2,\,\,a_i''\geqslant 2}\,
\zeta_{[{\sf w},{\sf d}+1,{\sf h}+2]}
\big(
\bfdots,1^{b_j'},2,1^{b_j''},\bfdots,
a_i',1^0,a_i'',\bfdots
\big)
\Bigg).}
\]
Mentally, one sees that their total number is:
\[
{\bf n}_{j\vert i}
:=
\sum_{1\leqslant j<i\leqslant{\sf h}}\,
b_j\,(a_i-2).
\]

\subsection{Final counting checking}
\label{final-counting-checking}
In order to convince ourselves that the formulas are completely free
of errors, let us verify the promised equality concerning the number
of obtained polyzetas. On one hand, we may firstly 
expand\,\,---\,\,disregarding temporarily the underlinings
whose r\^ole will be explained in a moment\,\,---\,:
\[
\footnotesize
\aligned
&
\frac{\big(a_1+b_1+\cdots+a_{\sf h}+b_{\sf h}+2\big)
\big(a_1+b_1+\cdots+a_{\sf h}+b_{\sf h}+1\big)}{2}
=
\\
&
\ \ \ \ \
=
\underline{
\sum_{1\leqslant i\leqslant{\sf h}}\,
\frac{a_i\,a_i}{2}
}_{\octagon\!\!\!\!\!1}
+
\underline{
\sum_{1\leqslant j\leqslant{\sf h}}\,
\frac{b_j\,b_j}{2}
}_{\octagon\!\!\!\!\!2}
+
\underline{
\sum_{1\leqslant i\leqslant{\sf h}
\atop
1\leqslant j\leqslant{\sf h}}\,
a_i\,b_j
}_{\octagon\!\!\!\!\!3}
+
\underline{
\sum_{1\leqslant i_1<i_2\leqslant{\sf h}}\,
a_{i_1}\,a_{i_2}
}_{\octagon\!\!\!\!\!4}
+
\\
&
\ \ \ \ \ \ \ \ \ \ 
+
\underline{
\sum_{1\leqslant j_1<j_2\leqslant{\sf h}}\,
b_{j_1}\,b_{j_2}
}_{\octagon\!\!\!\!\!5}
+
\underline{
\frac{3}{2}\,
\sum_{1\leqslant i\leqslant{\sf h}}\,
a_i
}_{\octagon\!\!\!\!\!6}
+
\underline{
\frac{3}{2}\,
\sum_{1\leqslant j\leqslant{\sf h}}\,
b_j
}_{\octagon\!\!\!\!\!7}
+
\underline{
\,1\,
}_{\octagon\!\!\!\!\!8}\:.
\endaligned
\]
On the other hand, let us collect the total number of polyzetas we
obtained in our six families:
\[
\footnotesize
\aligned
{\bf n}_{i\vert j}
+
{\bf n}_{j_1\vert j_2}
+
{\bf n}_{j\,j}
+
{\bf n}_{i\,i}
+
{\bf n}_{i_1\vert i_2}
+
{\bf n}_{j\vert i}
&
=
\underline{
\sum_{1\leqslant i\leqslant j\leqslant{\sf h}}\,
a_i\,b_j
}_{\octagon\!\!\!\!\!3}
+
\underline{
\sum_{1\leqslant i\leqslant j\leqslant{\sf h}}\,
2\,a_i
}_{\octagon\!\!\!\!\!6}
+
\\
&
\ \ \ \ \
+
\underline{
\sum_{1\leqslant j_1<j_2\leqslant{\sf h}}\,
b_{j_1}\,b_{j_2}
}_{\octagon\!\!\!\!\!5}
+
\zero{\sum_{1\leqslant j_1<j_2\leqslant{\sf h}}\,
2\,b_{j_1}}
+
\\
&
\ \ \ \ \
+
\sum_{1\leqslant j\leqslant{\sf h}}
\bigg(
\underline{
\frac{b_i\,b_j}{2}
}_{\octagon\!\!\!\!\!2}
+
\underline{
\frac{3\,b_j}{2}
}_{\octagon\!\!\!\!\!7}
+
\zerozero{\,1\,}
\bigg)
+
\\
&
\ \ \ \ \
+
\sum_{1\leqslant i\leqslant{\sf h}}
\bigg(
\underline{
\frac{a_i\,a_i}{2}
}_{\octagon\!\!\!\!\!1}
-
\underline{
\frac{a_i}{2}
}_{\octagon\!\!\!\!\!6}
-
\zerozero{\,1\,}
\bigg)
+
\\
&
\ \ \ \ \
+
\underline{
\sum_{1\leqslant i_1<i_2\leqslant{\sf h}}\,
a_{i_1}\,a_{i_2}
}_{\octagon\!\!\!\!\!4}
-
\underline{
\sum_{1\leqslant i_1<i_2\leqslant{\sf h}}\,
2\,a_{i_2}
}_{\octagon\!\!\!\!\!6}
+
\\
&
\ \ \ \ \
+
\underline{
\sum_{1\leqslant j<i\leqslant{\sf h}}\,
a_i\,b_j
}_{\octagon\!\!\!\!\!3}
-
\zero{\sum_{1\leqslant j<i\leqslant{\sf h}}\,
2\,b_j}.
\endaligned
\]
Thus, we have to compare these two series of sums that both yield a
nonnegative integer and to verify that the two resulting
integers are precisely equal. To
do this in a way that does not `leave' some non-immediate painful
computations to a reader, let us borrow from~\cite{Merker-2006} the
{\sl technique of numbered underlinings}, the purpose of which is to
enable a {\em pure eye-checking of computations}, without the need of
any pencil or of any extra sheet of paper.

To begin with, in the second series of sums just above, let us observe
that two pairs of terms auto-annihilate, and they are denoted by means
of the two pairs of underlinings $\underline{\,\,\,}_\circ$ and
$\underline{\,\,\,}_{\circ\circ}$.

Next, with the symbol $\octagon$, we employ numbered underlinings to
point out the coincidences of sum-terms between the two expressions
under study. Notice that, especially for number
$\octagon\!\!\!\!{\scriptstyle{6}}$\,\,, the summation reduction:
\[
\sum_{1\leqslant i\leqslant j\leqslant{\sf h}}\,
-
\sum_{1\leqslant i_1<i_2\leqslant{\sf h}}\,
=
\sum_{1\leqslant i\leqslant{\sf h}}\,
\]
must be used, while a last mental addition concludes the coincidence.
Thus, number of terms match, as was required by this test of coherency.
\qed

\smallskip

All the preceding reasonings enabled us to gain a complete expansion
of the shuffle product of $\zeta ( 2)$ with an arbitrary polyzeta.

\begin{Proposition}
For every height ${\sf h} \geqslant 1$, every entries $a_1 \geqslant
2$, \dots, $a_{\sf h} \geqslant 2$, every entries $b_1 \geqslant 0$,
\dots, $b_{\sf h} \geqslant 0$, if one sets:
\[
\aligned
a_1+b_1+\cdots+a_{\sf h}+b_{\sf h}
&
=:
{\sf w}-2,
\\
1+b_1+\cdots+1+b_{\sf h}
&
=:
{\sf d},
\endaligned
\]
then the $\shuffle$-shuffle product of $\zeta (2)$ with
any $\zeta_{[{\sf w}-2,{\sf d},{\sf h}]}
\big( a_1, 1^{ b_1}, \dots, a_{\sf h}, 1^{ b_{\sf h}} \big)$
of weight ${\sf w} - 2 \geqslant 2$, of depth ${\sf d}$ and of height ${\sf h}$
writes out, after complete finalization:
\[
\aligned
\zeta(2)
\shuffle
&
\zeta_{[{\sf w}-2,{\sf d},{\sf h}]}
\big(
a_1,1^{b_1},\dots,a_{\sf h},1^{b_{\sf h}}
\big)
=
\mathcal{F}_{i\vert j}
+
\mathcal{F}_{j_1\vert j_2}
+
\mathcal{F}_{j\,j}
+
\mathcal{F}_{i\,i}
+
\mathcal{F}_{i_1\vert i_2}
+
\mathcal{F}_{j\vert i}
=
\\
&
=\sum_{1\leqslant i\leqslant j\leqslant{\sf h}}\,
a_i\,(b_j+2)\,
\zeta_{[{\sf w},{\sf d+1},{\sf h}]}
\big(
\bfdots,a_i+1,\bfdots,1^{b_j+1},\bfdots
\big)
+
\\
&
\ \ \ \ \
+
\sum_{1\leqslant j_1<j_2\leqslant{\sf h}\atop b_{j_1}\geqslant 1}
\Bigg(
\sum_{b_{j_1}'+b_{j_1}''=b_{j_1-1}
\atop
b_{j_1}'\geqslant 0,\,\,b_{j_1}''\geqslant 0}\,
(b_{j_2}+2)
\cdot
\zeta_{[{\sf w},{\sf d}+1,{\sf h}+1]}
\big(
\bfdots,1^{b_{j_1'}},2,1^{b_{j_1}''},\bfdots,
1^{b_{j_2+1}},\bfdots
\big)
\Bigg)
+
\\
&
\ \ \ \ \
+
\sum_{1\leqslant j\leqslant{\sf h}}
\Bigg(
\sum_{b_j'+b_j''=b_j
\atop
b_j'\geqslant 0,\,\,b_j''\geqslant 0}\,
(b_j''+1)\cdot
\zeta_{[{\sf w},{\sf d}+1,{\sf h}+1]}
\big(\bfdots,1^{b_j'},2,1^{b_j''},\bfdots
\big)
\Bigg)
+
\\
&
\ \ \ \ \
+
\big(
2,a_1,1^{b_1},\dots,a_{\sf h},1^{b_{\sf h}}
\big)
+
\\
&
\ \ \ \ \
+
\sum_{1\leqslant i\leqslant{\sf h}\atop a_i\geqslant 5}
\Bigg(
\sum_{a_i'+a_i''=a_i+2
\atop
a_i'\geqslant 3,\,\,a_i''\geqslant 2}\,
(a_i'-1)\cdot
\zeta_{[{\sf w},{\sf d}+1,{\sf h}+1]}
\big(
\bfdots,a_i',1^0,a_i',\bfdots
\big)
\Bigg)
+
\\
&
\ \ \ \ \
+
\sum_{1\leqslant i_1<i_2\leqslant{\sf h}\atop a_{i_2}\geqslant 3}
\Bigg(
\sum_{a_{i_2}'+a_{i_2}''=a_{i_2}+1
\atop
a_{i_2}'\geqslant 2,\,\,a_{i_2}''\geqslant 2}\,
a_{i_1}\cdot
\zeta_{[{\sf w},{\sf d}+1,{\sf h}+1]}
\big(
\bfdots,a_{i_1}+1,\bfdots,a_{i_2}',1^0,a_{i_2}'',\bfdots
\big)
\Bigg)
+
\\
&
\ \ \ \ \
+
\sum_{1\leqslant j<i\leqslant{\sf h}
\atop
b_j\geqslant 1,\,\,a_i\geqslant 3}
\Bigg(
\sum_{b_j'+b_j''=b_j-1
\atop
b_j'\geqslant 0,\,\,b_j''\geqslant 0}\,
\sum_{a_i'+a_i''=a_i+1
\atop
a_i'\geqslant 2,\,\,a_i''\geqslant 2}\,
\zeta_{[{\sf w},{\sf d}+1,{\sf h}+2]}
\big(
\bfdots,1^{b_j'},2,1^{b_j''},\bfdots,
a_i',1^0,a_i'',\bfdots
\big)
\Bigg).
\endaligned
\]
\end{Proposition}

\subsection{The subtraction} Now that we have gained the two fundamental
expressions of the stuffle and of the shuffle product of $\zeta ( 2)$
with an arbitrary polyzeta of weight ${\sf w} - 2$, we can execute the
final subtraction. Not much terms collapse, no hint is really needed,
hence we plainly (re)state the result.

\begin{Theorem}
For every height ${\sf h} \geqslant 1$,
every entries $a_1 \geqslant 2$, \dots, $a_{\sf h} \geqslant 2$,
every entries $b_1 \geqslant 0$, \dots, $b_{\sf h} \geqslant 0$,
if one sets:
\[
\aligned
a_1+b_1+\cdots+a_{\sf h}+b_{\sf h}
&
=:
{\sf w}-2,
\\
1+b_1+\cdots+1+b_{\sf h}
&
=:
{\sf d},
\endaligned
\]
then the so-called double shuffle relation: 
\[
0
=
-\,(2)
\ast
\big(a_1,1^{b_1},\dots,a_{\sf h},1^{b_{\sf h}}\big)
+
01
\shuffle
0^{a_1-1}1\,1^{b_1}\,\cdots\,0^{a_{\sf h}-1}1\,1^{b_{\sf h}}
\]
between $\zeta (2)$ and any $\zeta_{[{\sf w}-2,{\sf d},{\sf h}]}
\big( a_1, 1^{ b_1}, \dots, a_{\sf h}, 1^{ b_{\sf h}} \big)$ of weight
${\sf w} - 2 \geqslant 2$, of depth ${\sf d}$ 
and of height ${\sf h}$ writes out after complete
finalization:
\[
\boxed{
\small
\aligned
0
&
=
-\,\sum_{1\leqslant i\leqslant{\sf h}}\,
\zeta_{[{\sf w},{\sf d},{\sf h}]}
\big(\bfdots,a_i+2,\bfdots\big)
-
\\
&
\ \ \ \ \
-\,
\sum_{1\leqslant j\leqslant{\sf h}\atop b_j\geqslant 1}
\Bigg(
\sum_{b_j'+b_j''=b_j-1
\atop
b_j'\geqslant 0,\,\,b_j''\geqslant 0}\,
\zeta_{[{\sf w},{\sf d},{\sf h}+1]}
\big(\bfdots,1^{b_j'},3,1^{b_j''},\bfdots\big)
\Bigg)
+
\\
&
\ \ \ \ \
+
\sum_{1\leqslant i\leqslant j\leqslant{\sf h}}\,
a_i\,(b_j+2)\cdot
\zeta_{[{\sf w},{\sf d}+1,{\sf h}]}
\big(\bfdots,a_i+1,\bfdots,1^{b_j+1},\bfdots\big)
+
\\
&
\ \ \ \ \
+
\sum_{1\leqslant j_1<j_2\leqslant{\sf h}\atop b_{j_1}\geqslant 1}
\Bigg(
\sum_{b_{j_1}'+b_{j_1}''=b_{j_1}-1
\atop
b_{j_1}'\geqslant 0,\,\,b_{j_1}''\geqslant 0}\,
(b_{j_2}+2)\cdot
\zeta_{[{\sf w},{\sf d}+1,{\sf h}+1]}
\big(\bfdots,1^{b_{j_1}'},2,1^{b_{j_1}''},\bfdots,1^{b_{j_2}+1},\bfdots\big)
\Bigg)
+
\\
&
\ \ \ \ \
+
\sum_{1\leqslant j\leqslant{\sf h}}
\Bigg(
\sum_{b_j'+b_j''=b_j
\atop
b_j'\geqslant 0,\,\,b_j''\geqslant 0}\,
b_j''\cdot
\zeta_{[{\sf w},{\sf d}+1,{\sf h}+1]}
\big(\bfdots,1^{b_j'},2,1^{b_j''},\bfdots\big)
\Bigg)
+
\\
&
\ \ \ \ \
+
\sum_{1\leqslant i\leqslant{\sf h}\atop a_i\geqslant 5}
\Bigg(
\sum_{a_i'+a_i''=a_i+2
\atop
a_i'\geqslant 3,\,\,a_i''\geqslant 2}\,
(a_i'-1)\cdot
\zeta_{[{\sf w},{\sf d}+1,{\sf h}+1]}
\big(\bfdots,a_i',1^0,a_i'',\bfdots\big)
\Bigg)
+
\\
&
\ \ \ \ \
+
\sum_{1\leqslant i_1<i_2\leqslant{\sf h}\atop a_{i_2}\geqslant 3}\,
\Bigg(
\sum_{a_{i_2}'+a_{i_2}''=a_{i_2}+1
\atop
a_{i_2}'\geqslant 2,\,\,a_{i_2''}\geqslant 2}\,
a_{i_1}
\cdot
\zeta_{[{\sf w},{\sf d}+1,{\sf h}+1]}
\big(\bfdots,a_{i_1}+1,\bfdots,a_{i_2}',1^0,a_{i_2}'',\bfdots\big)
\Bigg)
+
\\
&
\ \ \ \ \
+
\sum_{1\leqslant j<i\leqslant{\sf h}
\atop
b_j\geqslant 1,\,\,a_i\geqslant 3}
\Bigg(
\sum_{b_j'+b_j''=b_j-1
\atop
b_j'\geqslant 0,\,\,b_j''\geqslant 0}\,
\sum_{a_i'+a_i''=a_i+1
\atop
a_i'\geqslant 2,\,\,a_i''\geqslant 2}\,
\zeta_{[{\sf w},{\sf d}+1,{\sf h}+2]}
\big(\bfdots,1^{b_j'},2,1^{b_j''},\bfdots,a_i',1^0,a_i'',\bfdots\big)
\Bigg).
\endaligned}
\]
\end{Theorem}

\section{Computing $-(3) \ast \big(a_1,1^{b_1}, \dots,a_{\sf h},
1^{b_{\sf h}}\big) 
+ 001 \shuffle 0^{ a_1 -1}1\, 1^{b_1} \cdots
0^{a_{\sf h}-1}1\, 1^{b_{\sf h}}$}
\label{3-double-melange}
\HEAD{Computing $-(3) \ast \big(a_1,1^{b_1}, \dots,a_{\sf h},
1^{b_{\sf h}}\big) 
+ 001 \shuffle 0^{ a_1 -1}1\, 1^{b_1} \cdots
0^{a_{\sf h}-1}1\, 1^{b_{\sf h}}$}{
Jo\"el {\sc Merker}, Universit\'e Paris-Sud Orsay, France}

The principles are similar. As an appendix to this prepublication,
a detailed proof is expected to appear on {\footnotesize\sf
arxiv.org}.

\section{Computing $-(2,1) \ast \big(a_1,1^{b_1}, \dots,a_{\sf h},
1^{b_{\sf h}}\big) 
+ 011 \shuffle 0^{ a_1 -1}1\, 1^{b_1} \cdots
0^{a_{\sf h}-1}1\, 1^{b_{\sf h}}$}
\label{21-double-melange}
\HEAD{Computing $-(2,1) \ast \big(a_1,1^{b_1}, \dots,a_{\sf h},
1^{b_{\sf h}}\big) + 011 \shuffle 0^{ a_1 -1}1\, 1^{b_1} \cdots
0^{a_{\sf h}-1}1\, 1^{b_{\sf h}}$}{
Jo\"el {\sc Merker}, Universit\'e Paris-Sud Orsay, France}

See in the future 
the same appendix to this prepublication, which might be
twenty pages long. In fact, the point is not
just to set up these formulas (proofs
may then be forgotten), but mainly to {\em apply} them, 
to {\em handle} them, 
to {\em compute} with them.

\vfill\end{document}

%% file: 1-all-of.pstex_t
\begin{picture}(0,0)%
\includegraphics{1-all-of.pstex}%
\end{picture}%
\setlength{\unitlength}{4144sp}%
\begingroup\makeatletter\ifx\SetFigFont\undefined
\def\x#1#2#3#4#5#6#7\relax{\def\x{#1#2#3#4#5#6}}%
\expandafter\x\fmtname xxxxxx\relax \def\y{splain}%
\ifx\x\y   
\gdef\SetFigFont#1#2#3{%
  \ifnum #1<17\tiny\else \ifnum #1<20\small\else
  \ifnum #1<24\normalsize\else \ifnum #1<29\large\else
  \ifnum #1<34\Large\else \ifnum #1<41\LARGE\else
     \huge\fi\fi\fi\fi\fi\fi
  \csname #3\endcsname}%
\else
\gdef\SetFigFont#1#2#3{\begingroup
  \count@#1\relax \ifnum 25<\count@\count@25\fi
  \def\x{\endgroup\@setsize\SetFigFont{#2pt}}%
  \expandafter\x
    \csname \romannumeral\the\count@ pt\expandafter\endcsname
    \csname @\romannumeral\the\count@ pt\endcsname
  \csname #3\endcsname}%
\fi
\fi\endgroup
\begin{picture}(2373,1182)(1313,-1634)
\put(1352,-632){\makebox(0,0)[lb]{\smash{\SetFigFont{14}{16.8}{rm}{\color[rgb]{0,0,0}$\zsqov{a_1-1}\,\unsqov{b_1}$}%
}}}
\put(2434,-1571){\makebox(0,0)[lb]{\smash{\SetFigFont{14}{16.8}{rm}{\color[rgb]{0,0,0}$\gbf{1}$}%
}}}
\end{picture}

%% file: strictly-inside-0.pstex_t
\begin{picture}(0,0)%
\includegraphics{strictly-inside-0.pstex}%
\end{picture}%
\setlength{\unitlength}{4144sp}%
\begingroup\makeatletter\ifx\SetFigFont\undefined
\def\x#1#2#3#4#5#6#7\relax{\def\x{#1#2#3#4#5#6}}%
\expandafter\x\fmtname xxxxxx\relax \def\y{splain}%
\ifx\x\y   
\gdef\SetFigFont#1#2#3{%
  \ifnum #1<17\tiny\else \ifnum #1<20\small\else
  \ifnum #1<24\normalsize\else \ifnum #1<29\large\else
  \ifnum #1<34\Large\else \ifnum #1<41\LARGE\else
     \huge\fi\fi\fi\fi\fi\fi
  \csname #3\endcsname}%
\else
\gdef\SetFigFont#1#2#3{\begingroup
  \count@#1\relax \ifnum 25<\count@\count@25\fi
  \def\x{\endgroup\@setsize\SetFigFont{#2pt}}%
  \expandafter\x
    \csname \romannumeral\the\count@ pt\expandafter\endcsname
    \csname @\romannumeral\the\count@ pt\endcsname
  \csname #3\endcsname}%
\fi
\fi\endgroup
\begin{picture}(1188,1472)(908,-2892)
\put(908,-1613){\makebox(0,0)[lb]{\smash{\SetFigFont{14}{16.8}{rm}{\color[rgb]{0,0,0}$\zsqov{a_1'-1}\,\gbf{1}\,\zsqov{a_1''-1}\,\unsqov{b_1}$}%
}}}
\put(1988,-2824){\makebox(0,0)[lb]{\smash{\SetFigFont{14}{16.8}{rm}{\color[rgb]{0,0,0}$\gbf{1}$}%
}}}
\end{picture}

%% file: 1-agglutination-1.pstex_t
\begin{picture}(0,0)%
\includegraphics{1-agglutination-1.pstex}%
\end{picture}%
\setlength{\unitlength}{4144sp}%
\begingroup\makeatletter\ifx\SetFigFont\undefined
\def\x#1#2#3#4#5#6#7\relax{\def\x{#1#2#3#4#5#6}}%
\expandafter\x\fmtname xxxxxx\relax \def\y{splain}%
\ifx\x\y   
\gdef\SetFigFont#1#2#3{%
  \ifnum #1<17\tiny\else \ifnum #1<20\small\else
  \ifnum #1<24\normalsize\else \ifnum #1<29\large\else
  \ifnum #1<34\Large\else \ifnum #1<41\LARGE\else
     \huge\fi\fi\fi\fi\fi\fi
  \csname #3\endcsname}%
\else
\gdef\SetFigFont#1#2#3{\begingroup
  \count@#1\relax \ifnum 25<\count@\count@25\fi
  \def\x{\endgroup\@setsize\SetFigFont{#2pt}}%
  \expandafter\x
    \csname \romannumeral\the\count@ pt\expandafter\endcsname
    \csname @\romannumeral\the\count@ pt\endcsname
  \csname #3\endcsname}%
\fi
\fi\endgroup
\begin{picture}(2056,1332)(1184,-2821)
\put(2608,-2758){\makebox(0,0)[lb]{\smash{\SetFigFont{14}{16.8}{rm}{\color[rgb]{0,0,0}$\gbf{1}$}%
}}}
\put(1184,-1669){\makebox(0,0)[lb]{\smash{\SetFigFont{14}{16.8}{rm}{\color[rgb]{0,0,0}$\upbfdots\,\bbf{0}\,\bbf{0}\,\unsqov{b_j}\,\bbf{0}\,\bbf{0}\,\upbfdots$}%
}}}
\end{picture}

%% file: re-strictly-inside-0.pstex_t
\begin{picture}(0,0)%
\includegraphics{re-strictly-inside-0.pstex}%
\end{picture}%
\setlength{\unitlength}{4144sp}%
\begingroup\makeatletter\ifx\SetFigFont\undefined
\def\x#1#2#3#4#5#6#7\relax{\def\x{#1#2#3#4#5#6}}%
\expandafter\x\fmtname xxxxxx\relax \def\y{splain}%
\ifx\x\y   
\gdef\SetFigFont#1#2#3{%
  \ifnum #1<17\tiny\else \ifnum #1<20\small\else
  \ifnum #1<24\normalsize\else \ifnum #1<29\large\else
  \ifnum #1<34\Large\else \ifnum #1<41\LARGE\else
     \huge\fi\fi\fi\fi\fi\fi
  \csname #3\endcsname}%
\else
\gdef\SetFigFont#1#2#3{\begingroup
  \count@#1\relax \ifnum 25<\count@\count@25\fi
  \def\x{\endgroup\@setsize\SetFigFont{#2pt}}%
  \expandafter\x
    \csname \romannumeral\the\count@ pt\expandafter\endcsname
    \csname @\romannumeral\the\count@ pt\endcsname
  \csname #3\endcsname}%
\fi
\fi\endgroup
\begin{picture}(1733,1396)(62,-2881)
\put(1340,-2818){\makebox(0,0)[lb]{\smash{\SetFigFont{14}{16.8}{rm}{\color[rgb]{0,0,0}$\gbf{1}$}%
}}}
\put( 62,-1665){\makebox(0,0)[lb]{\smash{\SetFigFont{14}{16.8}{rm}{\color[rgb]{0,0,0}$\upbfdots\,\bbf{1}\,\bbf{1}\,\zsqov{a_i-1}\,\bbf{1}\,\bbf{1}\,\bbf{1}\,\upbfdots$}%
}}}
\end{picture}

%% file: 1-a1-b1.pstex_t
\begin{picture}(0,0)%
\includegraphics{1-a1-b1.pstex}%
\end{picture}%
\setlength{\unitlength}{4144sp}%
\begingroup\makeatletter\ifx\SetFigFont\undefined
\def\x#1#2#3#4#5#6#7\relax{\def\x{#1#2#3#4#5#6}}%
\expandafter\x\fmtname xxxxxx\relax \def\y{splain}%
\ifx\x\y   
\gdef\SetFigFont#1#2#3{%
  \ifnum #1<17\tiny\else \ifnum #1<20\small\else
  \ifnum #1<24\normalsize\else \ifnum #1<29\large\else
  \ifnum #1<34\Large\else \ifnum #1<41\LARGE\else
     \huge\fi\fi\fi\fi\fi\fi
  \csname #3\endcsname}%
\else
\gdef\SetFigFont#1#2#3{\begingroup
  \count@#1\relax \ifnum 25<\count@\count@25\fi
  \def\x{\endgroup\@setsize\SetFigFont{#2pt}}%
  \expandafter\x
    \csname \romannumeral\the\count@ pt\expandafter\endcsname
    \csname @\romannumeral\the\count@ pt\endcsname
  \csname #3\endcsname}%
\fi
\fi\endgroup
\begin{picture}(2332,1335)(908,-2821)
\put(908,-1666){\makebox(0,0)[lb]{\smash{\SetFigFont{14}{16.8}{rm}{\color[rgb]{0,0,0}$\zsqov{a_1-1}\,\unsqov{b_1}$}%
}}}
\put(2608,-2758){\makebox(0,0)[lb]{\smash{\SetFigFont{14}{16.8}{rm}{\color[rgb]{0,0,0}$\gbf{1}$}%
}}}
\end{picture}

%% file: 1-ins-1-in-00-00.pstex_t
\begin{picture}(0,0)%
\includegraphics{1-ins-1-in-00-00.pstex}%
\end{picture}%
\setlength{\unitlength}{4144sp}%
\begingroup\makeatletter\ifx\SetFigFont\undefined
\def\x#1#2#3#4#5#6#7\relax{\def\x{#1#2#3#4#5#6}}%
\expandafter\x\fmtname xxxxxx\relax \def\y{splain}%
\ifx\x\y   
\gdef\SetFigFont#1#2#3{%
  \ifnum #1<17\tiny\else \ifnum #1<20\small\else
  \ifnum #1<24\normalsize\else \ifnum #1<29\large\else
  \ifnum #1<34\Large\else \ifnum #1<41\LARGE\else
     \huge\fi\fi\fi\fi\fi\fi
  \csname #3\endcsname}%
\else
\gdef\SetFigFont#1#2#3{\begingroup
  \count@#1\relax \ifnum 25<\count@\count@25\fi
  \def\x{\endgroup\@setsize\SetFigFont{#2pt}}%
  \expandafter\x
    \csname \romannumeral\the\count@ pt\expandafter\endcsname
    \csname @\romannumeral\the\count@ pt\endcsname
  \csname #3\endcsname}%
\fi
\fi\endgroup
\begin{picture}(5191,1568)(626,-1317)
\put(626, 95){\makebox(0,0)[lb]{\smash{\SetFigFont{12}{14.4}{rm}{\color[rgb]{0,0,0}$\zsqov{a_1-1}\,\unsqov{b_1}\,\,\zsqov{a_2-1}\,\unsqov{b_2}\,\upbfdots\,\zsqov{a_{\sf h}-1}\,\unsqov{b_{\sf h}}$}%
}}}
\put(3466,-1259){\makebox(0,0)[lb]{\smash{\SetFigFont{12}{14.4}{rm}{\color[rgb]{0,0,0}$\gbf{1}$}%
}}}
\end{picture}

%% file: 1-ins-1-in-111-11.pstex_t
\begin{picture}(0,0)%
\includegraphics{1-ins-1-in-111-11.pstex}%
\end{picture}%
\setlength{\unitlength}{4144sp}%
\begingroup\makeatletter\ifx\SetFigFont\undefined
\def\x#1#2#3#4#5#6#7\relax{\def\x{#1#2#3#4#5#6}}%
\expandafter\x\fmtname xxxxxx\relax \def\y{splain}%
\ifx\x\y   
\gdef\SetFigFont#1#2#3{%
  \ifnum #1<17\tiny\else \ifnum #1<20\small\else
  \ifnum #1<24\normalsize\else \ifnum #1<29\large\else
  \ifnum #1<34\Large\else \ifnum #1<41\LARGE\else
     \huge\fi\fi\fi\fi\fi\fi
  \csname #3\endcsname}%
\else
\gdef\SetFigFont#1#2#3{\begingroup
  \count@#1\relax \ifnum 25<\count@\count@25\fi
  \def\x{\endgroup\@setsize\SetFigFont{#2pt}}%
  \expandafter\x
    \csname \romannumeral\the\count@ pt\expandafter\endcsname
    \csname @\romannumeral\the\count@ pt\endcsname
  \csname #3\endcsname}%
\fi
\fi\endgroup
\begin{picture}(6395,1568)(626,-1317)
\put(626, 95){\makebox(0,0)[lb]{\smash{\SetFigFont{12}{14.4}{rm}{\color[rgb]{0,0,0}$\zsqov{a_1-1}\,\unsqov{b_1}\,\,\zsqov{a_2-1}\,\unsqov{b_2}\,\upbfdots\,\zsqov{a_{\sf h}-1}\,\unsqov{b_{\sf h}}$}%
}}}
\put(3466,-1259){\makebox(0,0)[lb]{\smash{\SetFigFont{12}{14.4}{rm}{\color[rgb]{0,0,0}$\gbf{1}$}%
}}}
\end{picture}

%% file: 0-agglutination-0.pstex_t
\begin{picture}(0,0)%
\includegraphics{0-agglutination-0.pstex}%
\end{picture}%
\setlength{\unitlength}{4144sp}%
\begingroup\makeatletter\ifx\SetFigFont\undefined
\def\x#1#2#3#4#5#6#7\relax{\def\x{#1#2#3#4#5#6}}%
\expandafter\x\fmtname xxxxxx\relax \def\y{splain}%
\ifx\x\y   
\gdef\SetFigFont#1#2#3{%
  \ifnum #1<17\tiny\else \ifnum #1<20\small\else
  \ifnum #1<24\normalsize\else \ifnum #1<29\large\else
  \ifnum #1<34\Large\else \ifnum #1<41\LARGE\else
     \huge\fi\fi\fi\fi\fi\fi
  \csname #3\endcsname}%
\else
\gdef\SetFigFont#1#2#3{\begingroup
  \count@#1\relax \ifnum 25<\count@\count@25\fi
  \def\x{\endgroup\@setsize\SetFigFont{#2pt}}%
  \expandafter\x
    \csname \romannumeral\the\count@ pt\expandafter\endcsname
    \csname @\romannumeral\the\count@ pt\endcsname
  \csname #3\endcsname}%
\fi
\fi\endgroup
\begin{picture}(1893,1398)(1194,-2889)
\put(2671,-2826){\makebox(0,0)[lb]{\smash{\SetFigFont{14}{16.8}{rm}{\color[rgb]{0,0,0}$\gbf{0}$}%
}}}
\put(1194,-1671){\makebox(0,0)[lb]{\smash{\SetFigFont{14}{16.8}{rm}{\color[rgb]{0,0,0}$\upbfdots\,\bbf{1}\,\bbf{1}\,\zsqov{a_i-1}\,\bbf{1}\,\bbf{1}\,\upbfdots$}%
}}}
\end{picture}

%% file: re-strictly-inside-1.pstex_t
\begin{picture}(0,0)%
\includegraphics{re-strictly-inside-1.pstex}%
\end{picture}%
\setlength{\unitlength}{4144sp}%
\begingroup\makeatletter\ifx\SetFigFont\undefined
\def\x#1#2#3#4#5#6#7\relax{\def\x{#1#2#3#4#5#6}}%
\expandafter\x\fmtname xxxxxx\relax \def\y{splain}%
\ifx\x\y   
\gdef\SetFigFont#1#2#3{%
  \ifnum #1<17\tiny\else \ifnum #1<20\small\else
  \ifnum #1<24\normalsize\else \ifnum #1<29\large\else
  \ifnum #1<34\Large\else \ifnum #1<41\LARGE\else
     \huge\fi\fi\fi\fi\fi\fi
  \csname #3\endcsname}%
\else
\gdef\SetFigFont#1#2#3{\begingroup
  \count@#1\relax \ifnum 25<\count@\count@25\fi
  \def\x{\endgroup\@setsize\SetFigFont{#2pt}}%
  \expandafter\x
    \csname \romannumeral\the\count@ pt\expandafter\endcsname
    \csname @\romannumeral\the\count@ pt\endcsname
  \csname #3\endcsname}%
\fi
\fi\endgroup
\begin{picture}(1893,1398)(1194,-2889)
\put(2671,-2826){\makebox(0,0)[lb]{\smash{\SetFigFont{14}{16.8}{rm}{\color[rgb]{0,0,0}$\gbf{0}$}%
}}}
\put(1194,-1671){\makebox(0,0)[lb]{\smash{\SetFigFont{14}{16.8}{rm}{\color[rgb]{0,0,0}$\upbfdots\,\bbf{0}\,\bbf{0}\,\unsqov{b_j}\,\bbf{0}\,\bbf{0}\,\upbfdots$}%
}}}
\end{picture}

%% file: 0-goes-in-the-00-00.pstex_t
\begin{picture}(0,0)%
\includegraphics{0-goes-in-the-00-00.pstex}%
\end{picture}%
\setlength{\unitlength}{4144sp}%
\begingroup\makeatletter\ifx\SetFigFont\undefined
\def\x#1#2#3#4#5#6#7\relax{\def\x{#1#2#3#4#5#6}}%
\expandafter\x\fmtname xxxxxx\relax \def\y{splain}%
\ifx\x\y   
\gdef\SetFigFont#1#2#3{%
  \ifnum #1<17\tiny\else \ifnum #1<20\small\else
  \ifnum #1<24\normalsize\else \ifnum #1<29\large\else
  \ifnum #1<34\Large\else \ifnum #1<41\LARGE\else
     \huge\fi\fi\fi\fi\fi\fi
  \csname #3\endcsname}%
\else
\gdef\SetFigFont#1#2#3{\begingroup
  \count@#1\relax \ifnum 25<\count@\count@25\fi
  \def\x{\endgroup\@setsize\SetFigFont{#2pt}}%
  \expandafter\x
    \csname \romannumeral\the\count@ pt\expandafter\endcsname
    \csname @\romannumeral\the\count@ pt\endcsname
  \csname #3\endcsname}%
\fi
\fi\endgroup
\begin{picture}(4558,1299)(720,-1624)
\put(720,-481){\makebox(0,0)[lb]{\smash{\SetFigFont{12}{14.4}{rm}{\color[rgb]{0,0,0}$\zsqov{a_1-1}\,\unsqov{b_1}\,\upbfdots\,\zsqov{a_i-1}\,\upbfdots\,\zsqov{a_{\sf h}-1}\,\unsqov{b_{\sf h}}$}%
}}}
\put(3213,-1566){\makebox(0,0)[lb]{\smash{\SetFigFont{12}{14.4}{rm}{\color[rgb]{0,0,0}$\gbf{0}$}%
}}}
\end{picture}

%% file: 0-goes-in-the-111-11.pstex_t
\begin{picture}(0,0)%
\includegraphics{0-goes-in-the-111-11.pstex}%
\end{picture}%
\setlength{\unitlength}{4144sp}%
\begingroup\makeatletter\ifx\SetFigFont\undefined
\def\x#1#2#3#4#5#6#7\relax{\def\x{#1#2#3#4#5#6}}%
\expandafter\x\fmtname xxxxxx\relax \def\y{splain}%
\ifx\x\y   
\gdef\SetFigFont#1#2#3{%
  \ifnum #1<17\tiny\else \ifnum #1<20\small\else
  \ifnum #1<24\normalsize\else \ifnum #1<29\large\else
  \ifnum #1<34\Large\else \ifnum #1<41\LARGE\else
     \huge\fi\fi\fi\fi\fi\fi
  \csname #3\endcsname}%
\else
\gdef\SetFigFont#1#2#3{\begingroup
  \count@#1\relax \ifnum 25<\count@\count@25\fi
  \def\x{\endgroup\@setsize\SetFigFont{#2pt}}%
  \expandafter\x
    \csname \romannumeral\the\count@ pt\expandafter\endcsname
    \csname @\romannumeral\the\count@ pt\endcsname
  \csname #3\endcsname}%
\fi
\fi\endgroup
\begin{picture}(5770,1308)(720,-1633)
\put(720,-481){\makebox(0,0)[lb]{\smash{\SetFigFont{12}{14.4}{rm}{\color[rgb]{0,0,0}$\zsqov{a_1-1}\,\unsqov{b_1}\,\upbfdots\,\unsqov{b_j}\,\upbfdots\,\zsqov{a_{\sf h}-1}\,\unsqov{b_{\sf h}}$}%
}}}
\put(3180,-1575){\makebox(0,0)[lb]{\smash{\SetFigFont{12}{14.4}{rm}{\color[rgb]{0,0,0}$\gbf{0}$}%
}}}
\end{picture}

%% file: 1-goes-in-the-00-00.pstex_t
\begin{picture}(0,0)%
\includegraphics{1-goes-in-the-00-00.pstex}%
\end{picture}%
\setlength{\unitlength}{4144sp}%
\begingroup\makeatletter\ifx\SetFigFont\undefined
\def\x#1#2#3#4#5#6#7\relax{\def\x{#1#2#3#4#5#6}}%
\expandafter\x\fmtname xxxxxx\relax \def\y{splain}%
\ifx\x\y   
\gdef\SetFigFont#1#2#3{%
  \ifnum #1<17\tiny\else \ifnum #1<20\small\else
  \ifnum #1<24\normalsize\else \ifnum #1<29\large\else
  \ifnum #1<34\Large\else \ifnum #1<41\LARGE\else
     \huge\fi\fi\fi\fi\fi\fi
  \csname #3\endcsname}%
\else
\gdef\SetFigFont#1#2#3{\begingroup
  \count@#1\relax \ifnum 25<\count@\count@25\fi
  \def\x{\endgroup\@setsize\SetFigFont{#2pt}}%
  \expandafter\x
    \csname \romannumeral\the\count@ pt\expandafter\endcsname
    \csname @\romannumeral\the\count@ pt\endcsname
  \csname #3\endcsname}%
\fi
\fi\endgroup
\begin{picture}(4558,1299)(720,-1624)
\put(720,-481){\makebox(0,0)[lb]{\smash{\SetFigFont{12}{14.4}{rm}{\color[rgb]{0,0,0}$\zsqov{a_1-1}\,\unsqov{b_1}\,\upbfdots\,\zsqov{a_i-1}\,\upbfdots\,\zsqov{a_{\sf h}-1}\,\unsqov{b_{\sf h}}$}%
}}}
\put(3213,-1566){\makebox(0,0)[lb]{\smash{\SetFigFont{12}{14.4}{rm}{\color[rgb]{0,0,0}$\gbf{1}$}%
}}}
\end{picture}

%% file: 1-goes-in-the-111-11.pstex_t
\begin{picture}(0,0)%
\includegraphics{1-goes-in-the-111-11.pstex}%
\end{picture}%
\setlength{\unitlength}{4144sp}%
\begingroup\makeatletter\ifx\SetFigFont\undefined
\def\x#1#2#3#4#5#6#7\relax{\def\x{#1#2#3#4#5#6}}%
\expandafter\x\fmtname xxxxxx\relax \def\y{splain}%
\ifx\x\y   
\gdef\SetFigFont#1#2#3{%
  \ifnum #1<17\tiny\else \ifnum #1<20\small\else
  \ifnum #1<24\normalsize\else \ifnum #1<29\large\else
  \ifnum #1<34\Large\else \ifnum #1<41\LARGE\else
     \huge\fi\fi\fi\fi\fi\fi
  \csname #3\endcsname}%
\else
\gdef\SetFigFont#1#2#3{\begingroup
  \count@#1\relax \ifnum 25<\count@\count@25\fi
  \def\x{\endgroup\@setsize\SetFigFont{#2pt}}%
  \expandafter\x
    \csname \romannumeral\the\count@ pt\expandafter\endcsname
    \csname @\romannumeral\the\count@ pt\endcsname
  \csname #3\endcsname}%
\fi
\fi\endgroup
\begin{picture}(5770,1308)(720,-1633)
\put(720,-481){\makebox(0,0)[lb]{\smash{\SetFigFont{12}{14.4}{rm}{\color[rgb]{0,0,0}$\zsqov{a_1-1}\,\unsqov{b_1}\,\upbfdots\,\unsqov{b_j}\,\upbfdots\,\zsqov{a_{\sf h}-1}\,\unsqov{b_{\sf h}}$}%
}}}
\put(3180,-1575){\makebox(0,0)[lb]{\smash{\SetFigFont{12}{14.4}{rm}{\color[rgb]{0,0,0}$\gbf{1}$}%
}}}
\end{picture}

%% file: 0-goes-00-1-goes-111.pstex_t
\begin{picture}(0,0)%
\includegraphics{0-goes-00-1-goes-111.pstex}%
\end{picture}%
\setlength{\unitlength}{4144sp}%
\begingroup\makeatletter\ifx\SetFigFont\undefined
\def\x#1#2#3#4#5#6#7\relax{\def\x{#1#2#3#4#5#6}}%
\expandafter\x\fmtname xxxxxx\relax \def\y{splain}%
\ifx\x\y   
\gdef\SetFigFont#1#2#3{%
  \ifnum #1<17\tiny\else \ifnum #1<20\small\else
  \ifnum #1<24\normalsize\else \ifnum #1<29\large\else
  \ifnum #1<34\Large\else \ifnum #1<41\LARGE\else
     \huge\fi\fi\fi\fi\fi\fi
  \csname #3\endcsname}%
\else
\gdef\SetFigFont#1#2#3{\begingroup
  \count@#1\relax \ifnum 25<\count@\count@25\fi
  \def\x{\endgroup\@setsize\SetFigFont{#2pt}}%
  \expandafter\x
    \csname \romannumeral\the\count@ pt\expandafter\endcsname
    \csname @\romannumeral\the\count@ pt\endcsname
  \csname #3\endcsname}%
\fi
\fi\endgroup
\begin{picture}(3324,1444)(973,-1436)
\put(2620,-1369){\makebox(0,0)[lb]{\smash{\SetFigFont{14}{16.8}{rm}{\color[rgb]{0,0,0}$\gbf{0}\,\,\gbf{1}$}%
}}}
\put(973,-168){\makebox(0,0)[lb]{\smash{\SetFigFont{14}{16.8}{rm}{\color[rgb]{0,0,0}$\upbfdots\,\zsqov{a_i-1}\,\upbfdots\,\unsqov{b_j}\,\upbfdots$}%
}}}
\end{picture}

%% file: 0-j1-1-j2.pstex_t
\begin{picture}(0,0)%
\includegraphics{0-j1-1-j2.pstex}%
\end{picture}%
\setlength{\unitlength}{4144sp}%
\begingroup\makeatletter\ifx\SetFigFont\undefined
\def\x#1#2#3#4#5#6#7\relax{\def\x{#1#2#3#4#5#6}}%
\expandafter\x\fmtname xxxxxx\relax \def\y{splain}%
\ifx\x\y   
\gdef\SetFigFont#1#2#3{%
  \ifnum #1<17\tiny\else \ifnum #1<20\small\else
  \ifnum #1<24\normalsize\else \ifnum #1<29\large\else
  \ifnum #1<34\Large\else \ifnum #1<41\LARGE\else
     \huge\fi\fi\fi\fi\fi\fi
  \csname #3\endcsname}%
\else
\gdef\SetFigFont#1#2#3{\begingroup
  \count@#1\relax \ifnum 25<\count@\count@25\fi
  \def\x{\endgroup\@setsize\SetFigFont{#2pt}}%
  \expandafter\x
    \csname \romannumeral\the\count@ pt\expandafter\endcsname
    \csname @\romannumeral\the\count@ pt\endcsname
  \csname #3\endcsname}%
\fi
\fi\endgroup
\begin{picture}(4768,1577)(601,-1377)
\put(601, 44){\makebox(0,0)[lb]{\smash{\SetFigFont{12}{14.4}{rm}{\color[rgb]{0,0,0}$\upbfdots\,\unsqov{b_{j_1}'}\,\gbf{0}\,\unsqov{b_{j_1}''}\,\bbf{0}\,\bbf{0}\,\upbfdots\,\bbf{0}\,\bbf{0}\,\unsqov{b_{j_2}}\,\bbf{0}\,\bbf{0}\,\upbfdots$}%
}}}
\put(2111,-1319){\makebox(0,0)[lb]{\smash{\SetFigFont{12}{14.4}{rm}{\color[rgb]{0,0,0}$\gbf{0}$}%
}}}
\put(4904,-1316){\makebox(0,0)[lb]{\smash{\SetFigFont{12}{14.4}{rm}{\color[rgb]{0,0,0}$\gbf{1}$}%
}}}
\end{picture}

%% file: 0-j-1-j.pstex_t
\begin{picture}(0,0)%
\includegraphics{0-j-1-j.pstex}%
\end{picture}%
\setlength{\unitlength}{4144sp}%
\begingroup\makeatletter\ifx\SetFigFont\undefined
\def\x#1#2#3#4#5#6#7\relax{\def\x{#1#2#3#4#5#6}}%
\expandafter\x\fmtname xxxxxx\relax \def\y{splain}%
\ifx\x\y   
\gdef\SetFigFont#1#2#3{%
  \ifnum #1<17\tiny\else \ifnum #1<20\small\else
  \ifnum #1<24\normalsize\else \ifnum #1<29\large\else
  \ifnum #1<34\Large\else \ifnum #1<41\LARGE\else
     \huge\fi\fi\fi\fi\fi\fi
  \csname #3\endcsname}%
\else
\gdef\SetFigFont#1#2#3{\begingroup
  \count@#1\relax \ifnum 25<\count@\count@25\fi
  \def\x{\endgroup\@setsize\SetFigFont{#2pt}}%
  \expandafter\x
    \csname \romannumeral\the\count@ pt\expandafter\endcsname
    \csname @\romannumeral\the\count@ pt\endcsname
  \csname #3\endcsname}%
\fi
\fi\endgroup
\begin{picture}(2987,1626)(327,-1730)
\put(327,-260){\makebox(0,0)[lb]{\smash{\SetFigFont{12}{14.4}{rm}{\color[rgb]{0,0,0}$\upbfdots\,\bbf{0}\,\bbf{0}\,\unsqov{b_j'}\,\gbf{0}\,\unsqov{\widetilde{b_j}}\,\bbf{0}\,\bbf{0}\,\upbfdots$}%
}}}
\put(2111,-1672){\makebox(0,0)[lb]{\smash{\SetFigFont{12}{14.4}{rm}{\color[rgb]{0,.69,0}$\gbf{0\,\,1}$}%
}}}
\end{picture}

%% file: 0-i-1-i.pstex_t
\begin{picture}(0,0)%
\includegraphics{0-i-1-i.pstex}%
\end{picture}%
\setlength{\unitlength}{4144sp}%
\begingroup\makeatletter\ifx\SetFigFont\undefined
\def\x#1#2#3#4#5#6#7\relax{\def\x{#1#2#3#4#5#6}}%
\expandafter\x\fmtname xxxxxx\relax \def\y{splain}%
\ifx\x\y   
\gdef\SetFigFont#1#2#3{%
  \ifnum #1<17\tiny\else \ifnum #1<20\small\else
  \ifnum #1<24\normalsize\else \ifnum #1<29\large\else
  \ifnum #1<34\Large\else \ifnum #1<41\LARGE\else
     \huge\fi\fi\fi\fi\fi\fi
  \csname #3\endcsname}%
\else
\gdef\SetFigFont#1#2#3{\begingroup
  \count@#1\relax \ifnum 25<\count@\count@25\fi
  \def\x{\endgroup\@setsize\SetFigFont{#2pt}}%
  \expandafter\x
    \csname \romannumeral\the\count@ pt\expandafter\endcsname
    \csname @\romannumeral\the\count@ pt\endcsname
  \csname #3\endcsname}%
\fi
\fi\endgroup
\begin{picture}(1744,1622)(327,-1726)
\put(1798,-1668){\makebox(0,0)[lb]{\smash{\SetFigFont{12}{14.4}{rm}{\color[rgb]{0,.69,0}$\gbf{0\,\,1}$}%
}}}
\put(327,-260){\makebox(0,0)[lb]{\smash{\SetFigFont{12}{14.4}{rm}{\color[rgb]{0,0,0}$\upbfdots\,\bbf{1}\,\bbf{1}\,\zsqov{\widetilde{a_i}}\,\gbf{1}\,\zsqov{a_i''-1}\,\bbf{1}\,\bbf{1}\,\bbf{1}\,\upbfdots$}%
}}}
\end{picture}

%% file: 0-i1-1-i2.pstex_t
\begin{picture}(0,0)%
\includegraphics{0-i1-1-i2.pstex}%
\end{picture}%
\setlength{\unitlength}{4144sp}%
\begingroup\makeatletter\ifx\SetFigFont\undefined
\def\x#1#2#3#4#5#6#7\relax{\def\x{#1#2#3#4#5#6}}%
\expandafter\x\fmtname xxxxxx\relax \def\y{splain}%
\ifx\x\y   
\gdef\SetFigFont#1#2#3{%
  \ifnum #1<17\tiny\else \ifnum #1<20\small\else
  \ifnum #1<24\normalsize\else \ifnum #1<29\large\else
  \ifnum #1<34\Large\else \ifnum #1<41\LARGE\else
     \huge\fi\fi\fi\fi\fi\fi
  \csname #3\endcsname}%
\else
\gdef\SetFigFont#1#2#3{\begingroup
  \count@#1\relax \ifnum 25<\count@\count@25\fi
  \def\x{\endgroup\@setsize\SetFigFont{#2pt}}%
  \expandafter\x
    \csname \romannumeral\the\count@ pt\expandafter\endcsname
    \csname @\romannumeral\the\count@ pt\endcsname
  \csname #3\endcsname}%
\fi
\fi\endgroup
\begin{picture}(3810,1982)(-313,-2188)
\put(-313,-362){\makebox(0,0)[lb]{\smash{\SetFigFont{12}{14.4}{rm}{\color[rgb]{0,0,0}$\upbfdots\,\bbf{1}\,\bbf{1}\,\zsqov{a_{i_1}-1}\,\bbf{1}\,\bbf{1}\,\bbf{1}\,\upbfdots\,\bbf{1}\,\bbf{1}\,\zsqov{a_{i_2}'-1}\,\gbf{1}\,\zsqov{a_{i_2}''-1}\bbf{1}\,\bbf{1}\,\bbf{1}\,\upbfdots$}%
}}}
\put(3235,-2130){\makebox(0,0)[lb]{\smash{\SetFigFont{12}{14.4}{rm}{\color[rgb]{0,0,0}$\gbf{0\,\,1}$}%
}}}
\end{picture}

%% file: 0-j-1-i.pstex_t
\begin{picture}(0,0)%
\includegraphics{0-j-1-i.pstex}%
\end{picture}%
\setlength{\unitlength}{4144sp}%
\begingroup\makeatletter\ifx\SetFigFont\undefined
\def\x#1#2#3#4#5#6#7\relax{\def\x{#1#2#3#4#5#6}}%
\expandafter\x\fmtname xxxxxx\relax \def\y{splain}%
\ifx\x\y   
\gdef\SetFigFont#1#2#3{%
  \ifnum #1<17\tiny\else \ifnum #1<20\small\else
  \ifnum #1<24\normalsize\else \ifnum #1<29\large\else
  \ifnum #1<34\Large\else \ifnum #1<41\LARGE\else
     \huge\fi\fi\fi\fi\fi\fi
  \csname #3\endcsname}%
\else
\gdef\SetFigFont#1#2#3{\begingroup
  \count@#1\relax \ifnum 25<\count@\count@25\fi
  \def\x{\endgroup\@setsize\SetFigFont{#2pt}}%
  \expandafter\x
    \csname \romannumeral\the\count@ pt\expandafter\endcsname
    \csname @\romannumeral\the\count@ pt\endcsname
  \csname #3\endcsname}%
\fi
\fi\endgroup
\begin{picture}(3565,1995)(-298,-2195)
\put(-298,-356){\makebox(0,0)[lb]{\smash{\SetFigFont{12}{14.4}{rm}{\color[rgb]{0,0,0}$\upbfdots\,\bbf{0}\,\bbf{0}\,\unsqov{b_j}\,\bbf{0}\,\bbf{0}\,\upbfdots\,\bbf{1}\,\bbf{1}\,\zsqov{a_i-1}\,\bbf{1}\,\bbf{1}\,\bbf{1}\,\upbfdots$}%
}}}
\put(1780,-2137){\makebox(0,0)[lb]{\smash{\SetFigFont{12}{14.4}{rm}{\color[rgb]{0,0,0}$\gbf{0\,\,1}$}%
}}}
\end{picture}